\documentclass[
  preprint,
  longbibliography,
showkeys,
nofootinbib,
nobibnotes,
amsmath,
amssymb,
]{revtex4-1}

\usepackage{graphicx}
\usepackage{dcolumn}
\usepackage{bm}
\usepackage{todonotes}
\usepackage{subfig}

\usepackage{amsmath}
\usepackage{hyperref}
\usepackage{cleveref}
\usepackage{amsfonts}
\usepackage{amssymb}
\usepackage{algorithm}
\usepackage{algpseudocode}
\usepackage{float}
\usepackage{relsize}
\usepackage{epstopdf}
\usepackage{color}
\usepackage[mathscr]{eucal}
\usepackage{mathtools}

\usepackage[utf8]{inputenc}
\usepackage[english]{babel}
\usepackage{amsthm}

\usepackage{xcolor}
\definecolor{OliveGreen}{rgb}{0,0.4,0}

\newcommand{\nc}{\newcommand}
\nc\be{\begin{eqnarray}}
\nc\ee{\end{eqnarray}}
\nc\bs{\begin{eqnarray*}}
\nc\es{\end{eqnarray*}}
\nc\qu{\quad}
\nc\tends{\rightarrow}
\nc\dint{\int\!\int}

\newcommand{\bes} {\begin{eqnarray*}}
\newcommand{\ees} {\end{eqnarray*}}
\newcommand{\dd} \partial

\newcommand{\non} \nonumber

\nc\cA{{\cal{A}}}
\nc\cB{{\cal{B}}}
\nc\cZ{{\cal{Z}}}
\nc\cX{{\cal{X}}}
\nc\cS{{\cal{S}}}
\nc\cR{{\cal{R}}}
\nc\cW{{\cal W}}
\nc\cU{{\cal U}}
\nc\cP{{\cal P}}
\nc\cF{{\cal{F}}}
\nc\cG{{\cal{G}}}
\nc\bz{\bar{z}}
\nc\bZ{\bar{Z}}
\nc\baz{\bar{\zeta}}
\nc\bw{\bar{w}}
\nc\bW{\bar{W}}
\nc\bX{\bar{X}}
\nc\bet{\bar{\eta}}
\nc\bp{\bar{\phi}}
\nc\bP{\bar{\Phi}}
\nc\bc{\bar{\chi}}
\nc\baf{\bar{f}}
\nc\baF{\bar{F}}
\nc\ov{\overline}
\nc\oW{\ov{W}}

\nc\hA{\hat{A}}
\nc\hB{\hat{B}}
\nc\ho{\hat{O}}
\nc\hG{\hat{\Gamma}}
\nc\hS{\hat{S}}
\nc\hO{\hat{\Omega}}
\nc\bO{\ov{\Omega}}
\nc\gham{\hat{\gamma}}
\nc\gcam{\check{\gamma}}
\nc\z{\zeta}
\nc\s{\sigma}
\nc\ep{\epsilon}
\nc\up{\upsilon}
\nc\lam{\lambda}
\nc\Lam{\Lambda}
\nc\sig{\sigma}
\nc\om{\omega}
\nc\kap{\kappa}
\nc\gam{\gamma}
\nc\Om{\Omega}
\nc\pa{\partial}
\nc\dom{\pa\Omega}
\nc\domt{\pa\tilde{\Omega}}
\nc\doh{\pa\hat{\Omega}}
\nc\pad[2]{\frac{\pa #1}{\pa #2}}
\nc\padd[2]{\frac{\pa^2 #1}{\pa {#2}^2}}
\nc\pard[3]{\frac{\pa^2 #1}{\pa {#2}\pa {#3}}}
\nc\nd[2]{\frac{d #1}{d #2}}
\nc\ndd[2]{\frac{d^2 #1}{d {#2}^2}}
\nc\ds{\displaystyle}
\nc\del{\nabla}
\nc\lap{\nabla^2}
\nc\ud{|\z|\leq 1}

\nc\capil{\mbox{Ca}}

\nc\Rt{\tilde{R}}
\nc\pt{\tilde{a}}
\nc\ft{\tilde{f}}
\nc\vt{\tilde{v}}
\nc\wt{\tilde{w}}
\nc\nt{\tilde{\nu}}
\nc\xit{\tilde{\xi}}
\nc\xt{\tilde{x}}
\nc\etat{\tilde{\eta}}
\nc\nht{\tilde{\hat{\eta}}}
\nc\taut{\tilde{\tau}}
\nc\Taut{T}
\nc\zt{\tilde{\z}}
\nc\nut{\tilde{\nu}}

\nc\At{\tilde{\cA}}
\nc\qt{\tilde{B}}
\nc\Ct{\tilde{C}}
\nc\dt{\tilde{D}}
\nc\Dt{\tilde{\cU}}
\nc\Et{\tilde{L}}
\nc\Ft{\tilde{F}}
\nc\Ht{\tilde{H}}
\nc\Kt{\tilde{K}}
\nc\Lt{\tilde{K}}
\nc\Mt{\tilde{M}}
\nc\Nt{\tilde{N}}
\nc\Pt{\tilde{\Phi}}
\nc\Qt{\tilde{Q}}
\nc\St{\tilde{\cS}}
\nc\Tt{\tilde{\tau}}
\nc\Wt{\tilde{W}}
\nc\cWt{\tilde{\cW}}
\nc\Xt{\tilde{X}}

\nc\vs{\varsigma}
\nc\vp{\varpi}
\nc\ve{\varepsilon}
\nc\vepa{\ve_{\parallel}}
\nc\vepe{\ve_{\perp}}

\newlength{\fchartw}  
\newlength{\fdecw}  
\newlength{\gfchartw}  
\newlength{\gfdecw}  
\usepackage{tikz}
\usetikzlibrary{shapes.geometric, arrows}
\tikzstyle{startstop} = [rectangle, rounded corners, minimum width=3cm, minimum height=1cm, text centered, text width=\fchartw, draw=black, fill=red!30]
\tikzstyle{io} = [rectangle, minimum width=3cm, minimum height=1cm,text centered, text width=\fchartw, draw=black, fill=orange!30]
\tikzstyle{process} = [rectangle, minimum width=3cm, minimum height=1cm,text centered, text width=\fchartw, draw=black,fill=orange!30]
\tikzstyle{decision} = [rectangle, minimum width=3cm, minimum height=1cm,text centered, text width=\fchartw, draw=black, fill=green!30]
\tikzstyle{arrow} = [thick,->,>=stealth]

\usepackage{graphicx}
\usepackage{natbib}
\usepackage{color}
\usepackage{multirow}
\usepackage{amssymb}
\usepackage{amsmath}
\usepackage{empheq}

\begin{document}

\preprint{APS/123-QED}

\title{Filtration with multiple species of particles}

\author{Yixuan Sun, Lou Kondic, Linda J. Cummings}
\affiliation{Department of Mathematical Sciences and Center for Applied Mathematics and Statistics, New Jersey Institute of Technology, Newark, NJ 07102-1982, USA}

\begin{abstract}
Filtration of feed containing multiple species of particles is a common process in the industrial setting. In this work we propose a model for filtration of a suspension containing an arbitrary number of particle species, each with different affinities for the filter membrane. We formulate a number of optimization problems pertaining to effective separation of desired and undesired particles in the special case of two particle species and we present results showing how properties such as feed composition affect the optimal filter design (internal pore structure). In addition, we propose a novel multi-stage filtration strategy, which provides a significant mass yield improvement for the desired particles, and, surprisingly, higher purity of the product as well.

\end{abstract}

\maketitle

\section{Introduction\label{introduction}}

Membrane filtration is widely used in many technological applications \cite{sman2012,reis2007,yogarathinam2018,daufin2011,emami2018,sylvester2013} and in everyday life, for instance in coffee-making and air conditioning. Fouling of the membrane by particles in the feed is unavoidable in successful filtration and understanding of the fouling mechanism(s), critical for improving filtration performance and preventing filter failure, has therefore been the target of significant research effort (see for example~\cite{ives1970, spielman1977,tang2011,iritani2013,iritani2016}). 
Extensive experimental studies \cite{ho1999,tracey1994,jackson2014,lee2019, lee2017, iwasaki1937,lin2009} have been reported for a range of filtration scenarios, mostly focusing on a feed consisting of a single type of particle~\cite{ho1999,tracey1994,jackson2014, iwasaki1937,lin2009}, though possibly with a distribution of particle sizes ~\cite{lee2019,lee2017}. In reality however, filtration typically involves feed containing multiple species of particles (e.g., in gold extraction from ore \cite{ricci2015,acheampong2014}, vaccine extraction \cite{emami2018}, and other bio-product purification after fermentation~\cite{sman2012}), which interact with the membrane differently \cite{chen2004,debnath2019}. 

For feed containing multiple particle species, the goal may be to remove all suspended particles, but there are many applications in which the purpose of the filtration is to remove some particle species from the feed while recovering other species in the filtrate. For example, when producing vaccine by fermentation, one would want to filter the live virus out and retain the vaccine (detached protein shell of the virus, for example \cite{wickramasinghe2010}) in the filtrate. To our best knowledge, little theoretical study has been devoted to feed containing multiple species of particles. Some experimental results are available \cite{acheampong2014,ricci2015,debnath2019}, though the focus is mostly on the specific underlying application rather than mechanistic understanding of how the presence of different particle types affects the filtration process. 

Thanks to recent advances in the development of fast computational tools, numerical solution of the full Navier-Stokes equations and tracking of individual particles in the feed has become a feasible approach for modeling membrane filtration \cite{kloss2012}.
Several such computational fluid dynamics (CFD) studies, particularly focusing on particle deposition on the membrane, have been conducted \cite{wessling2001,kloss2012,bacchin2014,lohaus2018}. Such models may be very detailed, capable of tracking hundreds of millions of particles of arbitrary type and able to reproduce certain experimental data well. However, the computational demand for application-scale scenarios is extremely high; implementation of the CFD method is highly non-trivial and time consuming, and development of simpler models, which can treat different particle populations in an averaged sense, is desirable. 

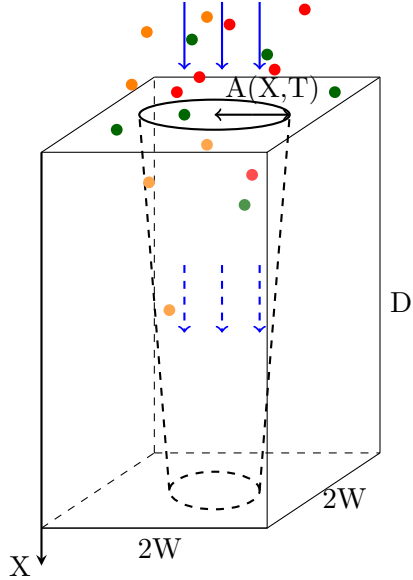
\begin{figure}

\begin{tikzpicture} [xscale=1,yscale=1]

\draw (0,0) rectangle (3,5);
\draw [black, very thin](0,0)--(2,0) node[anchor=north east] {2W} (2,0)--(3,0);
\draw[black] (0,5)--(1.5,6);

\draw [->, thick, blue] (1.9,7) -- (1.9,6.1) ;
\draw [->, thick, blue] (2.4,7) -- (2.4,6.1) ;
\draw [->, thick, blue] (2.9,7) -- (2.9,6.1) ;

\draw [->, dashed, thick, blue] (1.9,3.5) -- (1.9,2.6) ;
\draw [->, dashed, thick, blue] (2.4,3.5) -- (2.4,2.6) ;
\draw [->, dashed, thick, blue] (2.9,3.5) -- (2.9,2.6) ;

\draw [thick,->,>=stealth] (0,5) -- (0,-0.5) node[anchor=east] {X};
\draw[black] (3,5)--(4.5,6);
\draw (1.5,6)--(4.5,6);
\draw[black, thin, dashed] (0,0)--(1.5,1);
\draw[black] (3,0)--(3.6,0.4) node[anchor=west] {2W} (3.6,0.4)--(4.5,1);
\draw (4.5,1)--(4.5,3) node[anchor=west] {D} (4.5,3)--(4.5,6);
\draw [black, thin, dashed](1.5,1)--(4.5,1);
\draw [black, thin, dashed](1.5,1)--(1.5,6);
\draw [black, thick](2.3, 5.5) ellipse (1cm and 0.2cm);
\draw  [->, thick, black]  (2.8, 5.5)--(2.3, 5.5) node[anchor=south west] {A(X,T)};
\draw  [->, thick, black]  (2.8, 5.5) --(3.3, 5.5) ;
\draw [black, thick, dashed](2.3, 0.5) ellipse (0.6cm and 0.25cm);
\draw [black, thick, dashed](1.7,0.5)--(1.3,5.5);
\draw [black, thick, dashed](2.9,0.5)-- (3.3,5.5);

\draw [OliveGreen, ultra thick,fill=OliveGreen] (1,5.3) circle [radius=0.05];
\draw [OliveGreen, ultra thick,fill=OliveGreen] (1.9,5.5) circle [radius=0.05];
\draw [OliveGreen, ultra thick,fill=OliveGreen] (2,6.5) circle [radius=0.05];
\draw [OliveGreen!70!white, ultra thick,fill=OliveGreen!70!white] (2.7,4.3) circle [radius=0.05];
\draw [OliveGreen, ultra thick,fill=OliveGreen] (3.9,5.8) circle [radius=0.05];
\draw [OliveGreen, ultra thick,fill=OliveGreen] (3,6.3) circle [radius=0.05];

\draw [red, ultra thick,fill=red] (1.8,5.8) circle [radius=0.05];
\draw [red, ultra thick,fill=red] (2.1,6) circle [radius=0.05];
\draw [red, ultra thick,fill=red] (2.5,6.7) circle [radius=0.05];
\draw [red!70!white, ultra thick,fill=red!70!white] (2.8,4.7) circle [radius=0.05];
\draw [red, ultra thick,fill=red] (3.1,6.1) circle [radius=0.05];
\draw [red, ultra thick,fill=red] (3.5,6.9) circle [radius=0.05];

\draw [orange, ultra thick,fill=orange] (1.4,6.6) circle [radius=0.05];
\draw [orange, ultra thick,fill=orange] (2.2,6.8) circle [radius=0.05];
\draw [orange, ultra thick,fill=orange] (1.2,5.9) circle [radius=0.05];
\draw [orange!70!white, ultra thick,fill=orange!70!white] (1.43,4.6) circle [radius=0.05];
\draw [orange!70!white, ultra thick,fill=orange!70!white] (2.2,5.1) circle [radius=0.05];
\draw [orange!70!white, ultra thick,fill=orange!70!white] (1.7,2.9) circle [radius=0.05];

\end{tikzpicture}
\caption{\footnotesize{Sketch of a cylindrical pore of radius $A(X,T)$ and length $D$ inside a square prism, representing a basic building-block of the filter membrane (our model is relevant for any other regular tesselating pore-containing prism, {\it e.g.} a hexagonal or triangular prism). Blue arrows indicate the flow direction; colored dots indicate the different particle types present in the feed, and $W$ represents the maximum possible pore radius.}}
\label{prism}
\end{figure}

In earlier work \cite{sanaei2017}, Sanaei and Cummings proposed a simplified model for standard blocking (adsorption of particles, much smaller than the filter pores, onto the internal pore walls), derived from first principles. The model assumes the pore is of slender shape, with pore aspect ratio $\epsilon$ defined as typical width $W$ divided by the length $D$ of the pore, $\epsilon=W/D\ll 1$, see Fig.~\ref{prism}. This provides the basis for an asymptotic analysis of the advection-diffusion equation governing particle transport within the continuum framework, valid for a specific asymptotic range of particle P{\'e}clet numbers (details can be found in \cite{sanaei2017} Appendix A). The model is consistent with one proposed earlier by Iwasaki~\cite{iwasaki1937} based on experiments involving water filtration through sand beds, the validity of which was further confirmed in later experiments by Ison \& Ives~\cite{ison1969}. 

Building on that work, we recently proposed a filtration model focusing on standard blocking with quantitative tracking of particle concentration in the filtrate, which allowed for evaluation of the filtration performance of a given membrane in terms of its pore shape and particle capture characteristics, and for optimization of filtration of a homogeneous feed containing just one type of particles~\cite{sun2020}. In the present work we extend this approach to filtration with multiple species of particles in the feed. For simplicity, we consider dead-end filtration using a track-etched type of membrane. 
We study how the concentration ratio of the different types of particles in the feed, and the differences in membrane--particle interaction characteristics, affect the filtration process and we formulate optimization problems to determine the optimum pore shape (within a given class of shape functions) to achieve the desired objectives. To illustrate our model behavior and its application for design optimization, we explore some hypothetical scenarios of practical interest, in particular: When there are two compounds A and B in a mixture, which filter design will produce the maximum amount of purified compound B before the filter is completely fouled? Questions such as this lead naturally to constrained  optimization problems: how to design a filter such that a certain large fraction of type A particles is guaranteed to be removed, while retaining the maximum yield of type B particles in the filtrate, over the filtration duration? 

We propose new fast optimization methods to solve these problems, based on quantities evaluated at the beginning of the filtration, which are over 10 times faster than the method used in our earlier work \cite{sun2020}. Motivated by some of our findings, we also propose a new multi-stage filtration protocol, which can significantly increase the mass yield per filter of the desired compound, and simultaneously improve the purity of the final product. 

Many variations on the questions we address could be proposed, and the methods we present are readily adapted to a wide range of scenarios.
For brevity and simplicity, however, in the present work we focus chiefly on variants of the example outlined above to illustrate our methods.   
For the majority of the paper we present results for the case in which the feed contains just two particle  species, noting that (within the limitations of our modeling assumptions) our model is readily extended to any number of particle species (some sample results for feed containing more than two species are included in the Appendix \ref{results_multi_species}). 

The remainder of this paper is organized as follows.
We set up our two-species filtration model in \S \ref{model}, focusing attention on the filtration process within a representative pore of the membrane. 
We then outline a number of hypothetical filtration scenarios with multiple species of particles and formulate the corresponding optimization problems in \S\ref{opt_def_cost}.
Although our optimization criteria as defined rely on simulating filtration over the entire useful lifetime of the filter, we will demonstrate the feasibility of using data from the very early stages of our simulations as a reliable predictor of later behavior, offering a much faster route to optimization, discussed in \S \ref{opt_method}. Sample optimization results will be presented in \S \ref{results}.
Section \S \ref{sec:conclusion} is devoted to the summary and discussion.

\section{Filtration modeling with two particle species\label{model}}


In this paper, we focus on dead-end filtration feed solution, carrying multiple different particle species, through a membrane filter.
We first highlight some key modeling assumptions: we assume the particles are non-interacting (justifiable if the feed solution is sufficiently dilute); that the particles are much smaller than the pore radius; and that the pore is of slender shape, with length much larger than its width (this is the case for ``track-etched'' type membranes whose pores are straight and form a direct connection between upstream and downstream sides of the membrane; see, e.g., Apel~\cite{apel2001}). We consider only one type of fouling: the so-called standard blocking mechanism, in which particles (much smaller than pores) are adsorbed on the pore wall leading to pore shrinkage; and we inherit all the additional assumptions made in deriving the standard blocking model proposed by Sanaei and Cummings \cite{sanaei2017}.  Under these assumptions, we set up our model for constant pressure and constant flux conditions, in \S \ref{model_two_species_cp} and \S \ref{model_two_species_cflux} respectively.

\subsection{Solute at constant driving pressure \label{model_two_species_cp}}

We consider 
a feed solution containing two types (different physicochemical properties) of particles, type 1 and type 2, through a planar membrane filter under constant pressure. In the presentation that follows, we use uppercase fonts to denote dimensional quantities and lower case for nondimensional quantities, which will be defined in \S\ref{Nondimensionalization} when we introduce appropriate physical scalings.  We assume that the membrane is composed of identical pores of circular cross-section with radius
$A(X,T)$ (where $X$ is distance along the pore axis), periodically repeating in a regular (e.g., square or hexagonal) lattice arrangement. 
Each circular pore is contained within a regular tesselating polygonal prism, which accommodates a pore of maximum radius $W$ ($0<A\le W$) and height $D$ (see Fig. \ref{prism} for example), where $W\ll D$. We define the representative pore aspect ratio $\ep=W/D \ll 1$, which will be used in our particle deposition model discussed below.  The incompressible feed (assumed Newtonian with viscosity $\mu$) flows through the pore with cross-sectionally averaged axial velocity, $U_{\rm p}(X,T)$, given in terms of the pressure $P(X,T)$ by 
\be
U_{\rm p}(X,T)= -\frac{K_{\rm p} (X,T)}{\mu} \frac{\partial P}{\partial X},  \label{eqn_up}
\ee
where $K_{\rm p}=A^2(X,T)/8$ is the local permeability of an isolated pore (which follows from the Hagen-Poiseuille formula, see e.g., Probstein \cite{probstein1994}, consistent with our pore shape assumption). This is equivalent to a Darcy flow model with velocity $U(T)$ within the membrane related to $U_{\rm p}(X,T)$ via porosity 
$\Phi_{\rm m}=\pi A^2(X,T)/(2W)^2$,
\be
U(T)=\Phi_{\rm m} U_{\rm p}(X,T)=-\frac{K (X,T)}{\mu} \frac{\partial P}{\partial X}, \label{eqn_darcy}
\ee
where 
\be
K(X,T)=\frac{\pi A (X,T)^4}{32 W^2} \label{eqn_membrane_permeability}
\ee
is the membrane permeability.
The flow is driven by constant pressure drop $P_0$ across the membrane. Conservation of mass then closes the model, giving the equation and boundary conditions governing the pressure $P(X,T)$ within the membrane as 
\be
\frac{\partial}{\partial X} \bigl[ K(X,T) \frac{\partial P}{\partial X} \bigr] = 0, \qquad 0\leq X\leq D, \label{eqn_comass} \\
\quad P(0,T)=P_0, \quad P(D,T)=0.  \label{eqn_bc}
\ee


Extending the approach of Sanaei \& Cummings~\cite{sanaei2017}, we propose the following fouling model equations, which assume that the two particle types are transported independently by the solvent and do not interact with each other:
\be
U_{\rm p} \frac{\pa C_i}{\pa X}=-\Lambda_i \frac{C_i}{A}, \quad C_i(0,T)=C_{0i}, \quad i=1,2; \label{eqn_deposition}\\
\frac{\pa A}{\pa T}=-\sum_{i=1,2}\Lambda_i\alpha_iC_i, \quad A(X,0)=A_0(X),\label{eqn_fouling}
\ee
where $C_i(X,T)$ is the concentration (mass per unit volume of solution) of type $i$ particles; $\Lambda_i$ is a particle deposition coefficient for type $i$ particles; and $\alpha_i$ is an unknown (problem-dependent) constant, related inversely to the density of the material that comprises type $i$ particles.  Equations (\ref{eqn_deposition}) follow from a systematic asymptotic analysis (based on the small parameter $\ep$ defined above, see also Table \ref{2t:parameters2}) of advection-diffusion equations for each particle species. Equation (\ref{eqn_fouling}) assumes the rate of pore radius shrinkage (due to the particle deposition) is a linear function of the local particle concentrations at depth $X$, and derives from a mass-balance of the particles removed from the feed, consistent with (\ref{eqn_deposition}). Derivations of these results for filtration of a feed with just one particle type are given in Sanaei \& Cummings~\cite{sanaei2017} Appendix A and~\cite{sanaei2018}.

\subsection{Solute at constant flux \label{model_two_species_cflux}}

Here we briefly consider how the above model is modified for the same feed solution, supplied at constant flux  $U_0$. As fouling occurs the membrane resistance increases, hence the driving pressure must increase to maintain the same flux through the filter. Equation (\ref{eqn_darcy}) still holds for the superficial Darcy velocity, which is now held at constant value $U_0$ by adjusting the driving pressure $P(0,T)$, giving
\be
U= U_0 = -\frac{K (X,T)}{\mu} \frac{\partial P}{\partial X},  \label{eqn_darcy_cflux}
\ee
with just one boundary condition at the membrane outlet,
\be
 \quad P(D,T)=0.
\label{eqn_bc_cflux}
\ee
In this case, the incompressibility condition is satisfied automatically. Equations (\ref{eqn_deposition})-(\ref{eqn_fouling}) then close the model, as in the constant pressure case.

\subsection{Non-dimensionalization}\label{Nondimensionalization}

\subsubsection{Constant pressure}\label{Nondimensionalization_cpressure}

We non-dimensionalize our model (\ref{eqn_darcy})--(\ref{eqn_fouling}) using the following scalings, with lower-case fonts indicating the dimensionless variables: 
\be
p=\frac{P}{P_0}, \qquad u=U\frac{32 D \mu}{\pi W^2 P_0}, \qquad u_{\rm p}=U_{\rm p} \frac{32 D \mu}{\pi W^2 P_0},
\label{eqn_scale_cpressure_u_p}\\
 c_1=\frac{C_1}{C_{01}+C_{02}}, \qquad c_2=\frac{C_2}{C_{01}+C_{02}}, \qquad a=\frac{A}{W}, \label{eqn_scale_cpressure_x_c} \\
x=\frac{X}{D}, \qquad  t=\frac{T}{T_0}, ~~~\text{with}~~~ T_0=\frac{W}{\Lambda_1 \alpha_1 (C_{01}+C_{02})}, \label{eqn_scale_cpressure_a_t}
\ee
where the chosen timescale is based on the deposition rate of particle type 1. The resulting non-dimensionalized equations are listed below:
Eqs. (\ref{eqn_darcy})-(\ref{eqn_bc}) become \be
u= \frac{\pi a^2}{4} u_{\rm p} =-a^4 \frac{\pa p}{\pa x}, \label{darcy_box}  \\ 
\frac{\pa}{\pa x}\left(a^4 \frac{\pa p}{\pa x}\right)=0, \label{incompressibility_box}\\
p(0,t)=1,~~~ p(1, t)=0, \label{bc_box}
\ee
so that dimensionless permeability is just $a^4$ with the chosen scalings;
and Eqs. (\ref{eqn_deposition})-(\ref{eqn_fouling}) take the form
\be
u_{\rm p} \frac{\pa c_1}{\pa x}=-\lambda_1 \frac{c_1}{a}, \quad c_1(0,t)=\xi, \label{eqn_nd_cpressure_c_1}\\
u_{\rm p} \frac{\pa c_2}{\pa x}=-\lambda_2 \frac{c_2}{a}, \quad c_2(0,t)=1-\xi, \label{eqn_nd_cpressure_c_2}\\
\frac{\pa a}{\pa t}=-c_1-\beta c_2, ~a(x,0)=a_0(x), \label{eqn_nd_cpressure_fouling}
\ee
where $\lambda_i={32\Lambda_i D^2 \mu }/({\pi W^3 P_0})$ is the deposition coefficient for particle type $i$, $\xi={C_{01}}/{(C_{01}+C_{02})}$ is the concentration ratio between the two types of particles, $\beta=\Lambda_2 \alpha_2 / (\Lambda_1 \alpha_1)$ is the ratio for effective particle deposition coefficients between the two types of particles, and $0<a_0(x) \leq 1$ is the pore profile at initial time $t=0$. Since we consider scenarios where particle type 1 is to be removed by filtration while type 2 should be retained in the filtrate, only values $\beta\in (0,1)$ will be considered in this paper. The model parameters are summarized in Table \ref{2t:parameters2} for future reference.

\begin{table}
\centering
\begin{tabular}{|l|l|l|}
\hline
{\bf Parameter} & {\bf Description} & {\bf Typical value \& units}\\
\hline
$D$ &\mbox{Membrane thickness}& $300$ \mbox{$\mu$m}\\
$W$ &\mbox{Maximum possible pore radius}& 2~$\mu$m (very variable)\\
$P_0$ & \mbox{Pressure drop} &  unknown ${\rm N}/{\rm m}^2$ (Depends on application)\\
$K$ & \mbox{Representative membrane permeability} & 4$\times$\mbox{10$^{-13}$~{\rm m}$^2$} (very variable)\\
$C_{0i}$ &\mbox{Initial concentration of type $i$ particles in feed } & unknown \mbox{${\rm kg}/{\rm m}^{3}$}\\
$\Lambda_{i}$ &\mbox{Type $i$ particle deposition coefficient } & unknown \mbox{${\rm m}/{\rm s}$}\\
$\alpha_{i}$ &\mbox{Constant related to density of type $i$ particles} & unknown \mbox{${\rm m}^3/{\rm kg}$}\\
$\mu$ &\mbox{Dynamic viscosity} & unknown \mbox{${\rm Pa}\cdot{\rm s}$}\\
\hline
\end{tabular}
\caption{\footnotesize{Dimensional parameters, with approximate values (where known)~\citep{kumar2014}. Depending on the application, pore size may vary from 1 nm to 10 $\mu$m \citep{reis2007}.} 
}\label{2t:parameters1}
\end{table}

\begin{table}
\centering
\begin{tabular}{|c|c|m{7.5cm}|m{4.5cm}|}
\hline
{\bf Parameter} & {\bf Formula} & {\bf Description} & {\bf Value used in simulations}\\
\hline
$\lambda_i $ & ${32\Lambda_i D^2 \mu }/{(\pi W^3 P_0)}$ & Deposition coefficient for type $i$ particles & $\lambda_1 \in \{1, 0.1\}$, $\lambda_2 = \beta \lambda_1$\\
$\xi$ & ${C_{01}}/{(C_{01}+C_{02})}$& Concentration ratio of type 1 particles in feed & $\xi \in \{0.9, 0.5, 0.1\}$\\
$\beta$& $\Lambda_2 \alpha_2 / (\Lambda_1 \alpha_1)$ & Effective deposition coefficient ratio & $\beta \in \{0.1, 0.5, 0.7, 0.9\}$ \\
$\ep$& $W/D$ & Typical pore aspect ratio & Asymptotically small  \\
\hline
\end{tabular}
\caption{\footnotesize{Dimensionless parameters and descriptions (from Table~\ref{2t:parameters1}). 
}}
\label{2t:parameters2}
\end{table}

To solve this system numerically, we first note that 
Eqs. (\ref{darcy_box})--(\ref{bc_box}) 
can be solved to give 
\be
u(t)= \Big ( \int_0^1 \frac{1}{a^4(x,t)} dx \Big )^{-1}. \label{u_solved} 
\ee
Given $a_0(x)$, we compute $u (0)$ via Eq.~(\ref{u_solved}), which allows us to find $u_{\rm p} (x,0)$ via Eq. (\ref{darcy_box}). We then compute $c_1 (x,0), c_2 (x,0)$ via 
Eqs. (\ref{eqn_nd_cpressure_c_1}) and (\ref{eqn_nd_cpressure_c_2}) respectively. With $c_1, c_2$ determined we then compute the pore shape $a(x,t)$ for the next time step via Eq. (\ref{eqn_nd_cpressure_fouling}), then repeat the above  process until the chosen termination condition (based on flux falling below some minimum threshold) for the simulation is satisfied.

\subsubsection{Constant flux}\label{Nondimensionalization_cflux}

Most scales follow from the constant pressure case of \S\ref{Nondimensionalization_cpressure}; here we highlight only the differences for the constant flux scenario, again with lower case fonts indicating the non-dimensionalized variables:
\be
u=\frac{U}{U_0}, \qquad u_{\rm p}=\frac{U_{\rm p}}{U_0}, \qquad p={P} \frac{\pi W^2 }{32 U_0 D \mu}. 
\label{eqn_scale_cflux_u_p}
\ee
The remaining scalings are as in Eqs. (\ref{eqn_scale_cpressure_x_c}) and (\ref{eqn_scale_cpressure_a_t}), leading to the model
\be
u=1=-a^4 \frac{\pa p}{\pa x}, \quad p(1, t)=0, \label{eqn_nd_cflux_u}\\ 
u_{\rm p} = {4\over {\pi a^2}}, 
\label{up_cflux}\\
u_{\rm p} \frac{\pa c_1}{\pa x}=-\lambda_1 \frac{c_1}{a}, \quad c_1(0,t)=\xi, \label{eqn_nd_cflux_c_1} \\
u_{\rm p} \frac{\pa c_2}{\pa x}=-\lambda_2 \frac{c_2}{a}, \quad c_2(0,t)=1-\xi, \label{eqn_nd_cflux_c_2}\\
\frac{\pa a}{\pa t}=-c_1-\beta c_2, \quad a(x,0)=a_0(x). \label{eqn_nd_cflux_a}
\ee
To solve these equations numerically, we proceed as in the constant pressure case, with the simplification that $u=1$ and $u_p$ is a known function of $a$ (\ref{up_cflux}). Note that the inlet pressure $p(0,t)$ is given by
\be
p(0,t)=\int_0^1 \frac{dx}{a^4(x,t)}, \label{eqn_nd_cflux_p}
\ee 
from which it follows that, as the pore radius $a(x,t)$ decreases due to fouling, the driving pressure must increase to maintain the constant flux. We continue the simulation until the specified termination condition (based here on exhausting some fixed amount of feed, subject to a constraint on maximum inlet pressure $p(0,t)$) is reached.

\section{Optimization\label{sec:opt}}

In this section, we explore the specific scenarios introduced in \S \ref{introduction} to find the optimized initial pore shape $a_0(x)$ by defining a suitable objective function $J(a_0)$ with corresponding constraints. For the purpose of the mathematical formulation of the optimization problem, we assume $a_0(x) \in C([0,1])$ (the class of real-valued functions continuous on the real interval $[0,1]$); however, for practical purposes to obtain solutions within reasonable computing time we restrict the search space for the optimizer $a_0(x)$ to low degree polynomial functions (numerical implementation details will be given in \S \ref{opt_method}). In addition, we require $0<a_0(x)\le1$ so that the initial profile is contained within its unit prism (see Fig. \ref{prism}).  

In \ref{opt_def_cost}, we define key metrics that we use to measure the performance of the filter design and use these to set up the optimization problems. Three filtration scenarios will be studied: two under constant pressure conditions and the third under constant flux. In \ref{opt_method}, we outline our optimization methods: first a ``slow method'' (described in \ref{opt_method_slow}) based directly on the objective function defined in \ref{opt_def_cost} below; then we propose a ``fast method'' (in \ref{opt_method_fast}), motivated by results obtained using the slow method. We demonstrate the feasibility of using our model with fast optimization to predict and optimize for various filtration scenarios with multiple species of particles in the feed. 

\subsection{Definitions and objective functions\label{opt_def_cost}}

Adapting the approach taken in our earlier work \cite{sun2020}, we first define some key (dimensionless) quantities that will be used to measure the performance of the membrane. We define instantaneous flux through the membrane as $u(t)$, and cumulative throughput $j(t)$ 
as the time integral of the flux,
\be
 j(t)=\int_0^t u(\tau) d\tau.
 \label{eqn_throughput}
 \ee
 We denote the instantaneous concentration at the outlet ($x=1$) for each particle type $i$ in the filtrate, $c_i(1,t)$ as $c_{i,\rm ins}(t)$, and 
the {\it accumulative concentrations} of each particle type $i$ in the filtrate, $c_{i, \rm acm}$ as
 \be
c_{i, \rm acm} (t)=\frac{\int_0^t c_{i,\rm ins}(\tau) u(\tau) d\tau}{j(t)}.
\label{eqn_c_acm}
\ee
Let $t_{\rm f}$ denote the final time of the filtration process, when the termination condition is reached. For the constant pressure case, we define this to be when the flux drops below some specified fraction $\vartheta$ of its initial value (throughout our work here $\vartheta =0.1$, based on common industrial practice, see e.g., van Reis \& Zydney \cite{reis2007}); for the constant flux case, we consider $t_{\rm f}$ to be the fixed time at which the specified amount of feed is exhausted, assuming that the terminal driving pressure $p(0,t_{\rm f})$ is less than the maximum operating pressure $p_{\rm max}$ for all initial pore profile functions $a_0(x)$ in the searching space considered. 

{\color{black}To specify the particle removal requirement from the feed for each type of particles, we define the {\it instantaneous particle removal ratio} for type $i$ particles, $R_{i}(t) \in [0,1]$, as
\be
R_{i}(t)=1- \frac{c_{i,\rm ins}(t)}{c_i(0,t)},
\label{R}
\ee
where $c_{i,\rm ins}(t)$ is instantaneous concentration of particle type $i$ at the outlet and $c_i(0,t)$ is the type $i$ particle concentration in the feed at time $t$.\footnote{In the problems that we consider $c_i(0,t)=c_{i0}$ is fixed ($c_1(0,t)=c_{10}=\xi$ and $c_2(0,t)=c_{20}=1-\xi$), but if we wish to consider feed with time-varying particle concentrations, then $c_i(0,t)$ in (\ref{R_acm}) should be replaced by appropriate averaged concentrations, $\bar c_i(0,t):=({\int_0^t c_i(0,\tau) u(\tau) d\tau})/{j(t)}$.} Then the initial particle removal ratio $R_{i}(0)$ is the fraction of type $i$ particles removed after the feed passes through the clean filter. We also define the {\it cumulative particle removal ratio} for type $i$ particles, ${\bar R}_i (t)$, as 
\be
{\bar R}_i (t)=1- \frac{c_{i, \rm acm} (t)}{c_i(0,t)}, \qquad i=1,2, 
\label{R_acm}
\ee
where $c_{i, \rm acm} (t)$ is defined in (\ref{eqn_c_acm}). The final cumulative particle removal ratios at the end of the filtration are then ${\bar R}_i (t_{\rm f})$.
 
\begin{table}
\centering
\begin{tabular}{|l|l|l|}
\hline
{\bf Metric} & {\bf Description} & {\bf Range/value/definition}\\
\hline
$u(t)$ &\mbox{flux}& $\in (0,\infty)$ \\
$j(t)$ &\mbox{throughput}& $=\int_0^t u(\tau) d\tau$\\
$j(t_{\rm f})$ &\mbox{total throughput at final time $t_{\rm f}$}& $=\int_0^{t_{\rm f}} u(\tau) d\tau$\\
$c_{i,\rm ins}(t)$ & \mbox{instantaneous concentration at the outlet for each particle type $i$} &  $=c_i(1,t)\in (0,c_i(0,t))$ \\
$c_{i,\rm acm}(t)$ & \mbox{accumulative concentrations of each particle type $i$ in the filtrate} & $\in (0,c_i(0,t))$ (defined in Eq.(\ref{eqn_c_acm}))\\
$R_{i}(t)$ &\mbox{{instantaneous particle removal ratio} for type $i$ particles} & $\in [0,1]$\\
${\bar R}_i (t)$ &\mbox{{cumulative particle removal ratio} for type $i$ particles} & $\in [0,1]$\\
$\Rt$ &\mbox{{desired final cumulative particle removal ratio} for type 1 particles } & $0.99$\\
$\Upsilon$ &\mbox{{desired fraction of type 2 particles in filtrate (effective separation)}} & $0.5$\\
$\vartheta$ & flux fraction at termination (constant pressure filtration) & 0.1 \\
$k_i $ &\mbox{purity for type $i$ particles in the filtrate at the end of filtration } & $\in [0,1]$\\
$\gamma$ &\mbox{{effective physicochemical difference} between the two species} & $\in [0,1]$\\

\hline
\end{tabular}
\caption{\footnotesize{Key metrics defined in \S \ref{opt_def_cost} and  \S \ref{result_multi-stage} for measuring membrane performance and their ranges, values (where fixed across all simulations) or definitions.}}\label{performance_metrics}
\end{table}
Preliminary investigations for our  multi-species filtration model indicate that the particle removal capability of the filter improves, for the constant pressure scenarios, as the filtration proceeds and pores shrink, thus in our optimizations we impose the particle removal requirement only at the initial step, i.e., we require $R_1(0)$ to be greater than a specified number ($R_{}$)  between 0 and 1. Throughout this work we consider the desired final particle removal ratio for type 1 particles to be 0.99 and denote this fixed value by $\Rt$.  Other values of $R$ are used for ``intermediate'' filtration stages in our description of multi-stage filtration later, with the understanding that the final goal is to reach removal ratio $\Rt=0.99$. With these definitions, for the constant pressure case, we illustrate our methods by considering a number of membrane design optimization scenarios, outlined below. 

\noindent
{\bf Problem 1.} In many situations there are competing demands and it may be useful to consider objective functions that assign weights to different quantities of interest. Suppose we have a feed with known concentrations of type 1 and type 2 particles, where the goal is to remove type 1 particles from the feed, while retaining type 2 particles in the filtrate and simultaneously collecting as much filtrate as possible, until the termination time  $t_{\rm f}:={\rm inf} { \Big \{ t: u(t) \le \vartheta u(0) \Big \} }$ is reached. Which filter design $a_0(x)$ -- the initial pore profile, within our searching space -- will remove a specified fraction $R_{}\in [0,1]$ of type 1 particles and simultaneously maximize the objective function $J(a_0):=w_1 j(t_{\rm f})+ w_2 c_{2 \rm acm}(t_{\rm f})$ (where $w_1$ and $w_2$ are weights associated to the total throughput and final cumulative concentration of type 2 particles in the filtrate, respectively)? 
For example, in water purification \cite{hoslett2018}, type 1 particles could be toxins like lead (which we insist are removed), while type 2 particles are desirable minerals. 
In this application, it is of interest to retain type 2 particles, but the primary concern is to produce the purified water, so a larger value might be assigned to $w_1$ than $w_2$. 
This example motivates the following design optimization problem.

\noindent\fbox{%
\parbox{\textwidth}{%
\textbf{Optimization Problem 1}

\textbf{Maximize} 
\be
J(a_0):=  w_1 j(t_{\rm f})+ w_2 c_{2 \rm acm}(t_{\rm f})  \label{max_Jplusc2acm_slow}
\ee

\noindent \textbf{subject to} Eqs. (\ref{bc_box})-(\ref{eqn_nd_cpressure_fouling}),
\textbf{and}
\bs
0 < a_0 (x) \leq 1,  \quad  \forall x \in [0, 1], \\
R_1(0) \ge R_{} ,\\
t_{\rm f}={\rm inf} { \Big \{ t: u(t) \le \vartheta u(0) \Big \} }.
\es
    }%
}
Here, we seek the optimum pore shape $a_0(x)$ to maximize $J (a_0)$, a weighted combination of $j(t_{\rm f})$ and $c_{2 \rm acm}(t_{\rm f})$, subject to the flow and fouling rules dictated by our model (Eqs. (\ref{bc_box})-(\ref{eqn_nd_cpressure_fouling})), and the physical constraints that the pore is initially contained within the unit prism (so that adjacent pores cannot overlap), and the desired user-specified fraction $R_{}$ of type 1 particles is removed from the feed at the start of filtration. For example, if $w_1=1$, $w_2=0$, $R_{}=\Rt$ then we are maximizing the total throughput of filtrate, with a hard constraint that at least $99\%$ of type 1 particles are removed initially, and no concern for the proportion of type 2 particles retained in the filtrate. On the other hand, if  $w_1=0.5$, $w_2=0.5$, $R_{}=\Rt$ then (assuming the dimensionless quantities $j(t_{\rm f})$  and $c_{2 \rm acm}(t_{\rm f})$ are of similar magnitude) we care equally about total throughput and the proportion of type 2 particles retained in the filtrate, again with a hard constraint on removal of type 1 particles. 

\noindent
{\bf Problem 2.} Suppose we have a large quantity of feed containing known concentrations of type 1 and type 2 particles, where the goal is to remove type 1 particles and collect the maximum quantity of type 2 particles in the filtrate (e.g., for vaccine production after fermentation, one would want to filter out the live virus -- type 1 particles -- and retain as much vaccine -- type 2 particles -- as possible in the filtrate), until the termination time  $t_{\rm f}$ is reached. 
Which filter design $a_0(x)$, within our searching space, will remove a specified fraction $R_{}\in [0,1]$ of type 1 particles and simultaneously maximize the final yield of type 2 particles in the filtrate, $c_{2 \rm acm}(t_{\rm f}) j(t_{\rm f})=:J(a_0)$? This question leads to the following design optimization problem.

\noindent\fbox{%
\parbox{\textwidth}{%
\textbf{Optimization Problem 2}

\textbf{Maximize} 
\be
J(a_0):=  c_{2 \rm acm}(t_{\rm f})  j(t_{\rm f}) \label{max_Jtimesc2acm_slow}
\ee
\noindent \textbf{subject to} Eqs. (\ref{bc_box})-(\ref{eqn_nd_cpressure_fouling}),
\textbf{and}
\bs
0 < a_0 (x) \leq 1,  \quad  \forall x \in [0, 1], \\
R_1(0) \ge R_{} ,\\
t_{\rm f}={\rm inf} { \Big \{ t: u(t) \le \vartheta u(0) \Big \} }.
\es
   }%
}

Here we seek the optimum $a_0(x)$ that maximizes objective function $J(a_0)$, representing the final mass of type 2 particles in the filtrate, subject to the flow and fouling rules of our model, the physical constraints, and the desired particle removal requirement. 

For the case in which constant flux through the filter is specified, we consider the following illustrative scenario: \\
\noindent
{\bf Problem 3.}  
Given a fixed amount of feed, what is the best filter design
to maximize the yield of a purified compound of interest (e.g., gold in mining \cite{ladner2012}, vaccine extraction \cite{emami2018} or other bio-product purification), while removing an impurity?
In this scenario, we consider the optimization problem as finding the initial pore profile $a_0 (x)$ such that the filter removes a certain proportion ($R_{}$) of type 1 particles from the feed, while maximizing the amount of type 2 particles collected in the filtrate, until the feed is exhausted, i.e. the termination time $t_{\rm f}$ is reached. This problem statement motivates the following optimization problem:

\noindent\fbox{%
\parbox{\textwidth}{%
\textbf{Optimization Problem 3}

\textbf{Maximize} 
\be
J(a_0):=  c_{2 \rm acm}(t_{\rm f}) j(t_{\rm f})  \label{max_Jtimes2acm_slow_cflux}
\ee

\noindent \textbf{subject to}
Eqs. (\ref{eqn_nd_cflux_u})-(\ref{eqn_nd_cflux_a}), \textbf{and}
\bs
0 < a_0 (x) \leq 1,  \quad  \forall x \in [0, 1], \\
R_1(0) \ge R_{} ,\\
t_{\rm f} \quad\mbox{specified (time at which feed is exhausted)}.
\es
    }%
}
The objective function $J (a_0)$ in (\ref{max_Jtimes2acm_slow_cflux}) represents choosing the optimum $a_0(x)$ to maximize the final mass of purified type 2 particles obtained at the end of the filtration. 

\subsection{Optimization methodology\label{opt_method}}

The design optimization problems outlined above are mathematically challenging and computationally expensive in general due to non-convexity \cite{boyd2004} (of both the objective function and the constraints), large number of design variables (in our case the number of possible design variables is infinite, as our searching space for the pore shape $a_0(x)$ is the infinite-dimensional function class $C([0,1])$) and the computational cost of evaluating the objective function (which requires that we solve the flow and transport equations until the termination time $t_{\rm f}$) many times. 
For simplicity and efficiency we therefore restrict our searching space for $a_0(x)$ to be (low degree) polynomial functions, the coefficients of which represent our design variables (the class of searchable functions could be expanded without difficulty but with commensurate increase in computational cost). In the following two subsections we outline our two optimization routines: first the slow method, which arises naturally from the problems posed and which relies on running many simulations over the entire lifetime of the filter; then the proposed new fast method, which uses data from only the very earliest stages of filtration to predict the optimum over the filter lifetime. The two methods are compared in \S\ref{results}.  

\subsubsection{Slow method\label{opt_method_slow}}
 
We vary the coefficients of polynomials $a_0(x)$ to find the values that maximize the objective functions defined in \textbf{Problems 1-3}, under the constraints specified in each case. The functions $a_0(x)$ are referred to as shape functions in the shape optimization literature \cite{ta'asan1992}. In the interests of reducing computation time, for the purpose of the demonstration simulations presented here, we restrict our searching space to be the linear pore profile, i.e. we consider initial profiles of the form $a_0(x)=d_1 x+d_0$, where $d_1, d_0$ are the design variables to be optimized for each specific scenario, with searching range $(d_1, d_0) \in [-1,1]\times [0,1]$.  We use the  {\tt MultiStart} method with {\tt fmincon}
 as local solver from the {\tt  MATLAB}\textregistered~{\tt Global Optimization} toolbox for this optimization. Since the routine finds a minimizer while we want to  maximize $J (a_0)$, we work with the cost function $-J (a_0)$; more details of the implementation of the cost function and constraints can be found in our earlier work \cite{sun2020}. 
 
We specify a starting point $(d_1^0, d_0^0) \in [-1,1]\times [0,1]$  (initial guess for running the local solver {\tt fmincon}), 
cost function (based on our objective functions and constraints), 
design variable searching range, 
and number of searching points $n$ (the number of points in $(d_1,d_0)$-space that will be explored) for the {\tt MultiStart} method. 
With the user-specifed starting point $(d_1^0, d_0^0)$, an additional $(n-1)$ starting points $(d_1^i, d_0^i) \in [-1,1]\times [0,1], i=1, 2, ...n-1$ are generated by the {\tt MultiStart} algorithm. 
The resulting $n$ points are then used to run the local solver {\tt fmincon} (based on a gradient descent method) to find a list of local minimizers. 
We use the best minimizer from the list as the coefficients for our optimized linear pore profile. 
Note that there is no guarantee the method will find the global minimizer due to the nature of gradient descent methods applied to non-convex problems (the result found depends on the starting-points); however, local minimizers can be systematically improved
and for practical purposes may be useful if they provide significant improvement over current practice (see e.g., the study by Hicks {\it et al.} on airfoil design \cite{hicks1974}). 
Simulations using this method are presented in Figures \ref{P1_f_s_pore_evl}, \ref{P1_vary_W}, \ref{P2_vary_xi_beta}, \ref{fast_method_vary_search_pt}, described in \S \ref{results} later.
 
\subsubsection{Fast method\label{opt_method_fast}}

The ``slow'' optimization method described above is straightforward and easy to implement, but reliable results require that many ($n$ large) individual model simulations be run through to the termination time $t_{\rm f}$. The results presented in this paper are restricted to optimizing membrane structure within the class of linear pore profiles only, but in any real application it may be desirable to optimize over wider function classes, e.g., polynomials of higher order. We find (empirically) that each unit 
increase of the degree of polynomial $a_0(x)$ requires roughly a 10-fold increase in the number of searching points to reach the best local optimum, with a corresponding increase in the run time. Run time will also increase if more than two particle species are considered, or if some of the constraints are removed or inactive (e.g., a less strict particle removal requirement) and the feasible region becomes very large. Maximum computational efficiency in practical situations is therefore critical. Motivated by the idea that imposing carefully-chosen conditions on the initial state of a system can, in many cases, guarantee certain features of later states, we propose a fast method based on simulations of the very early stages of filtration. We note that similar ideas have been used to estimate filter capacity (the total amount of feed processed during a filtration)  using a method called $V_{\rm max}$, which essentially predicts the filter capacity using only the first 10-15 minutes of filtration data \cite{zydney2002}.

Extensive preliminary simulations for {\bf Problem 1}, with $w_2=0$,\footnote{In which we optimize for total throughput only 
in constant pressure-driven flow for two particle species with $\lambda_1=1$, $\beta \in (0.1,0.9)$ and $\xi \in (0.1,0.9)$, using the slow method outlined above.} indicate that the $(u(t), j(t))$ flux-throughput graph at optimum is initially flat and high (see Fig.~\ref{P1_vary_W}(a) for example), with small gradient $|u'(0)|$ and large vertical intercept $u(0)$ in comparison with graphs for sub-optimal solutions. 
Moreover, fouling shrinks the pore and increases resistance, thus flux decreases in time and $u'(0)<0$ for all model solutions. We therefore expect that, at optimum, $u(0)$ should be as large as possible and $u'(0)$ as close to zero as possible. 
Similar ideas apply to the case $w_2>0$, where we wish also to maximize $c_{2 \rm acm}(t_{\rm f})$, the cumulative concentration of type 2 particles in the feed at the final time: we propose instead to maximize a function based on the initial state of the system as characterized by $c_{2\rm ins}(0)$ and $c_{2\rm ins}'(0)$ (respectively the initial concentration and initial concentration gradient, with respect to $t$, of type 2 particles at the membrane outlet $x=1$).
Again, preliminary simulations indicate that at optimum $c_{2\rm ins}'(0)$ is close to zero and negative,\footnote{The occasional increase in particle concentration at the membrane outlet that was observed at early times in our previous work~\cite{sun2020} for single-particle-species filtration was never seen here.} while $c_{2\rm ins}(0)$ is large, compared to sub-optimal solutions. With these observations, we expect to maximize $c_{2 \rm acm}(t_{\rm f})$ by insisting on large $c_{2\rm ins}(0)$ and small initial gradient $c_{2\rm ins}'(0)$. 

With these motivations, we now define modified objective functions for our fast method. In place of
(\ref{max_Jplusc2acm_slow}) in Optimization \textbf{Problem 1}, we propose the following fast objective function, which uses data from only the initial stage of the model solution:
\be
J_{1, \rm fast}(a_0)=w_1u(0)+w_1u'(0)+ w_2c_{2\rm ins}(0)+w_2c_{2\rm ins}'(0), \label{max_Jplusc2acm_fast}
\ee
in which the terms in $w_1$ act to maximize total throughput and those in $w_2$ maximize concentration of type 2 particles in the filtrate, where $w_1$ and $w_2$ can be tuned depending on the relative importance of the two quantities. Note that the weights assigned to $u(0)$ ($c_{2\rm ins}(0)$) and $u'(0)$ ($c_{2\rm ins}'(0)$) do not have to be the same; we could allow four independent weights for the four quantities in (\ref{max_Jplusc2acm_fast}).  However, for the simple application scenarios we considered we found just two independent weights $w_1$, $w_2$ to be sufficient to give reliable results in an efficient manner.

To replace (\ref{max_Jtimesc2acm_slow}) in Optimization \textbf{Problem 2}, we propose the following fast objective function
\be
J_{2, \rm fast}(a_0)=u(0)c_{2\rm ins}(0), \label{max_Jtimes2acm_fast}
\ee
in which $u(0)c_{2\rm ins}(0)$ captures the initial collection of the particle 2 in the filtrate.
Other forms involving $u'(0)$ and $c_{2\rm ins}'(0)$ were tested, 
but found to confer no improvements, hence we opt for the simplest effective objective.

In the following section we demonstrate that our fast optimization method always gives results at least as good as those for the slow method, and then use it to investigate various model features and predictions.

\section{Results\label{results}}
In this section we present our simulation results for several two-species filtration scenarios. 
We focus on the effects of $\xi$, the concentration ratio of the two particle types in the feed, and $\beta =\Lambda_2\alpha_2/(\Lambda_1\alpha_1)$, the ratio of the effective particle deposition coefficients for the two particle types (both these parameters are unique to multi-species filtration, having no counterparts in single-species models). 
For most of our simulations, we fix $\lambda_1=1$ (particle type 1 has fixed affinity for the membrane throughout) and the initial fraction of type 1 particles to be removed is fixed at $\Rt$. In \S \ref{result_two_species_cpressure}, we first present sample comparison results between the fast and slow methods for \textbf{Problem 1} and \textbf{Problem 2}, noting that many more tests than are presented here were conducted to verify that the fast method reliably finds optima as good as or superior to those found by the slow method, under a wide range of conditions. 
We then use the fast method to study the effects of varying parameters $\beta$ and $\xi$ for two-species filtration under constant pressure conditions.
Based on these observations, we propose a multi-stage filtration strategy that will increase the mass yield of particle we wish to recover in \S \ref{result_multi-stage}.
We present sample results for the constant flux case in \S \ref{result_two_species_cflux}, focusing on \textbf{Problem 3}.

\subsection{Optimization of constant pressure filtration}\label{result_two_species_cpressure}

\subsubsection{Efficacy of the fast optimization method}

\begin{figure}
{\scriptsize (a)}\rotatebox{0}{\includegraphics[scale=.36]{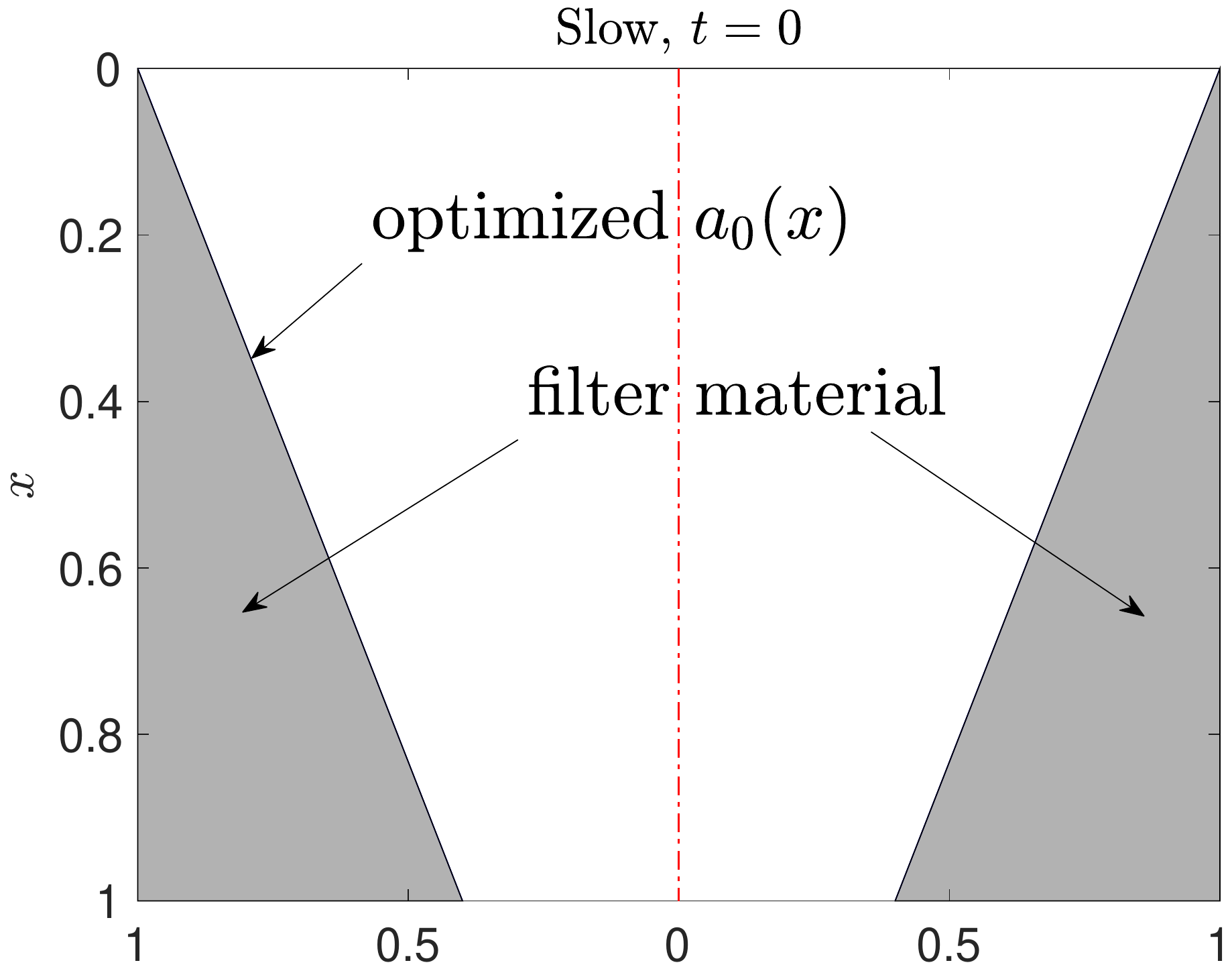}}
{\scriptsize (d)}\includegraphics[scale=.36]{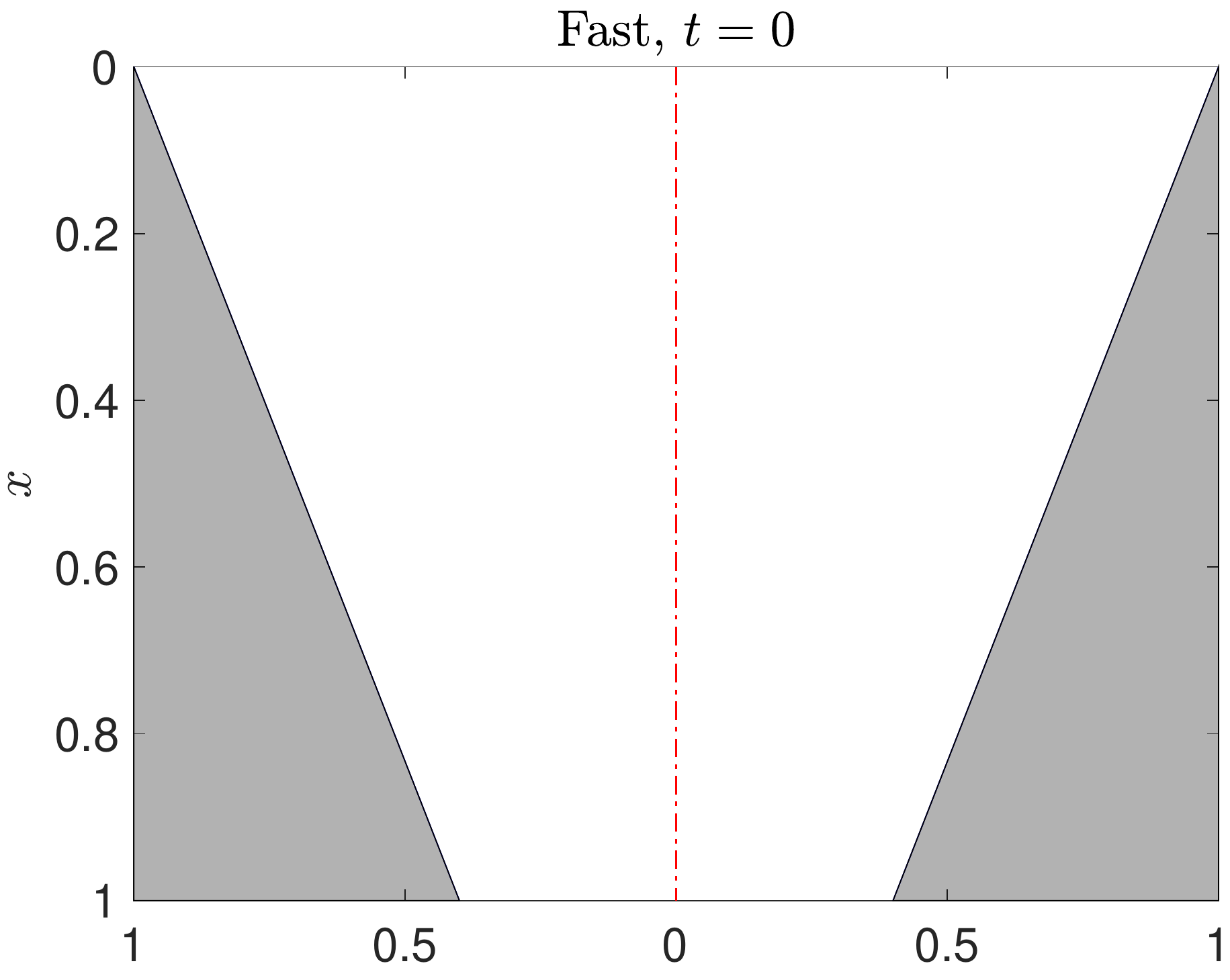}
{\scriptsize (b)}\includegraphics[scale=.36]{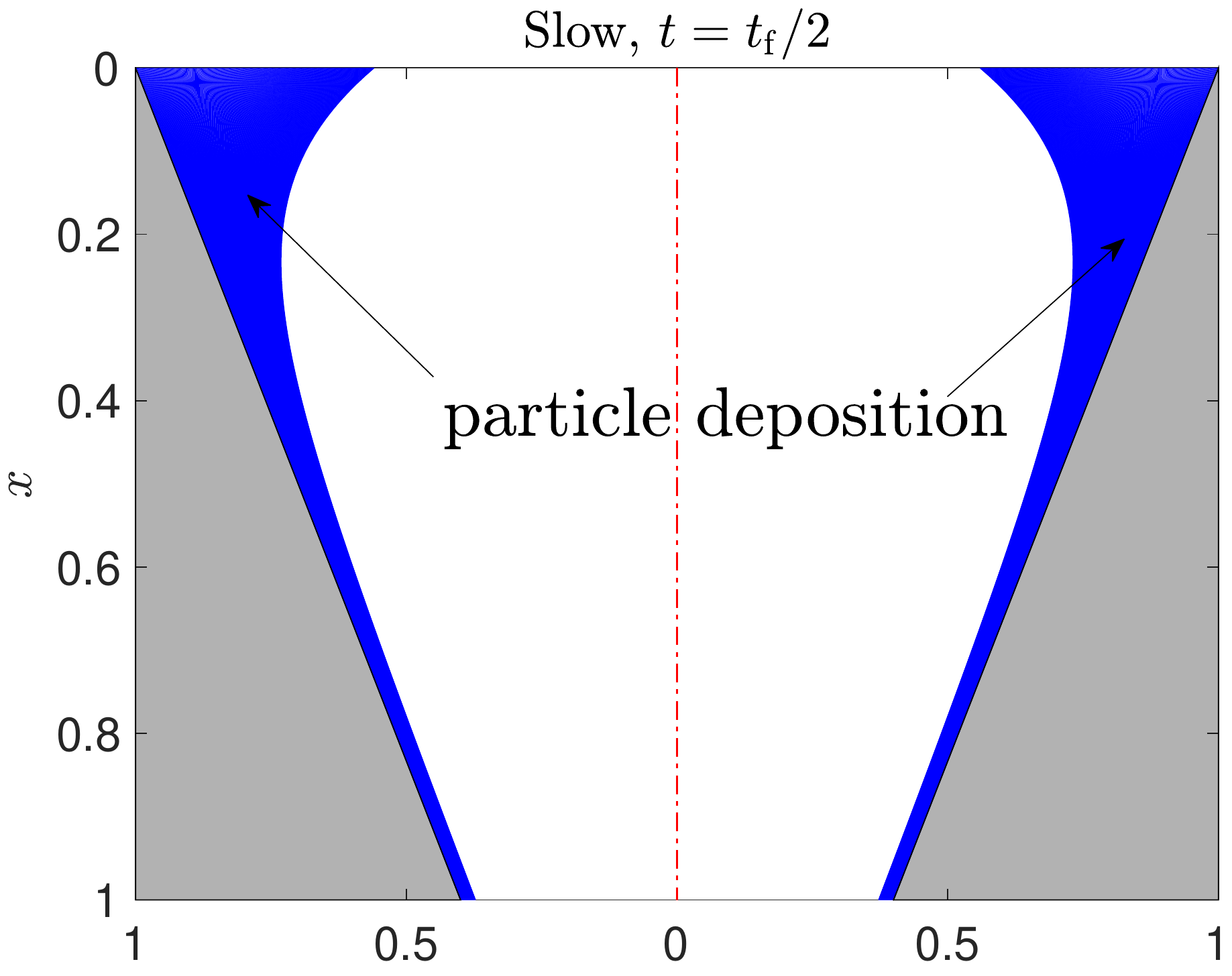}
{\scriptsize (e)}\includegraphics[scale=.36]{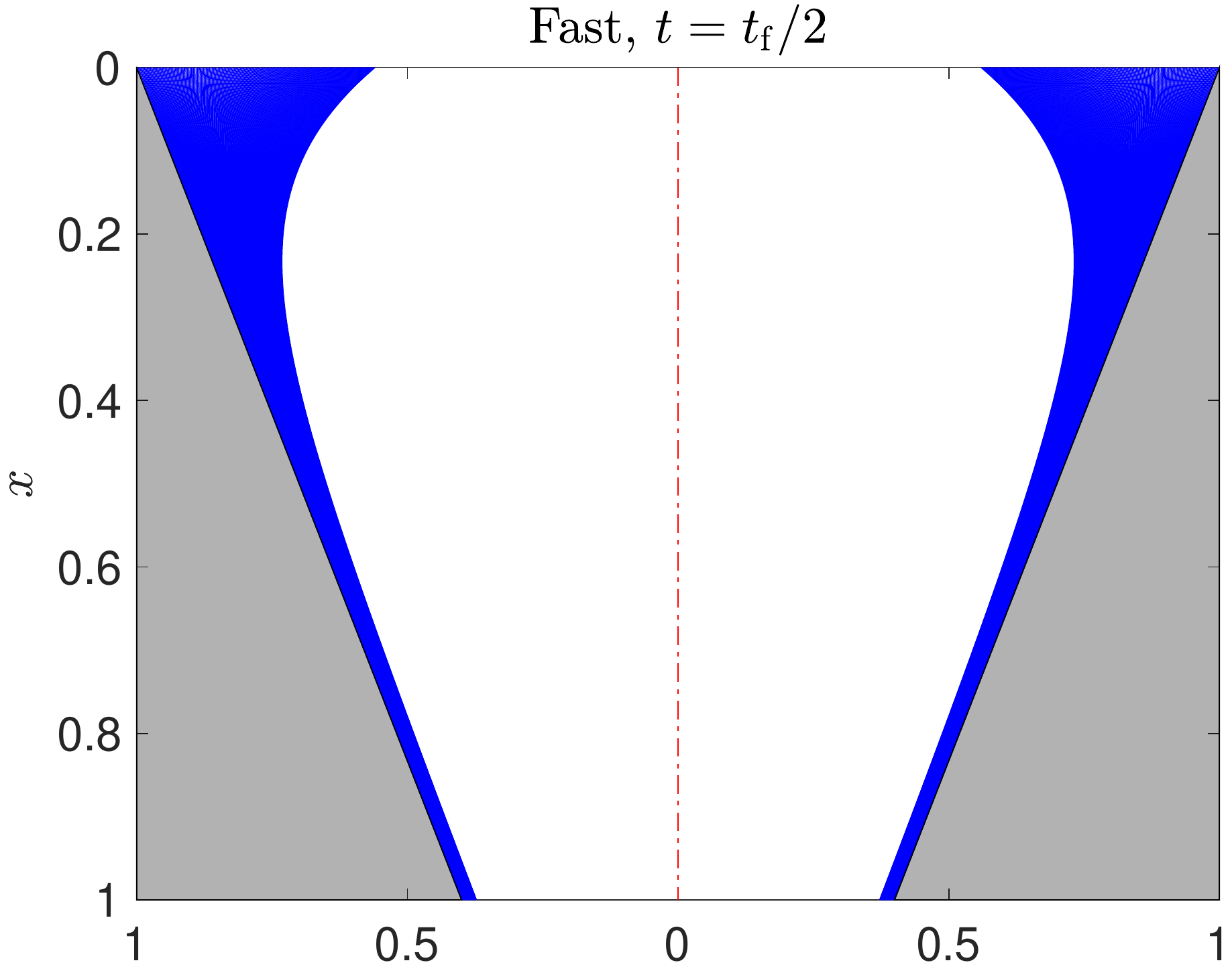}
{\scriptsize (c)}\rotatebox{0}{\includegraphics[scale=.36]{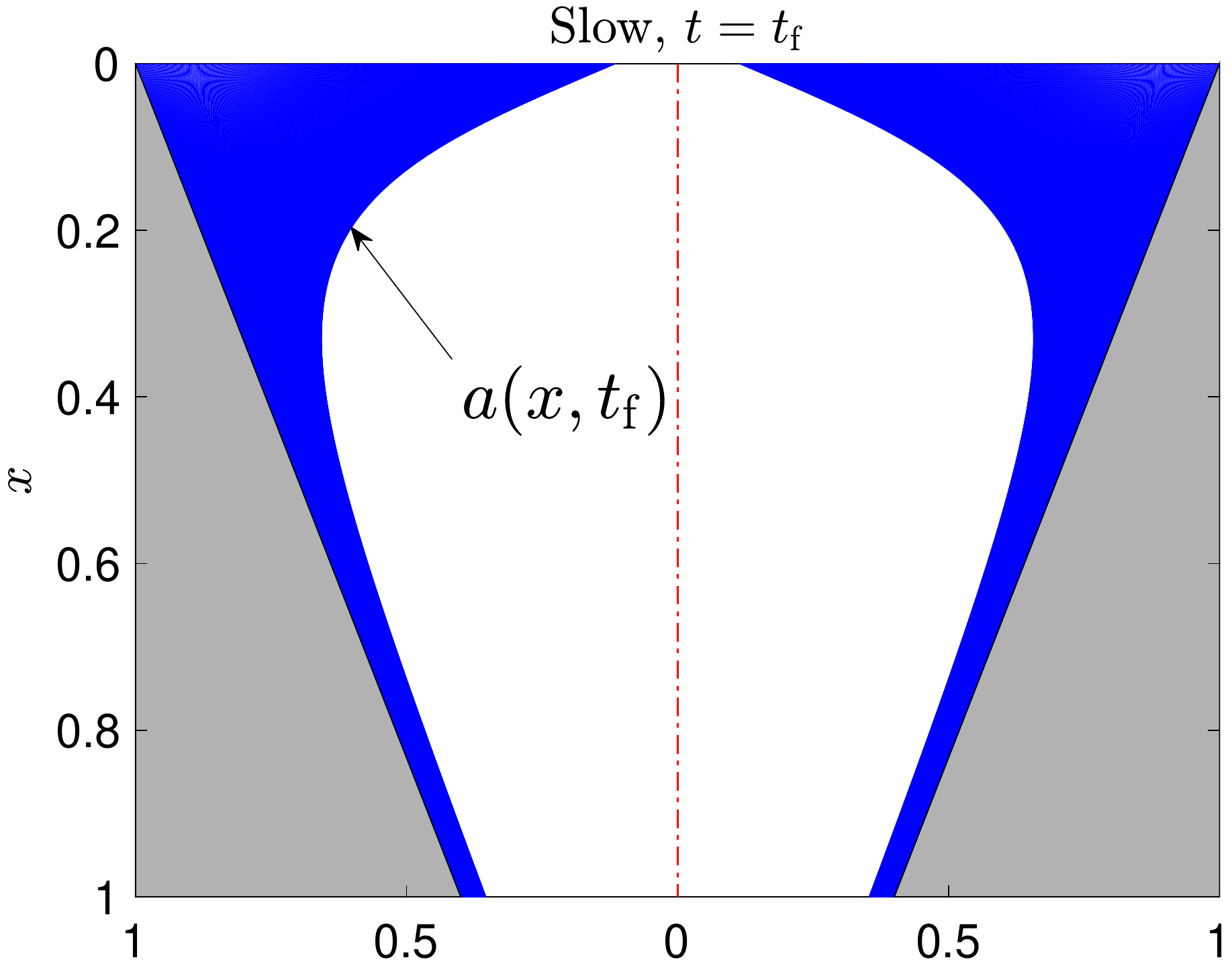}}
{\scriptsize (f)}\includegraphics[scale=.36]{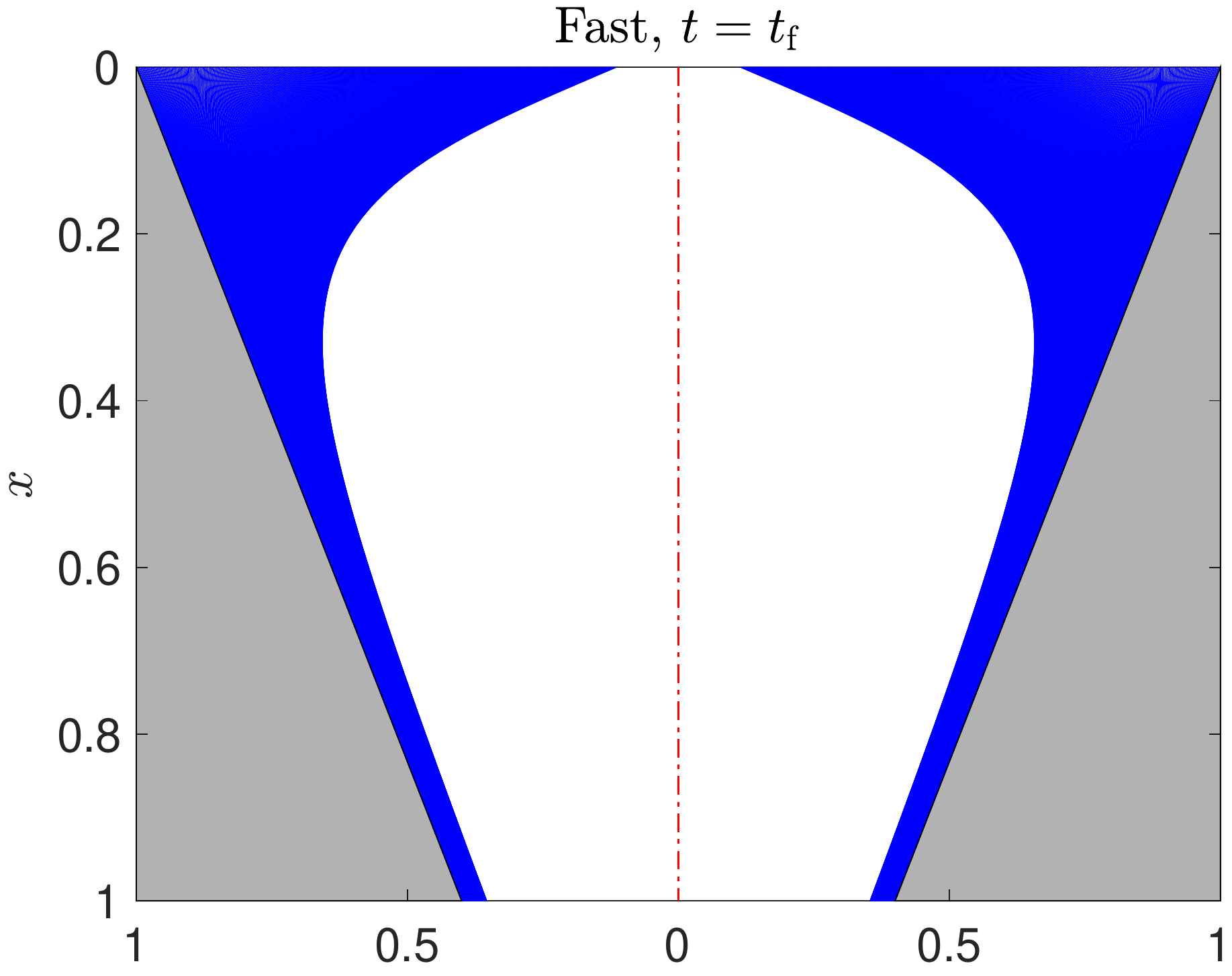}
\caption{
\footnotesize{ 
Fouling  evolution of the optimized membrane pore $a_0(x)$ for {\bf Problem 1} with $ w_1=1, w_2=0, \xi=0.5, \beta=0.1, \alpha_1=\alpha_2, \lambda_1=1$; 
(a-c) show evolution for
$a_0(x)=-0.6001x+0.9998$
optimized using slow method and (d-f) show evolution of 
$a_0(x)=-0.6002x+0.9999$ 
for the corresponding fast method, at $t=0$ (unfouled; (a) and (d)), $t=t_{\rm f}/2$ (halfway through filtration; (b) and (e)) and $t=t_{\rm f}$ (end of filtration; (c) and (f)).}}
\label{P1_f_s_pore_evl}
\end{figure}

\begin{figure}
{\scriptsize (a)}\includegraphics[scale=.39]{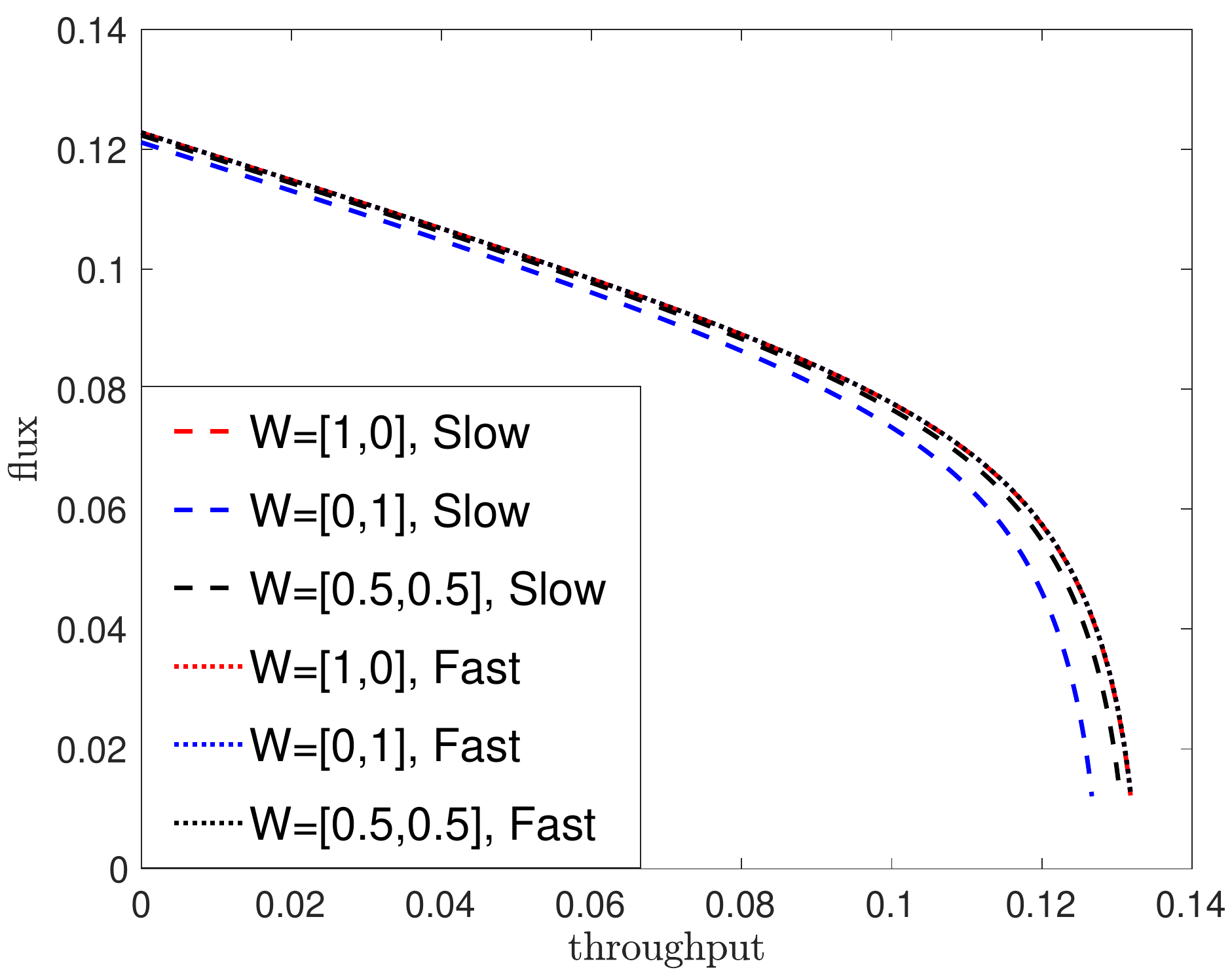}

{\scriptsize (b)}\includegraphics[scale=.39]{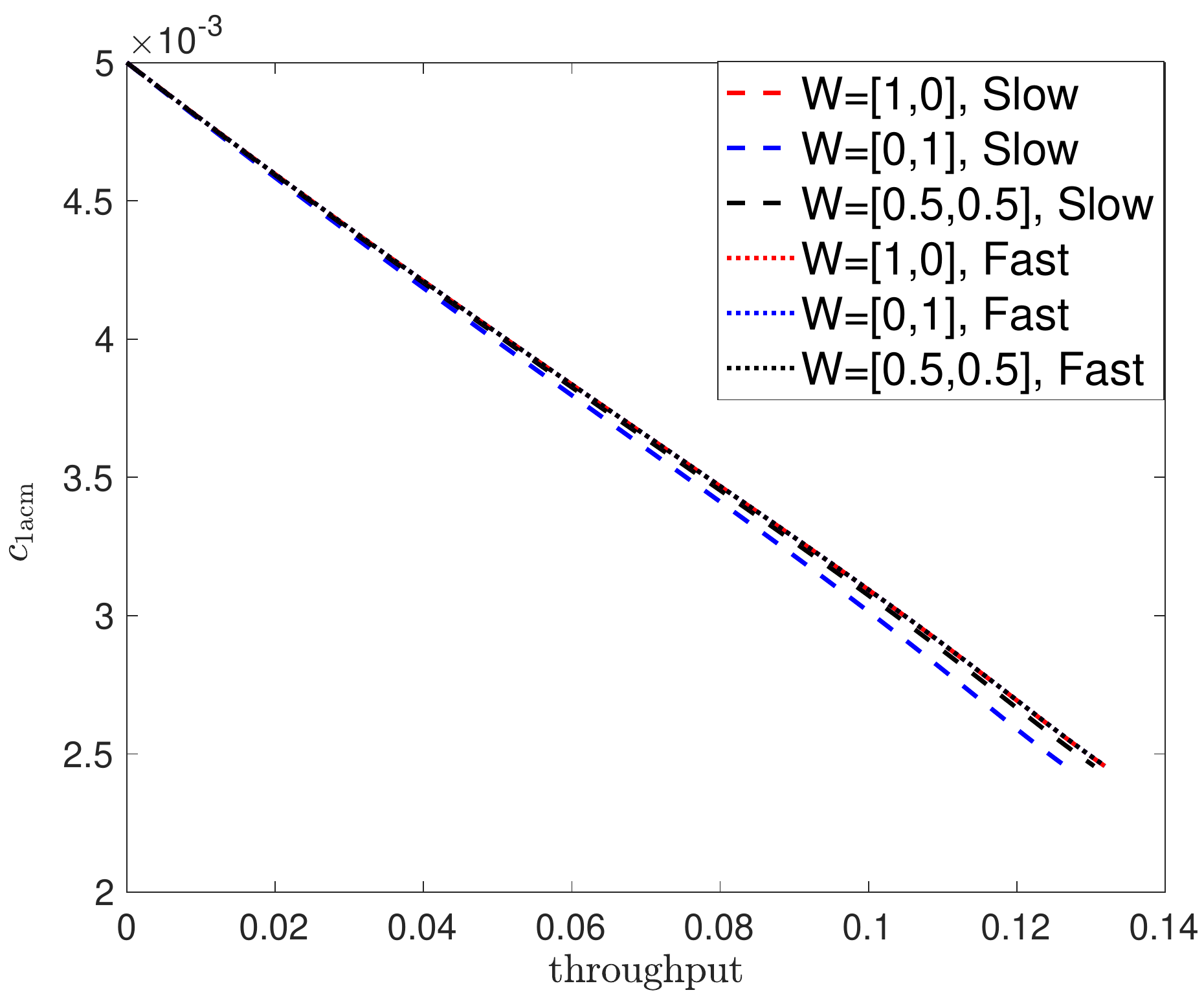}
{\scriptsize (c)}\includegraphics[scale=.39]{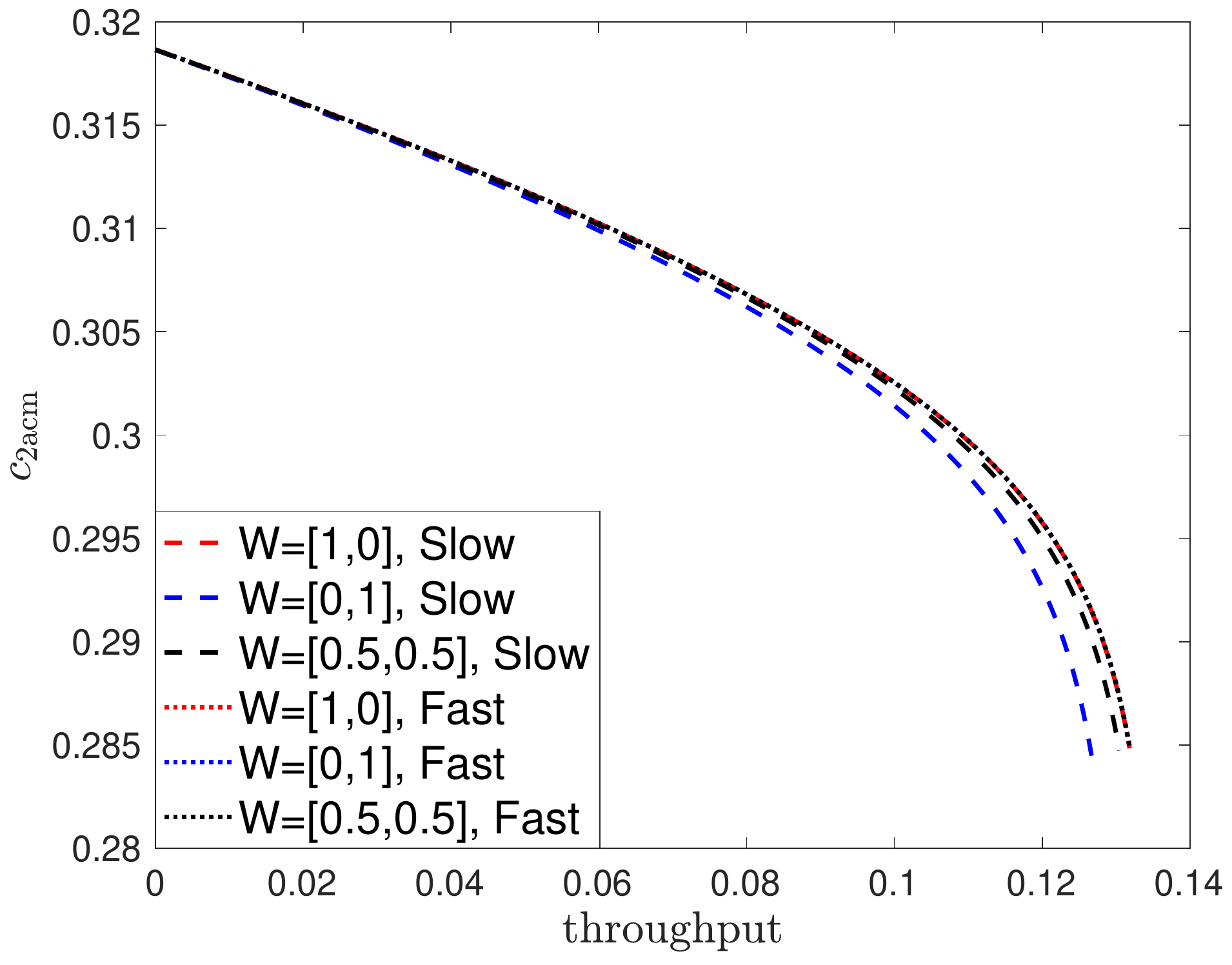}

\caption{\footnotesize{ Comparison of slow method $J(a_0)=w_1 j(t_{\rm f})+w_2c_{2\rm acm}$ (dashed curves) with fast method $J_{1, \rm fast}(a_0)=w_1 u(0)+w_1 u'(0)+ w_2c_{2\rm ins}(0)+ w_2c_{2\rm ins}'(0)$ (dotted curves)  for various weights W$=[w_1,w_2]$, with $\beta=0.1, \xi=0.5, {\color{black}\lambda_1=1}$: (a) flux vs. throughput ($u, j$) plot, (b) cumulative concentration of type 1 particles vs. throughput ($c_{1\rm acm}, j$) plot, (c) cumulative concentration of type 2 particles vs. throughput ($c_{2\rm acm}, j$) plot.}}
\label{P1_vary_W}
\end{figure}

We begin by demonstrating both slow and fast optimization methods described in \S\ref{opt_method} above. Figure~\ref{P1_f_s_pore_evl} shows the fouling evolution of the optimized membrane pores $a_0(x)$ for {\bf Problem 1} with $ w_1=1, w_2=0, \xi=0.5, \beta=0.1, \alpha_1=\alpha_2, \lambda_1=1$, and $R_1(0) \ge \Rt$, using the slow method (left panel) and fast method (right panel). The top row shows the clean, unfouled optimized pore profiles at $t=0$ (Figs.~\ref{P1_f_s_pore_evl}~(a) and (d)); the center row shows the fouling of these pores at $t=t_{\rm f}/2$ halfway through the filtration (Figs.~\ref{P1_f_s_pore_evl}~(b) and (e)); and the bottom row shows the fouled pores at termination time $t=t_{\rm f}$ (Figs.~\ref{P1_f_s_pore_evl}~(c) and (f)).  The gray region is the filter material, and the dark blue color indicates the fouling by deposited particles. The white area denotes the open pore (void), and the red center line is the axis of symmetry of the pore (which has circular cross-section). This figure illustrates that the optima $a_0(x)$ found by fast and slow methods are indistinguishable.

Figure~\ref{P1_vary_W} further compares the slow method (objective function (\ref{max_Jplusc2acm_slow}), dashed curves) and the corresponding fast method (objective function (\ref{max_Jplusc2acm_fast}), dotted curves) for \textbf{Problem 1}, with various weights $[w_1, w_2]$, weighting total throughput and cumulative type-2 particle concentration respectively, indicated in the legend. The results presented here for $[w_1,w_2]=[1,0]$ correspond to the optimized profiles presented in Fig.~\ref{P1_f_s_pore_evl}.
Results are plotted as functions of filtrate throughput over the duration of the filtration, $0\leq t\leq t_{\rm f}$. The quantities shown in Figure~\ref{P1_vary_W}, for the optima obtained using both methods, are: (a) flux vs. throughput ($u, j$) plot; (b) accumulative concentration of type 1 particles vs. throughput, ($c_{1\rm acm}, j$) plot; and
(c) accumulative concentration of type 2 particles vs. throughput, ($c_{2\rm acm}, j$) plot; all simulated with $\xi=0.5, \beta=0.1, \alpha_1=\alpha_2, \lambda_1=1$, and $R_1(0) \ge \Rt$. The figure shows that for all three sets of weights considered, $[w_1, w_2]= [1,0], [0.5, 0.5]$ and $[0,1]$, the fast method finds an optimized $a_0(x)$ as good as or better than that found by the slow method (larger or the same values for $j(t_{\rm f})$ and $c_{2\rm acm}(t_{\rm f})$, while always satisfying the removal criterion $\Rt$ for particle type 1).

These results, as well as many others not discussed here, demonstrate that with the same number of searching points the fast method converges to an optimizer that in all cases is as good as or better than that obtained using the slow method, with considerably shorter running time (a typical optimization for the slow method takes 40 minutes with 10,000 searching points, while the fast method takes only 4 minutes). We also observe that varying the weights $[w_1, w_2]$ does not change the optimized profile significantly, indicating that maximizations of  $j(t_{\rm f})$ and of  $c_{2\rm acm}(t_{\rm f})$ are correlated for the parameter values considered. One possible explanation for this correlation is that, provided the type 1 particle removal constraint $R_1(0) \ge \Rt$ is met, the initial concentration of type 2 particles in the filtrate $c_{2}(1,0)$ should be maximized by maximizing the initial flux $u(0)$, since the higher the flux, the more type 2 particles will escape capture by the filter.

Figure~\ref{P2_vary_xi_beta} presents direct comparisons of the slow and fast methods for \textbf{Problem 2}.  
\begin{figure}
{\scriptsize (a)}\rotatebox{0}{\includegraphics[scale=.39]{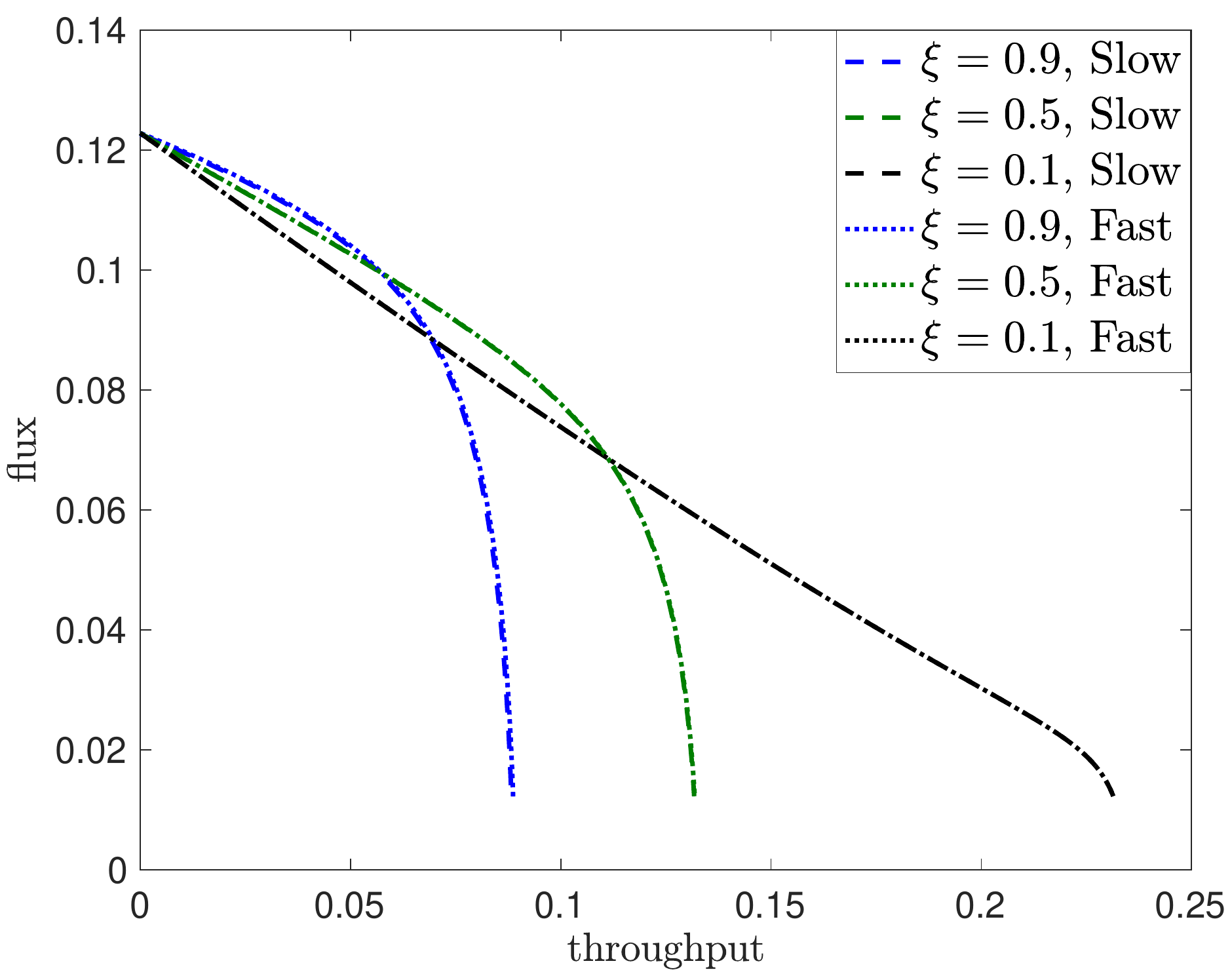}}
{\scriptsize (d)}\includegraphics[scale=.37]{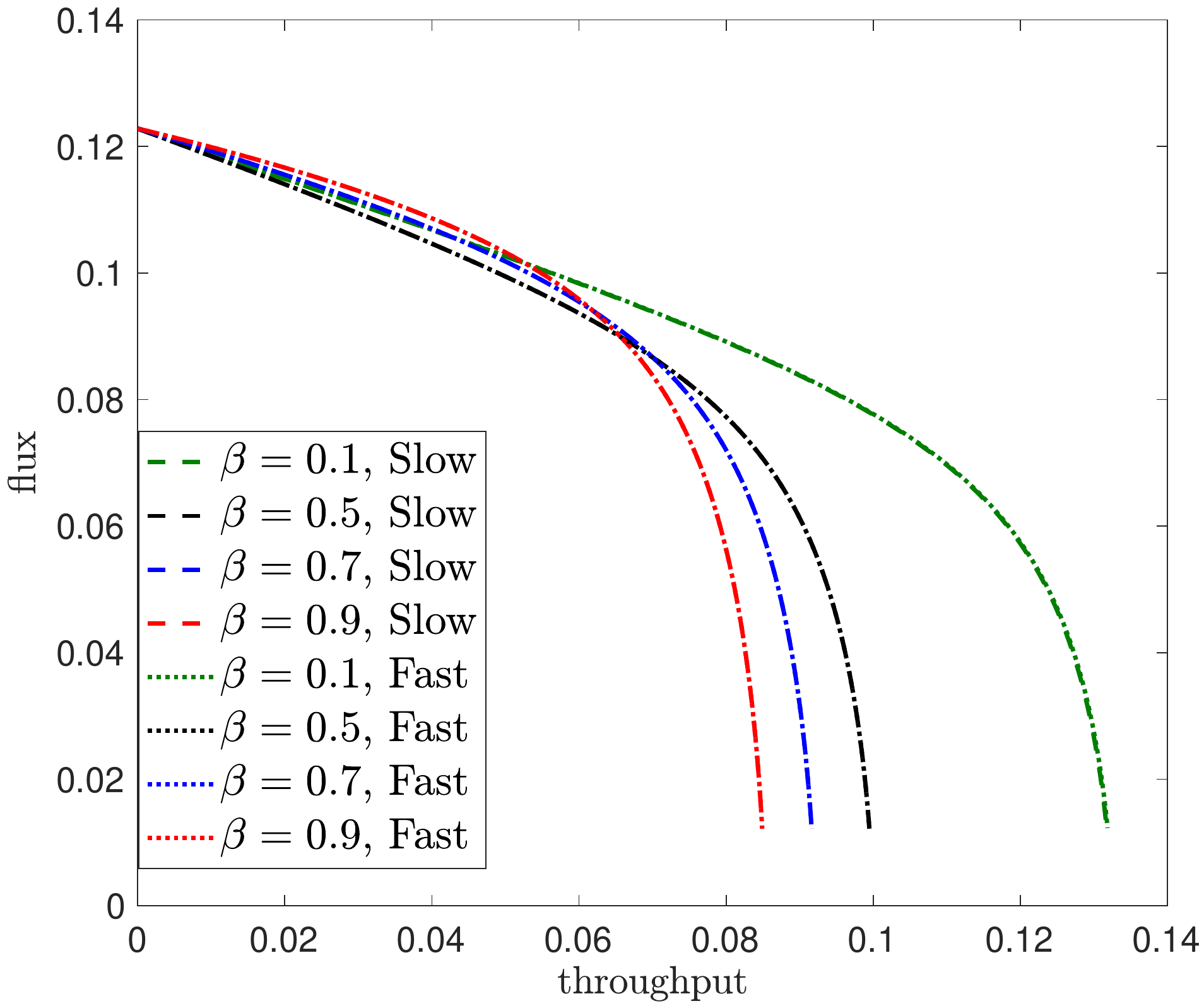}
{\scriptsize (b)}\includegraphics[scale=.39]{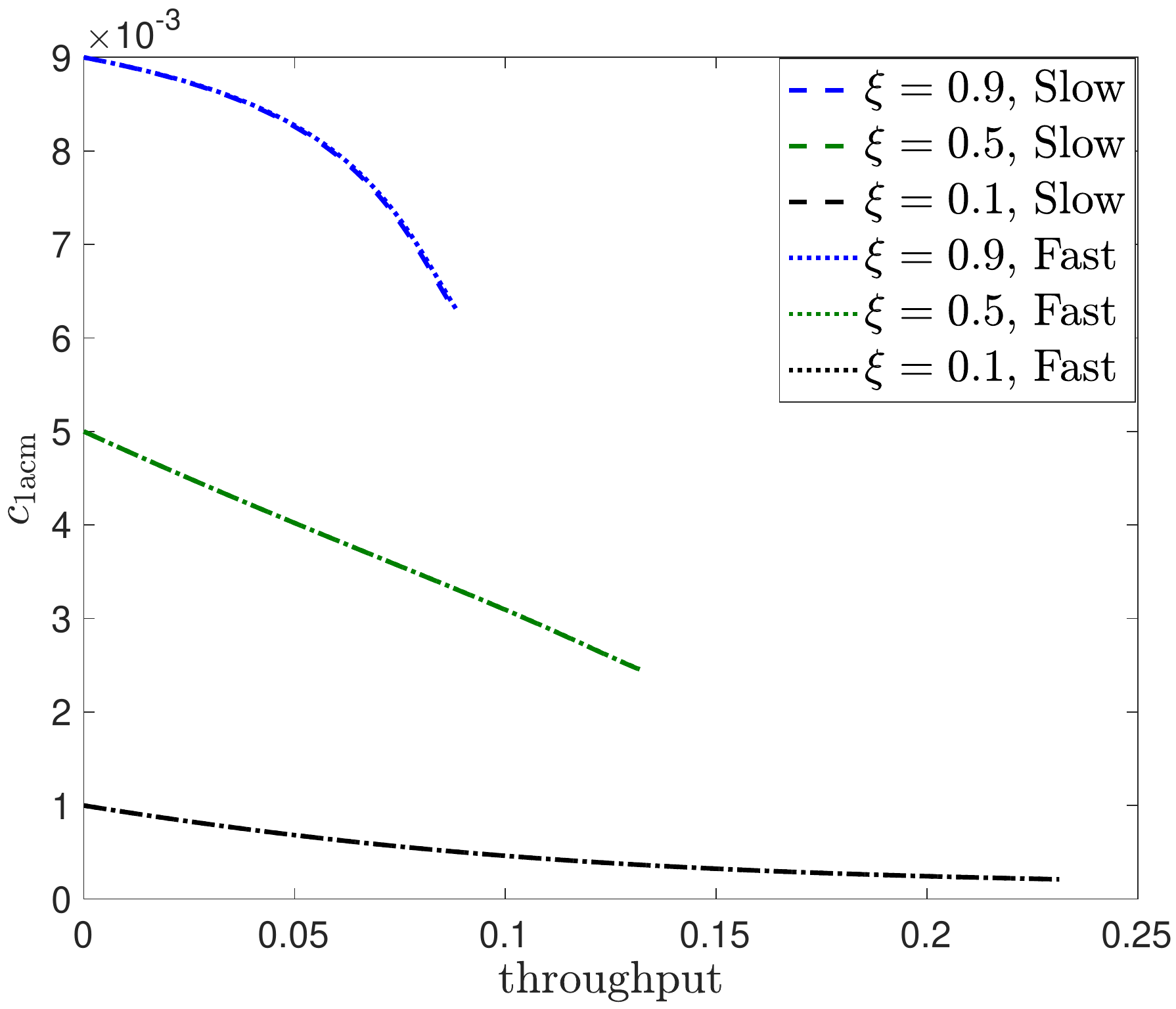}
{\scriptsize (e)}\includegraphics[scale=.39]{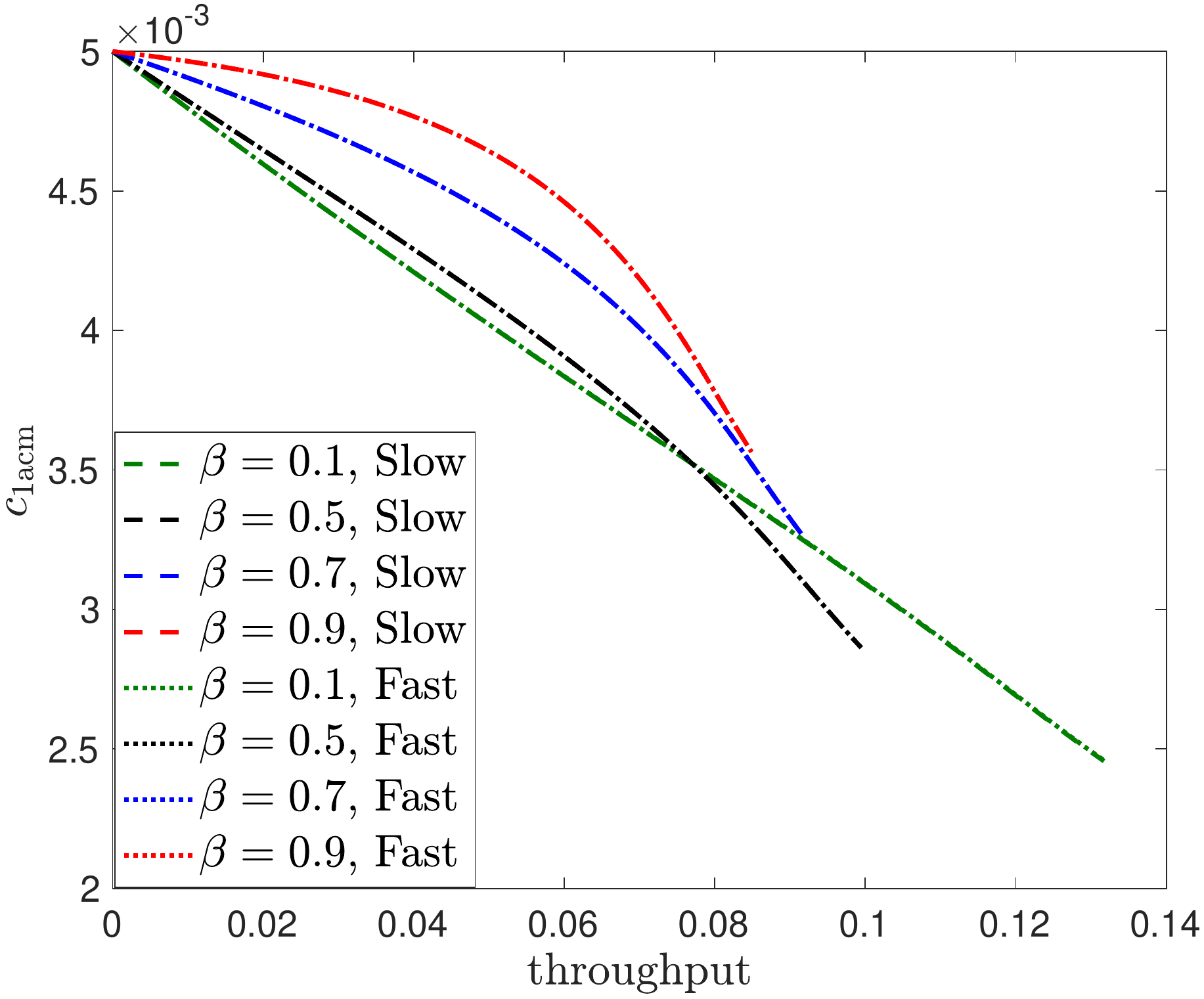}
{\scriptsize (c)}\includegraphics[scale=.39]{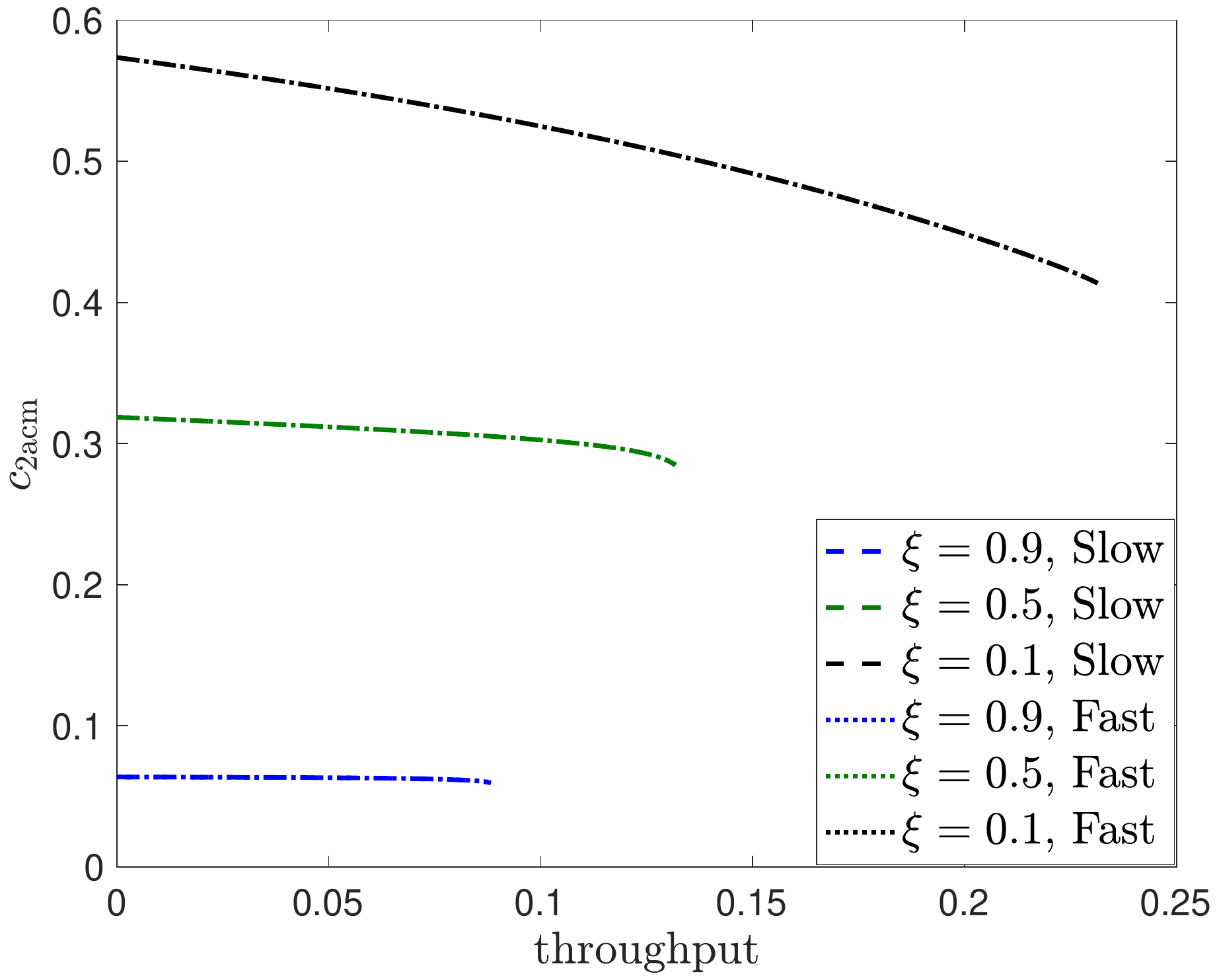}
{\scriptsize (f)}\includegraphics[scale=.39]{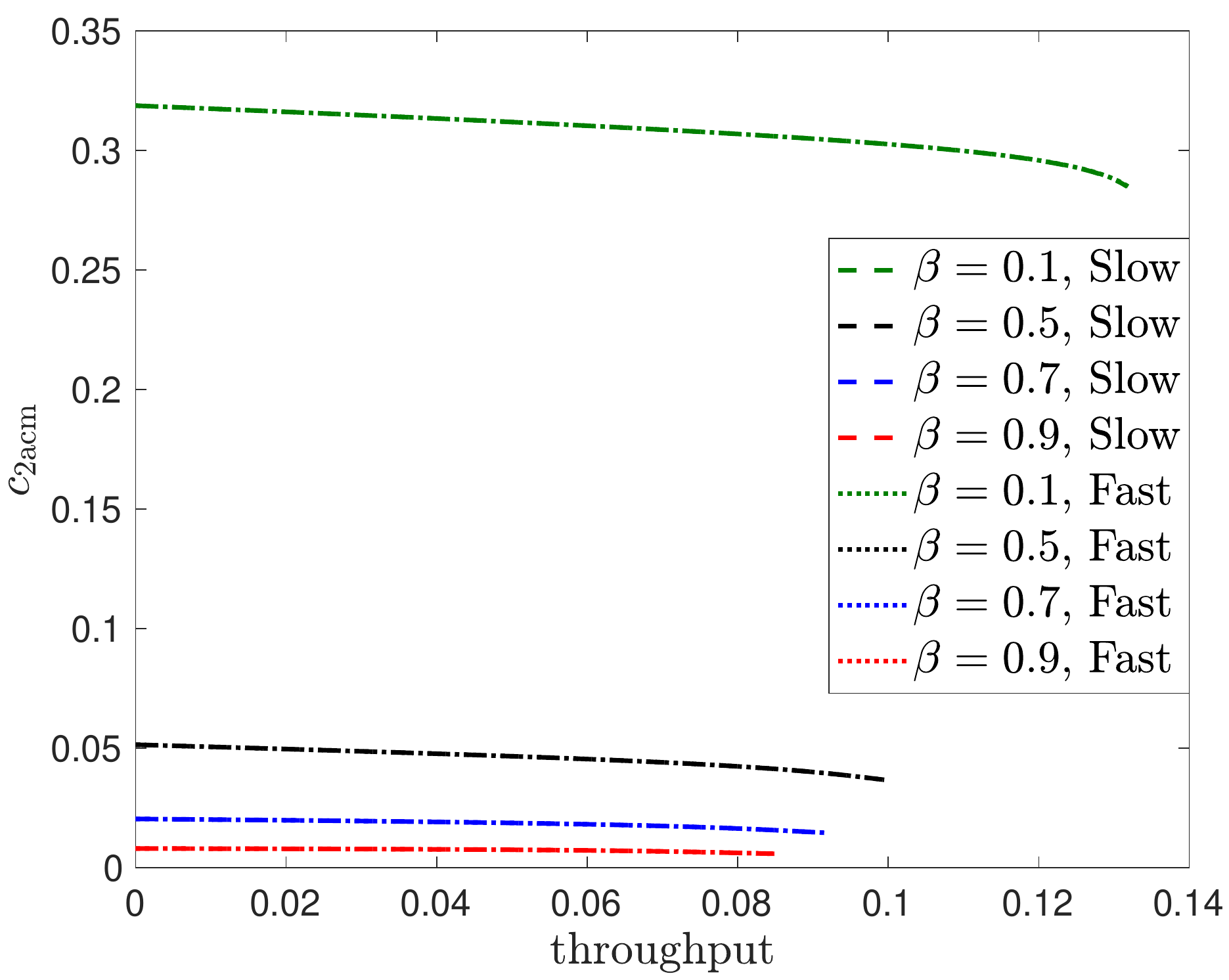}

\caption{\footnotesize{ Comparison of slow method with objective function $J(a_0)=j(t_{\rm f})c_{2\rm acm}(t_{\rm f})$ (dashed curves) and fast method with objective $J_{2, \rm fast}(a_0)=u(0)c_{2\rm ins}(0)$ (dotted curves)    (a-c): with $\xi=0.9, 0.5, 0.1$ and $\beta=0.1$, $\alpha_1=\alpha_2$, $\lambda_1=1$: (a) ($u, j$) plot, (b) ($c_{1\rm acm}, j$) plot, (c) ($c_{2\rm acm}, j$) plot. (d-f): with $\beta\in [0.1,0.9]$ varying and $\xi=0.5, \lambda_1=1$: (d) ($u, j$) plot, (e) ($c_{1\rm acm}, j$)  plot,  (f) ($c_{2\rm acm}, j$) plot.
}}
\label{P2_vary_xi_beta}
\end{figure}
Results for the slow method, with objective function (\ref{max_Jtimesc2acm_slow}), are indicated by dashed curves; and those for the corresponding fast method, with objective function (\ref{max_Jtimes2acm_fast}), by dotted curves. The left panel, Figs.~\ref{P2_vary_xi_beta}~(a-c), shows results for various feed particle-composition ratios $\xi$ (other parameters as in Fig.~\ref{P1_vary_W}); while the right panel, Figs.~\ref{P2_vary_xi_beta}~(d-f), compares results for various effective particle-membrane interaction ratios $\beta$, with $\xi=0.5$.
The flux through the membrane and the cumulative particle concentrations of each particle type in the filtrate are plotted as functions of filtrate throughput over the duration of the filtration, $0\leq t\leq t_{\rm f}$.
In all cases, for the same number of searching points, 
the fast method converges to the same optimal pore profile as the slow method across all $\xi$ and $\beta$ values considered (though the optima obtained are different for each parameter set). 
Similar to Problem 1, the computational speedup is considerable using the fast method. 

In addition to demonstrating the efficacy of the fast optimization method, the results also illustrate some general features of the model. When a feed contains a larger fraction (higher $\xi$-value) of particles to be removed (type 1 particles here), our model predicts shorter filter lifetime (due to faster fouling) when compared to a feed with lower $\xi$-value, leading to less total throughput and lower final accumulative particle concentration of type 2 particles in the filtrate, see e.g., Figs.~\ref{P2_vary_xi_beta}(a) and (c). This is not desirable if we want to maximize total collection of type 2 particles; we will present one possible way to circumvent this issue in \S\ref{result_multi-stage}, where a multi-stage filtration is proposed.  

Figure~\ref{P2_vary_xi_beta}(b) is the ($c_{1\rm acm}, j$) plot. 
Note that in all cases, the constraint for removal of particle type 1 is tight at the optimum, with the exact specified proportion $\Rt$ of particles removed from the feed at time $t=0$. Figures \ref{P2_vary_xi_beta} (d-f) show that larger $\beta$ values (meaning that the two particle types are more physicochemically similar; recall $\beta\in (0,1)$ throughout our study, and if $\beta=1$ both particle types interact identically with the membrane) lead to faster fouling of the filter, with lower total throughput and lower total yield of type 2 particles in the filtrate. This confirms our expectation that the more similar the particle types are, the more challenging it is to separate them by filtration. To achieve effective separation, a sufficient physicochemical difference $\gamma=1-\beta$ is required. 

Encouraged by the excellent results and significant speedup obtained when using the fast method with the same number of searching points as the slow method, we next investigate its performance with fewer searching points. Figure \ref{fast_method_vary_search_pt} shows the comparison between the slow method (dashed curves) with 10,000 searching points (found, empirically, to be the minimum number required for reliable results) and the fast method (dotted curves) with decreasing number of search points (10,000, 1,000, 100). Model parameters are fixed at $\xi=0.5$ and $ \beta=0.1$; other parameters are as in Fig.~\ref{P1_vary_W}. These results (as well as many other tests, not shown here) indicate that the fast method produces reliable results with just 1,000 searching points (blue dotted curves; this optimum even appears superior to the slow method with 10,000 search points, providing slightly higher total throughput). Even with as few as 100 search-points the fast method produces reasonable, though suboptimal, results (black dotted curve). In all cases the particle removal constraint on $c_1$ is again tight at the optima found. Since run time for the optimization routine scales in direct proportion to the number of searching points, a 10-fold reduction in the number of search points needed represents a significant additional computational saving: the fast method with 1,000 search points is approximately 100 times faster than the slow method with 10,000 points. 
  
\begin{figure}
{\scriptsize (a)}\rotatebox{0}{\includegraphics[scale=.4]{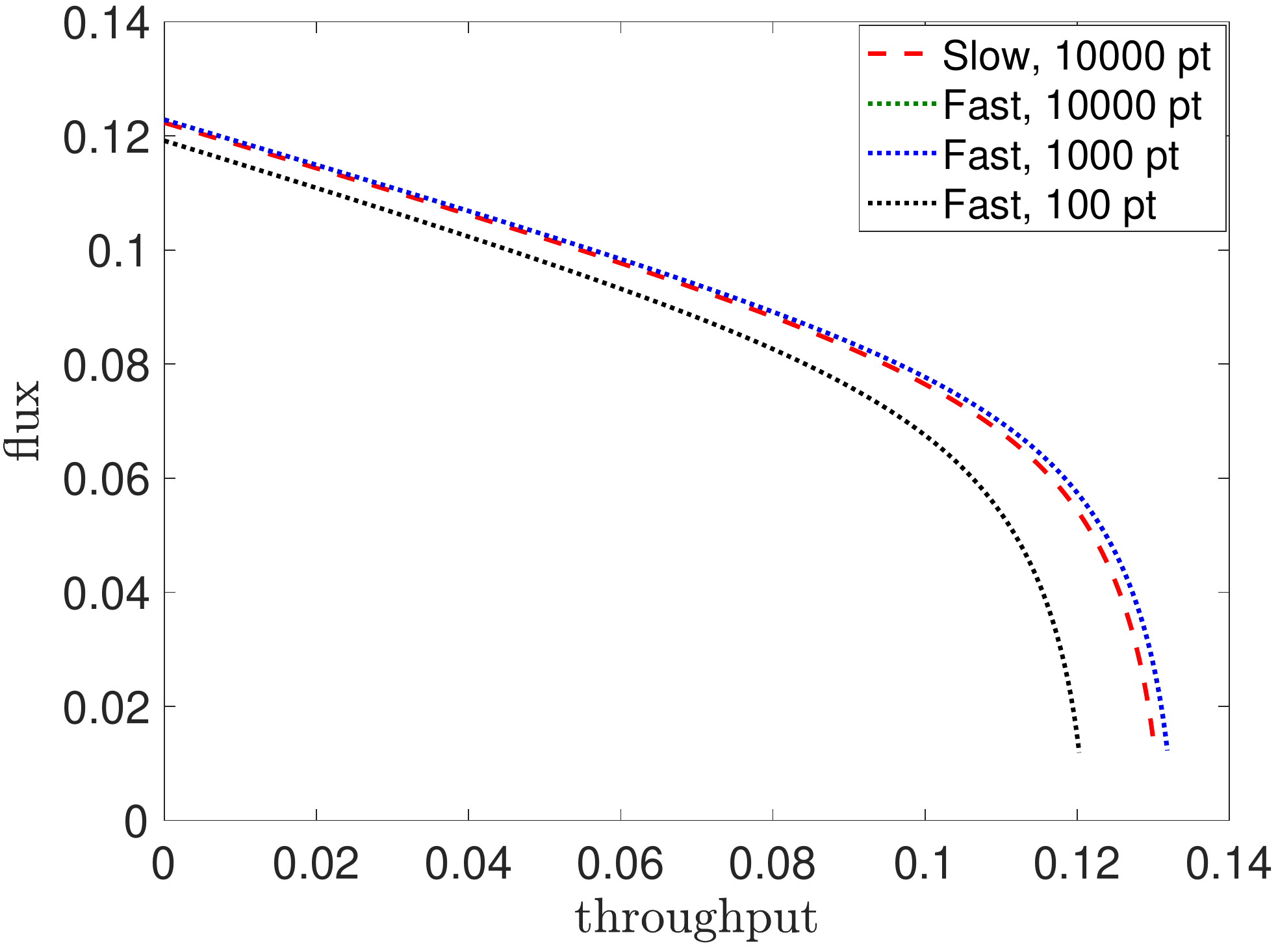}}
{\scriptsize (b)}\includegraphics[scale=.38]{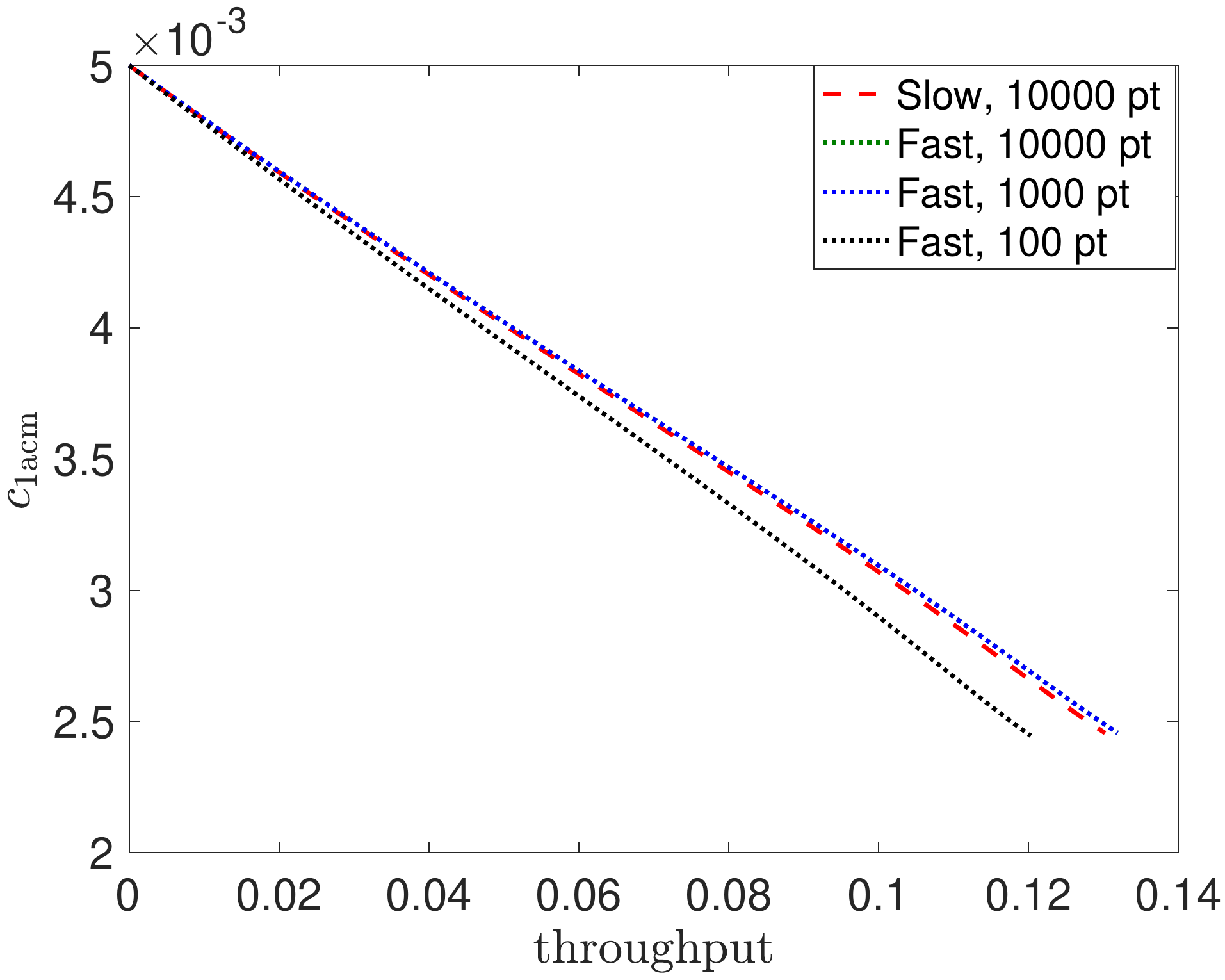}
{\scriptsize (c)}\includegraphics[scale=.39]{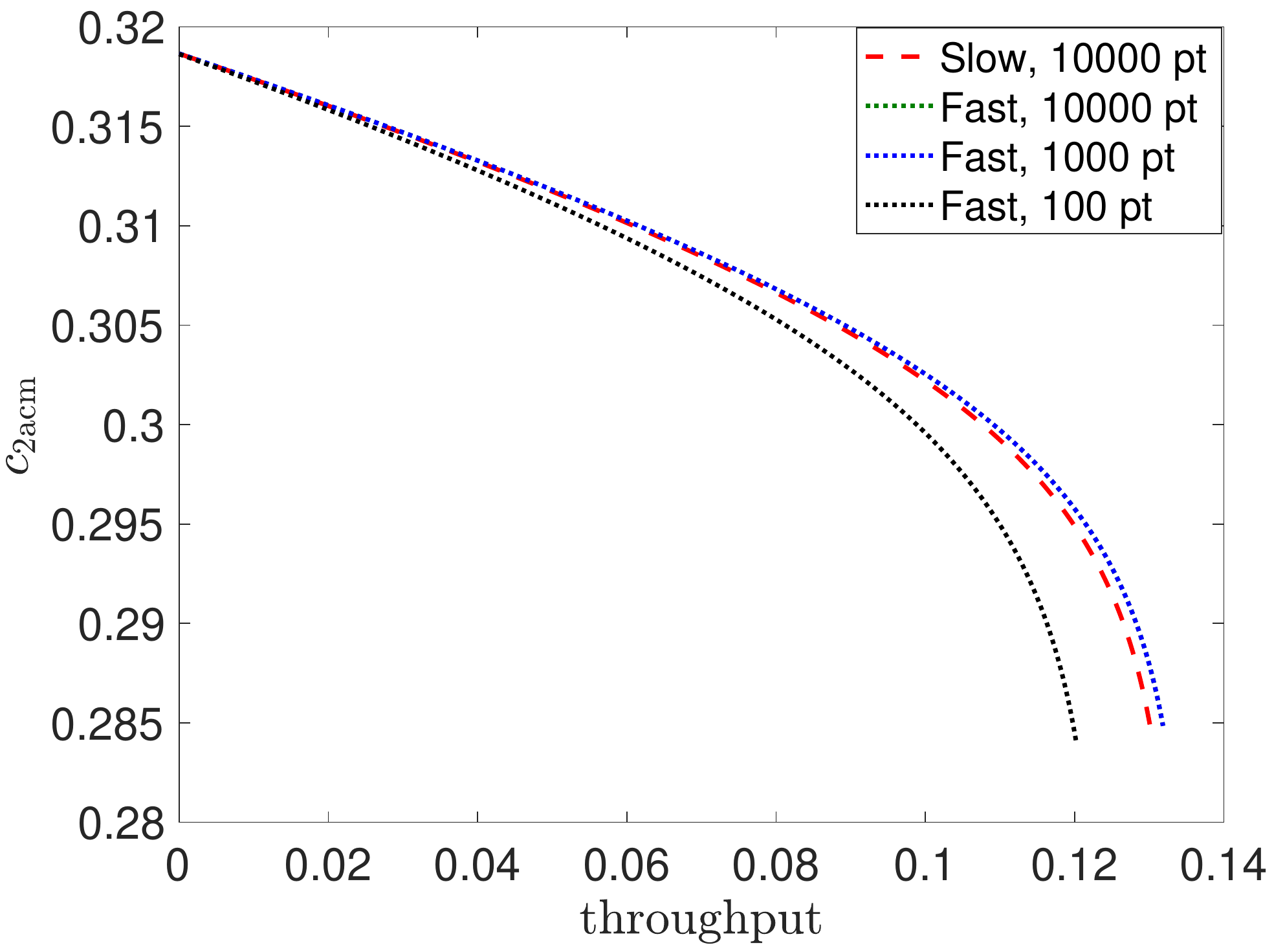}


\caption{\footnotesize{Comparison of slow method $J_{\rm slow}(a_0)=j(t_{\rm f})c_{2\rm acm}(t_{\rm f})$ (dashed curves) with 10,000 start points, with fast method $J_{\rm fast}(a_0)=u(0)c_{2\rm ins}(0)$ (dotted curves) using variable number of searching points, with $\xi=0.5, \lambda_1=1$ and $ \beta=0.1$: (a) ($u, j$) plot; (b) ($c_{1\rm acm}, j$) plot, and (c)  ($c_{2\rm acm}, j$) plot.
}}
\label{fast_method_vary_search_pt}
\end{figure}

\subsubsection{Multi-stage filtrations}\label{result_multi-stage}

In Fig. \ref{P2_vary_xi_beta} (a) we observed that, with a higher concentration ratio of type 1 particles in the feed, the optimized filter for a single-stage filtration tends to be fouled faster, which leads to lower total throughput per filter. This makes sense as the filter needs to remove a higher mass of impurity (type 1 particles) to achieve the initial particle removal threshold $R_{1}(0) \ge \Rt$ when $\xi$ is larger. Our simulations also reveal that the fouling is largely confined to a narrow region adjacent to the upstream surface of the filter at optimum, with the majority of the downstream portion of the filter remaining unused (see Fig. \ref{P1_f_s_pore_evl}). In this section we propose a multi-stage filtration scenario that has the potential to alleviate some of these inefficiencies.  Heuristically, we would like to process more feed per filter by increasing the membrane porosity 
and simultaneously make more efficient use of the membrane material by fouling a higher proportion of the pore (void) volume.  However, increasing the porosity in general decreases the particle removal efficiency, so the filtrate will require further purification to satisfy the particle removal requirement. 
This provides the motivation for the proposed multi-stage filtration strategy: we will lower the initial particle removal requirement to increase the amount of feed processed per filter and try to satisfy the particle removal requirement by filtering the collected filtrate again, possibly more than once (multi-stage). The multi-stage filtration will be cost-effective if the increase in feed processed can offset the increase in the number of additional filters required to meet the particle removal requirement. From the optimization point of view, we increase the feasible searching space by relaxing the initial particle removal constraint so that a better optimizer might be found. 

In the following discussion we focus on the optimization {\bf Problem 2}, where the goal is to maximize the yield of type 2 particles per filter used, while achieving {\it effective separation},\footnote{The term ``effective separation'' has been used in the literature, though without clear quantitative definition,
e.g.~\cite{shi2020, jin2020, ma2020}.} which for definiteness we here define 
as removing the desired fraction $\Rt$ of particle type 1 from the feed ($\bar R_1(\rm t_f) \ge \Rt$) while simultaneously
recovering a minimum desired yield fraction $\Upsilon$ of type 2 particles in the filtrate (${\bar R}_2 (t_{\rm f}) \le 1-\Upsilon$). For all of our simulations, $\Rt=0.99$ and $\Upsilon=0.5$.
We define the {\it purity} for type $i$ particles in the filtrate, $k_i \in [0,1]$, as 
\be
k_i =\frac{c_{i, \rm acm} (t_{\rm f})}{ \sum_{i=1,2} c_{i, \rm acm} (t_{\rm f})},
\qquad i=1,2, 
\label{purity}
\ee
where $c_{i, \rm acm} (t_{\rm f})$ is the accumulative concentration of the type $i$ particle in the filtrate at the end of the filtration. With our hypothesized scenario of feed containing desired (type 2) and undesired (type 1) particles in mind, we note a simple relationship between the purity of type 2 particles and the final cumulative removal ratios:
\be
k_2=\frac{(1-\xi)[1-{\bar R}_2 (t_{\rm f})]} {\xi [1-{\bar R}_1 (t_{\rm f})]+(1-\xi) [1-{\bar R}_2 (t_{\rm f})]}.
\label{k_R}
\ee
We will return to these definitions in our discussion of the multi-stage filtration results below. 

The basic idea behind our multi-stage filtration is to first optimize the filter with a less strict initial type 1 particle removal requirement (i.e., we require $R_1(0) \ge R <\Rt $) and filter the feed solution two or more times to achieve a larger total {\it yield per filter} of purified type 2 particles than in a single-stage filtration, with the effective separation condition satisfied at the end of the multi-stage filtration. 
We determine the {\it stage} of filtration by how many times the solution has passed through clean filters: for example, the clean stage 1 filter will take feed directly and be used to exhaustion; the filtrate collected from the stage 1 filter will then be sent through a new (clean) stage 2 filter (which may be used more than once within stage 2).

We propose the following two-stage or multi-stage filtration strategy for {\bf Problem 2} (also summarized as a flow chart in Figure \ref{flowchart_multistage}): {\bf 1.} Optimize the filter with a less strict initial particle removal requirement than desired ($R_1(0) \ge R <\Rt$); denote the optimized filter as  $F_{R}$ (e.g., for $R=0.5$, we denote the optimized filter as  $F_{0.5}$). {\bf 2.} (Stage 1) Run the filtration simulation using $F_{R}$ until the filter is completely fouled; collect the filtrate. In scenarios to be considered later we allow stage 1 to use several filters simultaneously. {\bf 3.} Re-filter the collected filtrate through another clean $F_{R}$; collect the new filtrate. {\bf 4.} Test the filtrate from step 3. Does it meet required type 1 particle final removal requirement ${\bar R}_1 (t_{\rm f}) \ge \Rt$? If yes, we are done; if no, repeat step 3 using the same filter until the requirement is met, or until this  $F_{R}$ is completely fouled. {\bf 5.} {\bf (multi-stage)} If $F_{R}$ is completely fouled before ${\bar R}_1 (t_{\rm f}) \ge \Rt$, use another clean $F_{R}$ and repeat step 3 until ${\bar R}_1 (t_{\rm f}) \ge \Rt$.
{\bf 6.} Once the threshold ${\bar R}_1 (t_{\rm f}) \ge \Rt$ is met, record the total mass of compound 2 in the filtrate and the number of $F_{R}$ used, to compute the mass yield per filter. 

In order to keep track of the number of filters used and how many times each is reused, we identify each (stage) filter used by $F_{R,m}$, e.g., the second (stage) $F_{0.5}$ filter will be denoted $F_{0.5,2}$, and we track how many times each stage filter has been used by $n(m)$, e.g., $n(2)$ denotes the number of times $F_{0.5,2}$ is used. In cases where stage 1 involves more than one clean filter (used simultaneously) we also use notation $l_m$ to denote the total number of $m$th stage filters used, e.g., $l_2$ is the number of stage 2 filters used; and we denote the total number of filters used at the end of multi-stage filtration by $M=\sum_{m} l_m$. After the multi-stage filtration is concluded we calculate the total mass of type 2 particles collected per filter used, $c_{2 \rm acm}(t_{\rm f}) j(t_{\rm f})/M$, and compare it with the collected mass from the filter $F_{R}$ optimized for single-stage filtration. 
If the mass collected per filter used is larger for the multi-stage filtration, then the process is deemed more cost-effective. 
\begin{figure}
	\begin{tikzpicture}[node distance=2cm]
		\node (start) [startstop] {Start};
		\node (in0) [io, below of=start] {Optimize to obtain $F_{R}$ with $R<\Rt$};
		\node (in1) [io, below of=in0] {Stage 1: set $m=1$; run filtration with  $F_{R , m}$ until termination; collect filtrate; set $n=1$.};
		\node (pro1) [process, below of=in1,yshift=-0.5cm] {Record $(m,n)$; set $n = 0$, let $m\to m+1$; take a new $F_{R}$ and register it as $F_{R , m}$.};
		\node (pro2) [process, below of=pro1,yshift=-0.5cm] {Filter the filtrate with $F_{R ,m}$; let $n\to n+1$; collect filtrate.};
		\node (dec1) [decision, below of=pro2] {${\bar R}_1 (t_{\rm f}) \ge \Rt$? };
		\node (pro3) [process, below of=dec1, yshift=0.3cm] {Record $(m,n)$, $c_{i, \rm acm}(t_{\rm f})$ and $j(t_{\rm f})$.};
		\node (end) [startstop, below of=pro3, yshift=0.3cm] {End};
		\node (dec2) [decision, right of=pro2, xshift=4.5cm] {Is  $F_{R, m}$ completely fouled? };
		
		\draw [arrow] (start) -- (in0);
		\draw [arrow] (in0) -- (in1);
		\draw [arrow] (in1) -- (pro1);
		\draw [arrow] (pro1) -- (pro2);
		\draw [arrow] (pro2) -- (dec1);		
		\draw [arrow] (pro3) -- (end);
		\draw [arrow] (dec1) -- node[anchor=east] {yes} (pro3);
		\draw [arrow] (dec1) -- node[anchor=north] {no} (dec2);
		\draw [arrow] (dec2) -- node[anchor=north] {no} (pro2);
		\draw [arrow] (dec2) -- node[anchor=north] {yes} (pro1);

	\end{tikzpicture}
	\caption{Flow chart of multi-stage filtration. $F_{R,m}$ signifies the $m$th stage filter used; for each $F_{R ,m}$, $n(m)$ records how many times the filter is used.}
	\label{flowchart_multistage}
\end{figure}
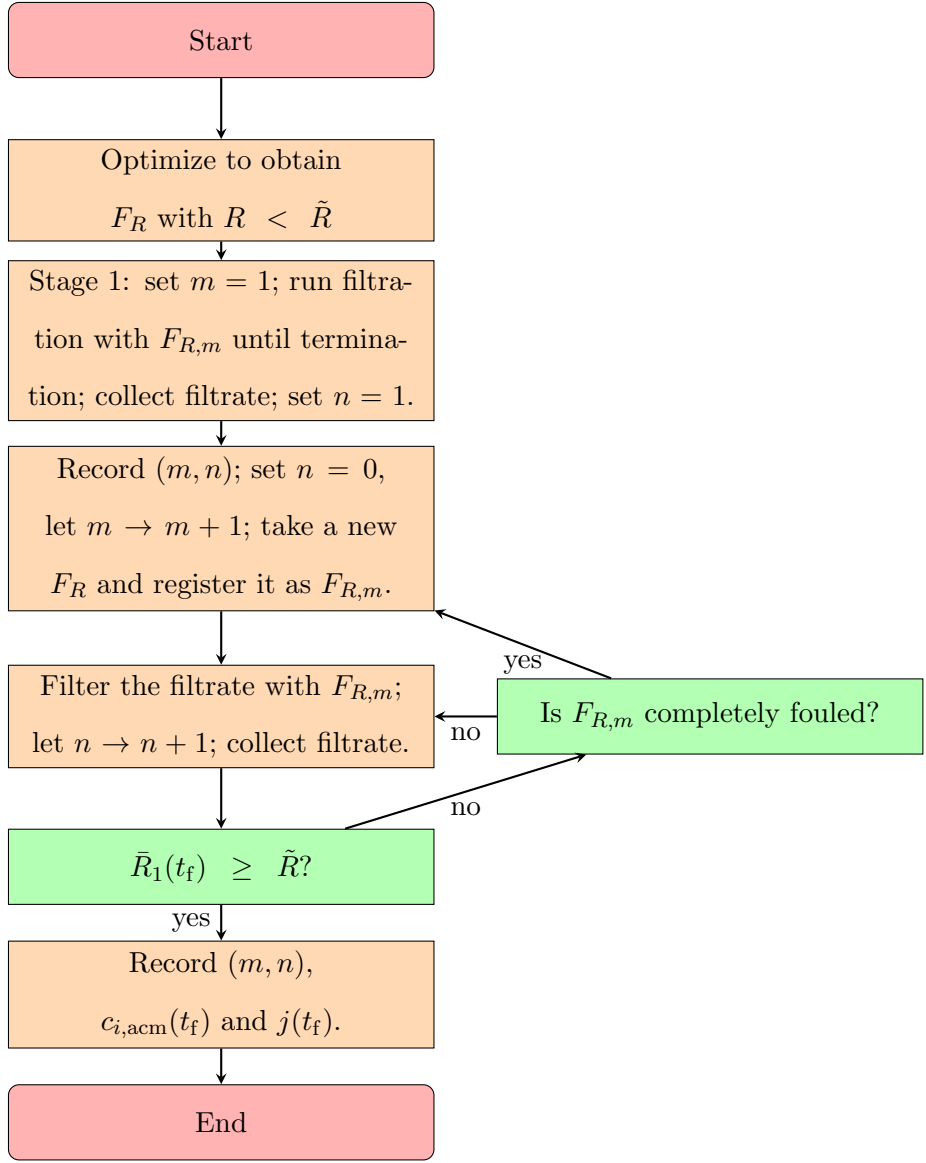

\begin{table}
\centering
\begin{tabular}{|m{1.5cm}|m{1cm}|m{1cm}|c|c|c|c|c|}
\hline
{\bf $R$} & {\bf $M$}  & {\bf $n(2)$} & {\bf $\bar R_1(\rm t_f)$} & {\bf $\bar R_2(\rm t_f)$} & {\bf $k_2$}  & {\bf $j(\rm t_f)$}  & {\bf mass yield/filter}\\
\hline
$ \Rt (=0.99) $ & $1$ &n/a &0.993 & 0.404 &0.904  &0.089 &  0.00528\\
$ 0.7 $ & $2$  & 3 &0.997 & 0.455 &0.954  &0.217 & 0.00591\\
$ 0.5 $ & $2$ &4  &0.995 & 0.427 &0.935  & 0.316 & 0.00905 \\
\hline
\end{tabular}
\caption{\footnotesize{Comparisons of single-stage filtration ($R=\Rt$) with 2-stage filtrations ($R =0.7$ and $R=0.5$). We record: $M$, the total number of filters used for each filtration process; $n(2)$, the number of times the $2$nd filter is used for each multi-stage filtration process; $\bar R_1(\rm t_f)$ and $\bar R_2(\rm t_f)$, the final cumulative particle removal ratios for particle types 1 and 2 respectively; $k_2$, the purity of type 2 particles in the final filtrate; $j(\rm t_f)$, total throughput; and type 2 particle mass yield per filter.}
}
\label{2t:2stage}
\end{table}

In table \ref{2t:2stage}, we list results comparing a single-stage filtration using a filter $F_{\Rt}$ (optimized for particle type 1 initial removal threshold set at the desired value $\Rt$), with two separate two-stage filtrations using filters $F_{0.7}$ and $F_{0.5}$ (optimized for lower thresholds $R=0.7$ and $R=0.5$, respectively). In all cases fast optimization was carried out using objective function $J_{2, \rm fast}(a_0)=u(0)c_2(0)$ with $\xi=0.9, \beta=0.1, \alpha_1=\alpha_2$, and $\lambda_1=1$.  The quantities listed in table \ref{2t:2stage} are: $R$, initial type 1 particle removal threshold; $M$, total number of filters used in each case; $n(2)$, the number of times the $2$nd (stage 2) filter is used; ${\bar R}_1 (t_{\rm f}) $ and ${\bar R}_2 (t_{\rm f}) $, the
final cumulative particle removal ratios  for particle types 1 and 2 respectively; $k_2$, the purity of type 2 particles in the final collected filtrate; $j(\rm t_f)$, the total throughput; and the total mass yield of purified type 2 particles per filter (all relevant quantities are defined in Table~\ref{performance_metrics}). 

These preliminary results show that, when our multi-stage filtration protocol is applied, we can achieve the same final particle removal requirement ${\bar R}_1 (t_{\rm f}) \ge \Rt$ as the single stage filtration but with much higher yield per filter of particle type 2; the third row of Table~\ref{2t:2stage} shows that, with $R =0.5$ the yield of purified type 2 particles per filter is almost doubled when compared to the single-stage filtration optimized for $R=\Rt$. From table \ref{2t:2stage} we also observe that the multi-stage filtrations improve the purity of the filtrate as indicated by the $k_2$ values. We note that all three filtrations achieve effective separation according to our (somewhat arbitrary) definition, which corresponds to purity $k_2\ge0.847$ for the cases considered in table \ref{2t:2stage}, i.e. $\xi=0.9$. 
If higher purity is desired to consider a separation effective,  the removal ratios ${\bar R}_1 (t_{\rm f})$ and ${\bar R}_2(t_{\rm f})$ can be adjusted accordingly based on Eq. (\ref{k_R}).

\begin{figure}
{\scriptsize (a)}\rotatebox{0}{\includegraphics[scale=.38]{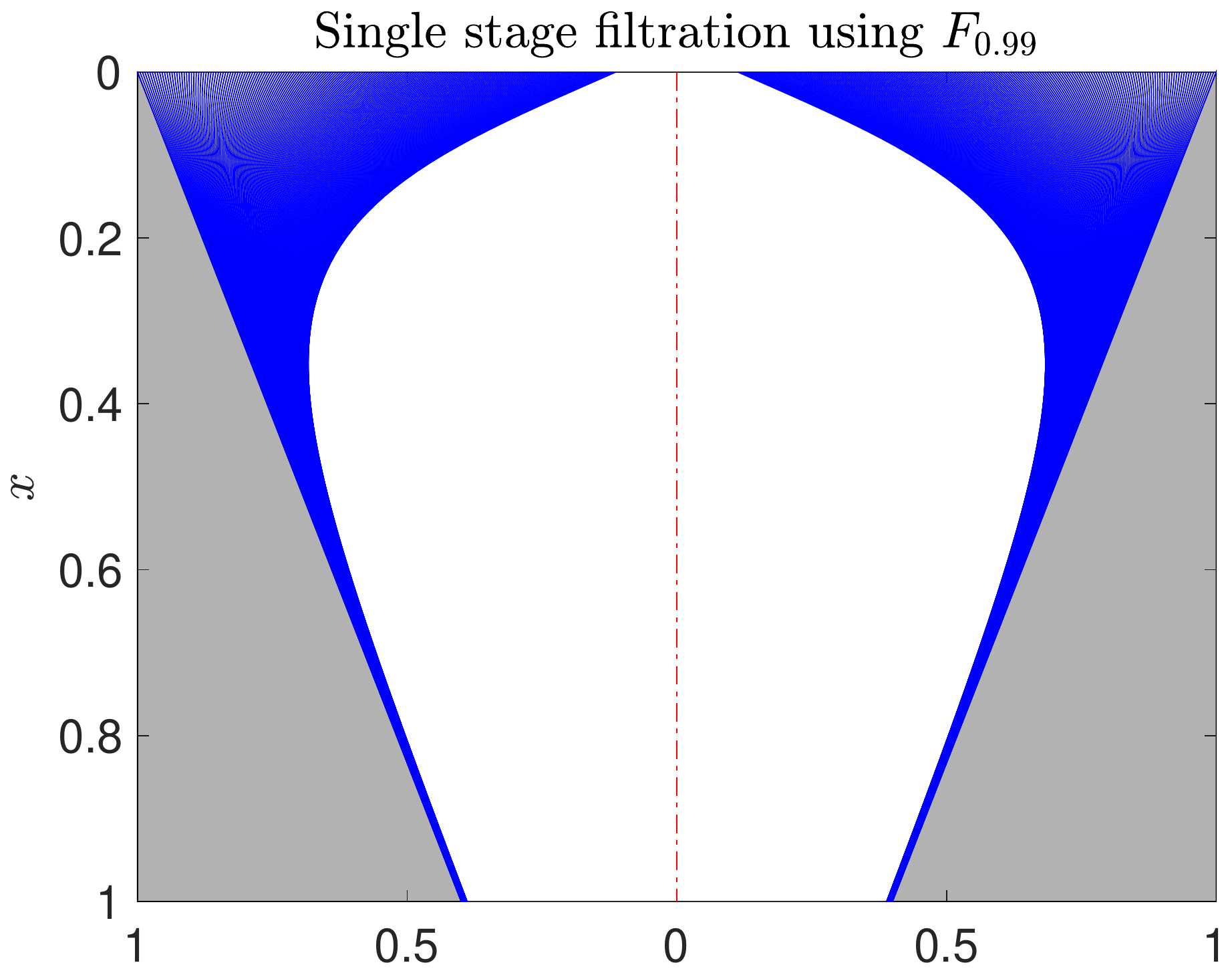}}\hspace{3cm}
{\scriptsize (b)}\includegraphics[scale=.39]{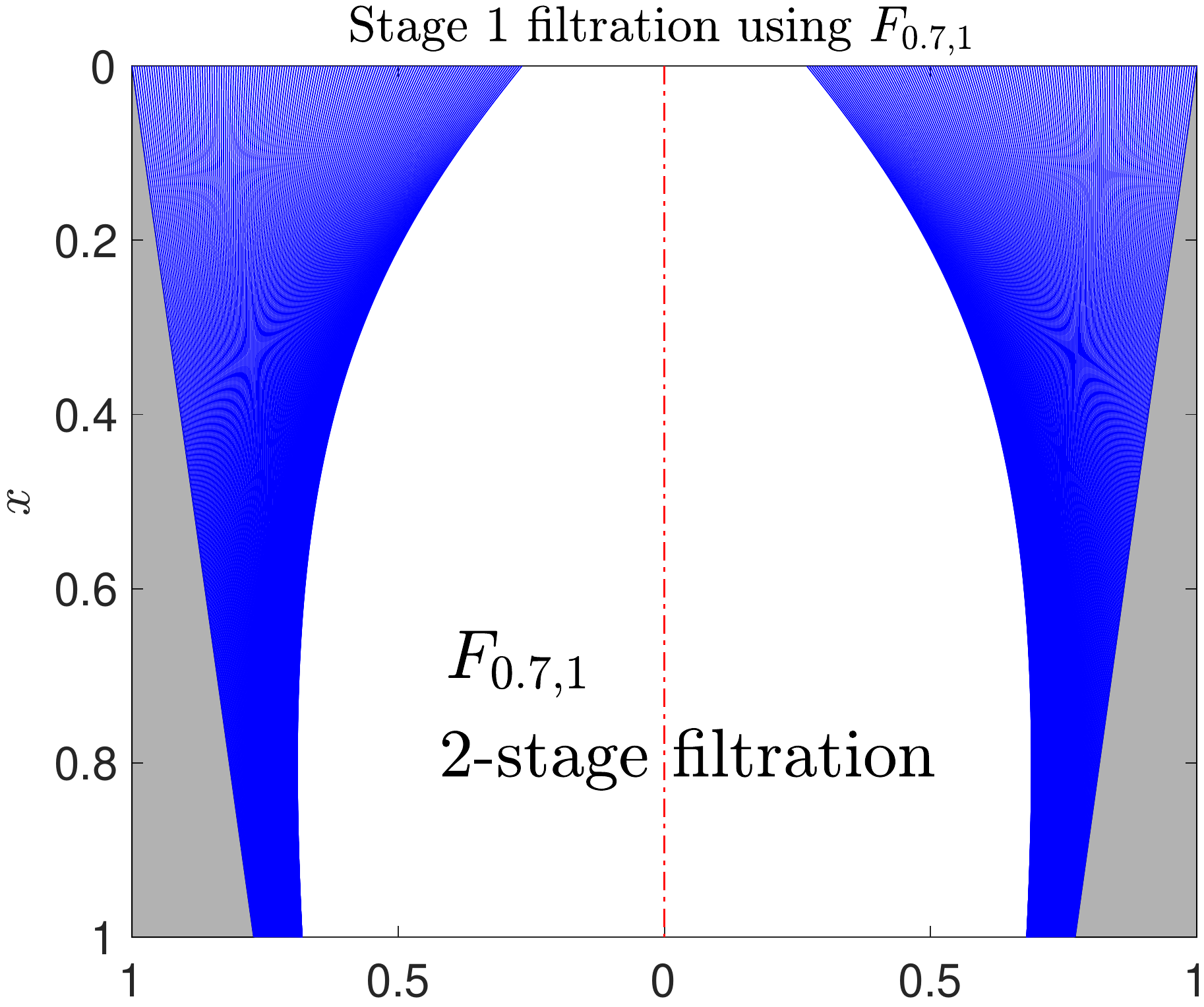}
{\scriptsize (d)}\includegraphics[scale=.39]{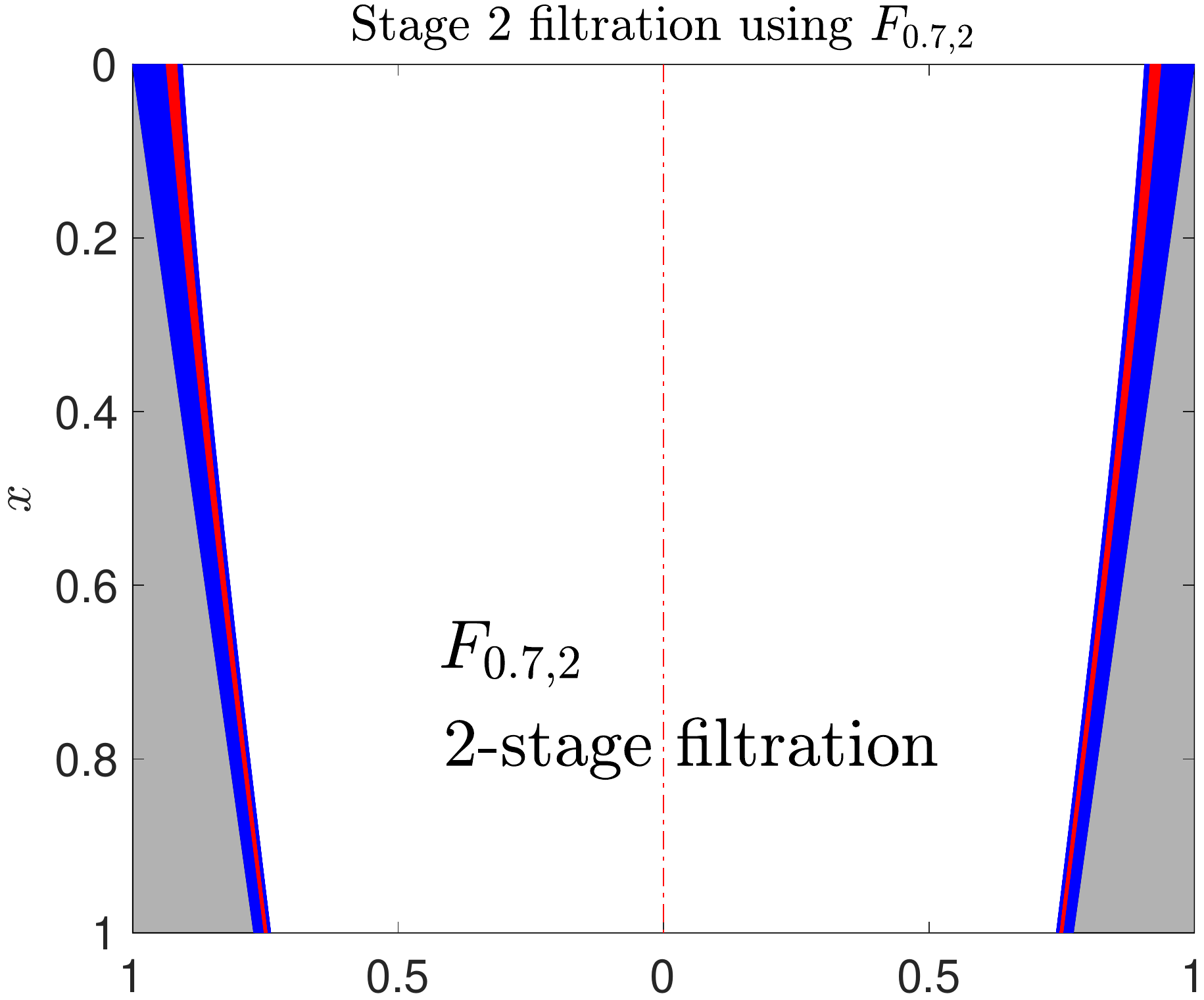}
{\scriptsize (c)}\includegraphics[scale=.39]{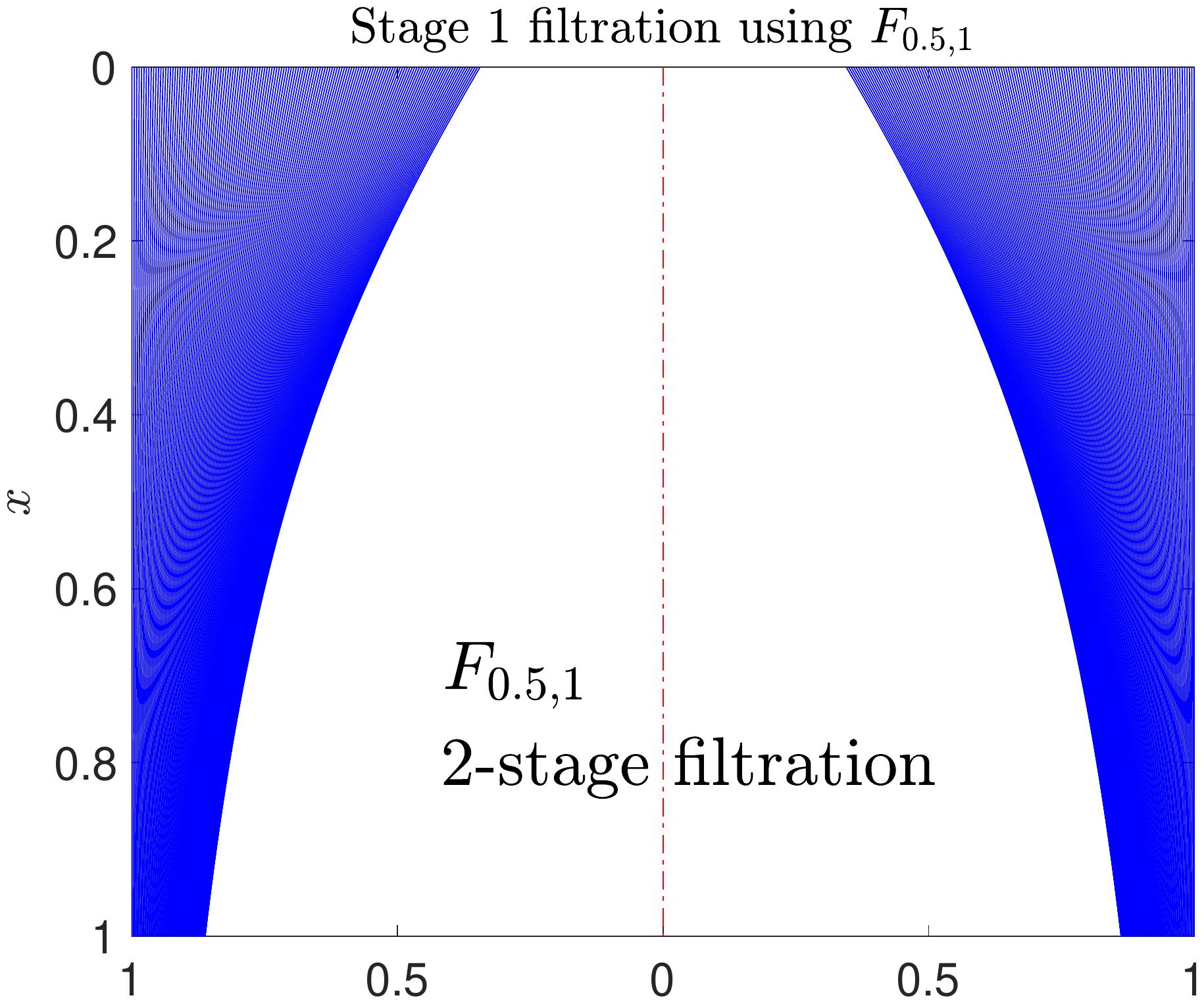}
{\scriptsize (e)}\includegraphics[scale=.39]{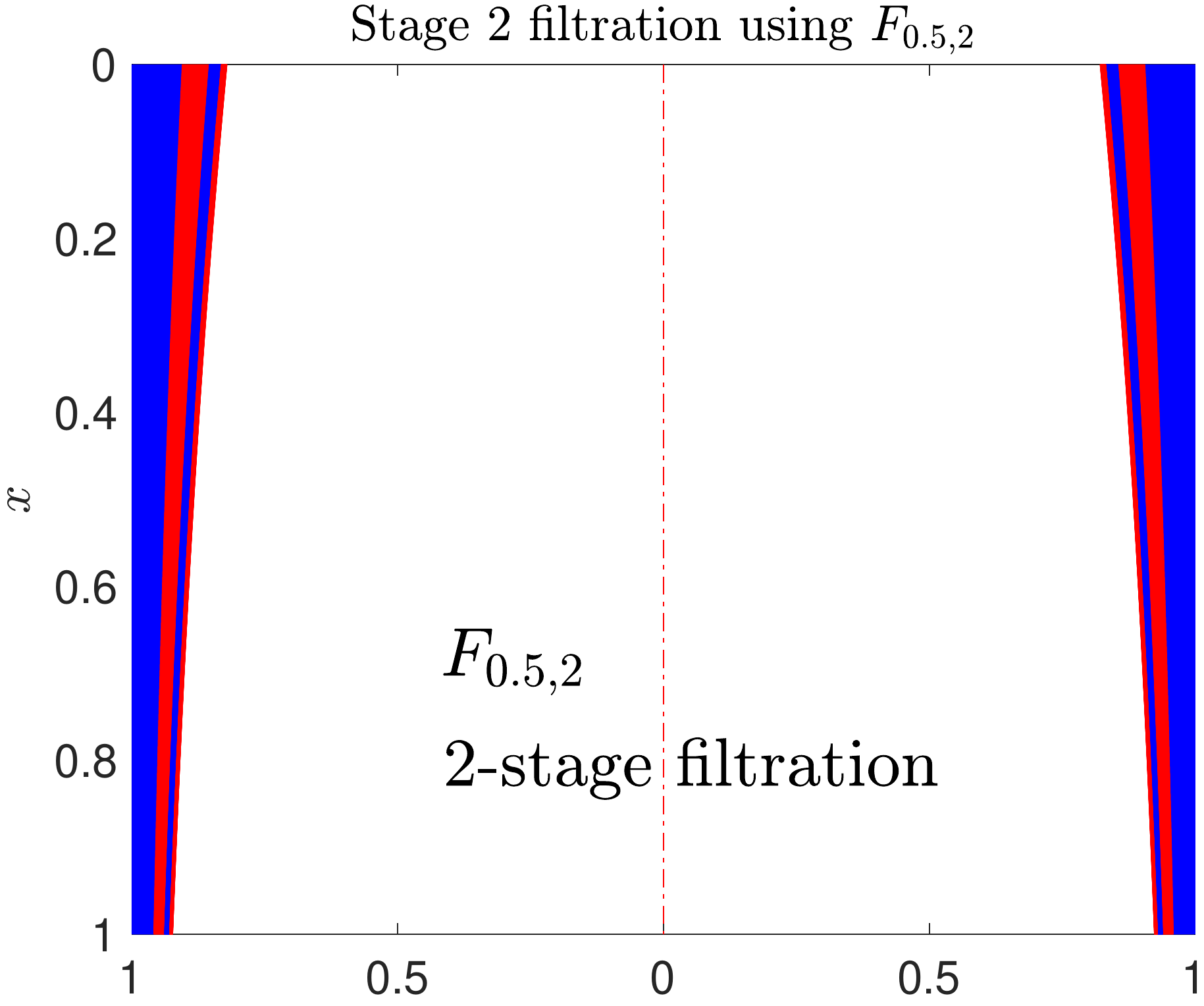}
\caption{\footnotesize{Comparison of single stage filtration and two-stage filtrations. (a-c) show completely fouled filters $F_{R_{},1}$ optimized for: (a) $R =\Rt$ (single-stage filtration), (b) $R =0.7$ and (c) $R =0.5$, with other parameters fixed at $\xi=0.9, \beta=0.1, \alpha_1=\alpha_2, \lambda_1=1$. Gray color indicates membrane material, blue is deposited particles, and white is void. (d, e) show the fouling sequence for the second filtration stages, required when $R < \Rt$: (d) filtrate from (b) is passed repeatedly through $F_{0.7,2}$ and (e) filtrate from (c) is passed repeatedly through $F_{0.5,2}$, with alternating blue and red indicating deposited particles from the successive filtrations (filter reuse). Full details in text.
}}
\label{multi-stage_vary_R}
\end{figure}

Figure~\ref{multi-stage_vary_R} illustrates the results for the optimized filters summarized in table \ref{2t:2stage} and described above via the fouling evolution of the filter pores. 
Figures \ref{multi-stage_vary_R} (a-c) show the filters from the first filtration stage, optimized for particle removal thresholds $R_{}=\Rt$ (a), $R= 0.7$ (b), and $R= 0.5$ (c), at time $t=t_{\rm f}$ (when the flux is reduced to the fraction $\vartheta=0.1$ of its initial value). The blue and red colors indicate deposited particles; a change of color indicates reuse of the filter. 
Figures \ref{multi-stage_vary_R} (d) and (e) show the fouling of the second stage filters, $F_{0.7,2}$ and $F_{0.5,2}$ respectively. 
We can see that when the initial removal threshold $R$ is decreased, the fouling of the pore becomes more uniform along its depth, and the porosity of the corresponding optimized filter $F_{R}$ increases. 
For the case $R=0.5$, the optimized pore profile is almost as wide as possible; the gray colored region corresponding to the membrane material is too thin to be visible. The high mass yield per filter and small quantity of membrane material required to produce $F_{0.5}$ indicates that if the membrane material has good selectivity (higher $\lambda_1$ value in our model), it might be advantageous to focus on maximizing filter porosity as a design approach to increasing the mass yield per filter, while simultaneously reducing the membrane material cost per filter and achieving effective separation using multi-stage filtration.

\begin{figure}
{\scriptsize (a)}\rotatebox{0}{\includegraphics[scale=.39]{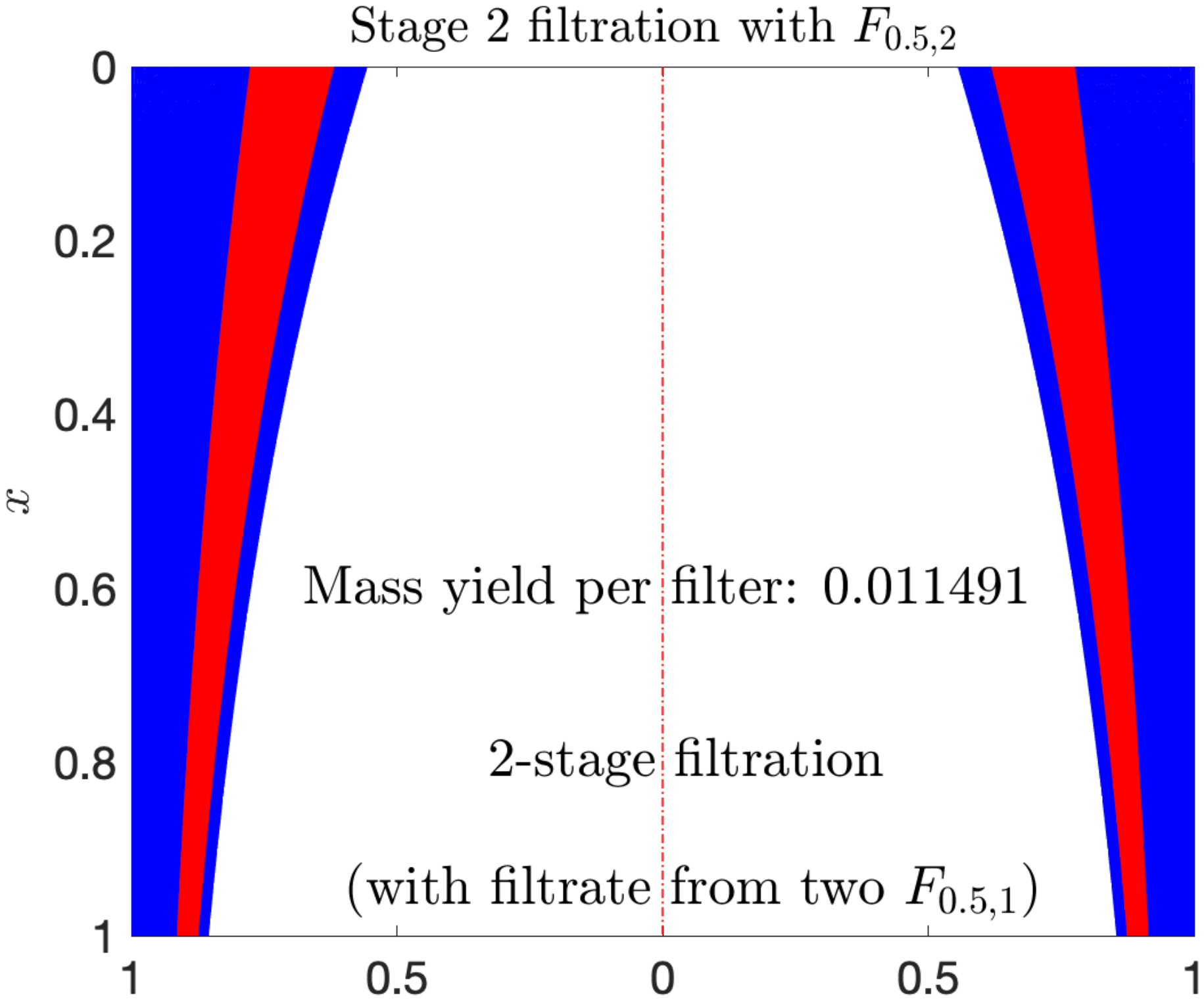}}
{\scriptsize (b)}\includegraphics[scale=.39]{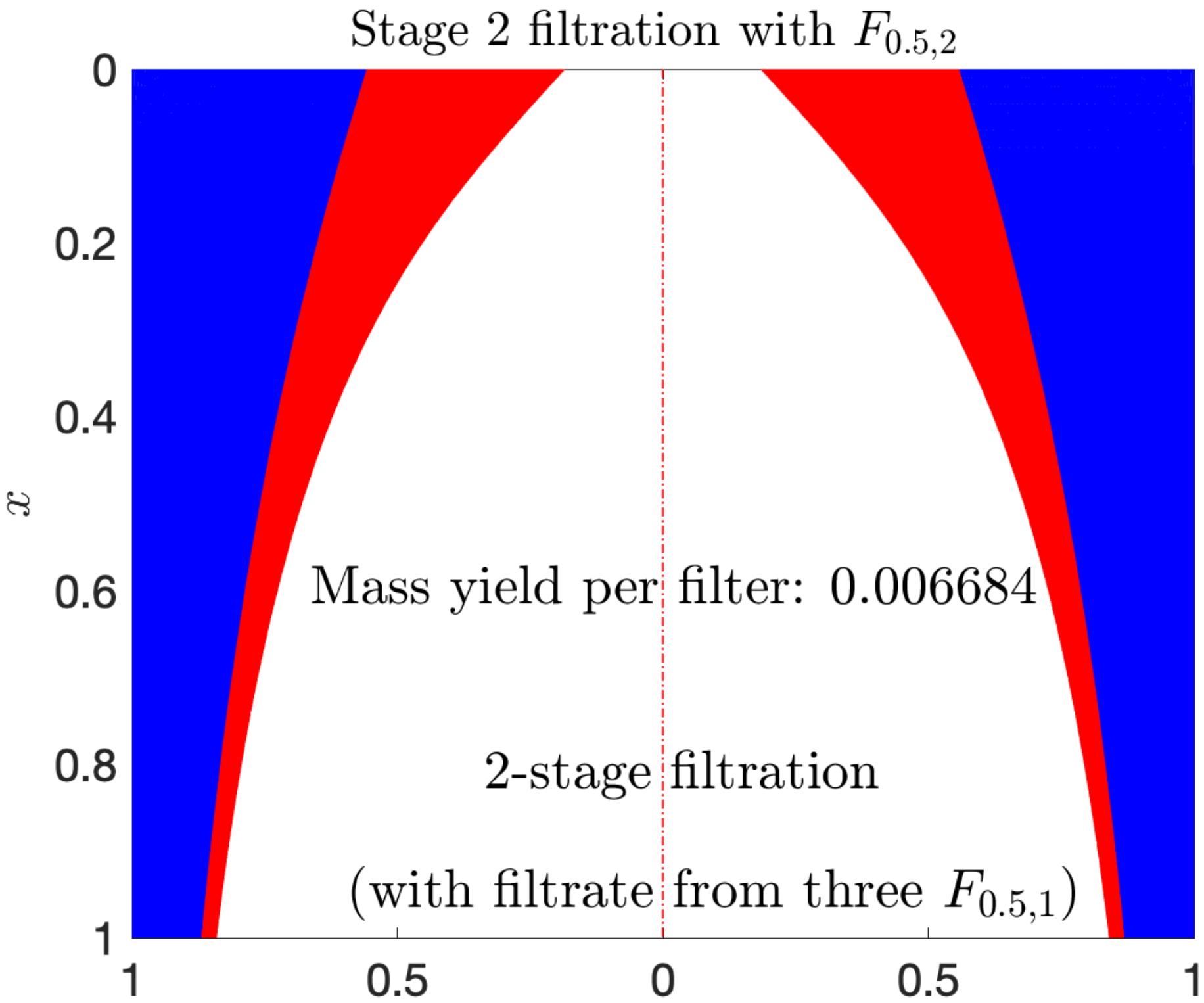}
\rotatebox{0}{\includegraphics[scale=.46]{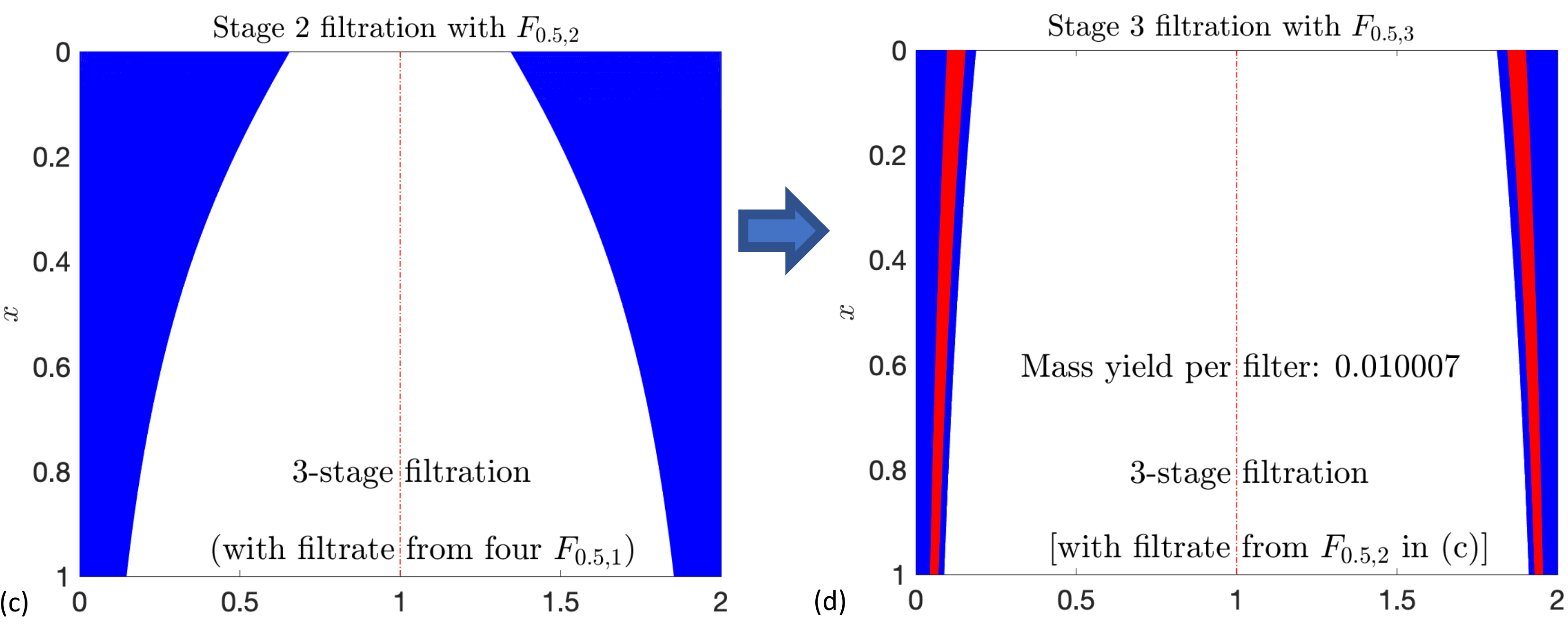}}

\caption{\footnotesize{Multi-stage filtrations: (a,b) show second stage of 2-stage filtrations; (c,d) show 2nd and 3rd stages of a 3-stage filtration. (a) fouling of $F_{0.5,2}$ by filtering filtrate collected from two $F_{0.5,1}$ filters. (b) fouling of $F_{0.5,2}$ by filtering filtrate collected from three $F_{0.5,1}$ filters. (c) and (d): 3-stage filtration with (c) fouling of $F_{0.5,2}$ by filtering filtrate collected from four $F_{0.5, 1}$ filters; (d) fouling of $F_{0.5,3}$ by filtering filtrate collected from $F_{0.5,2}$ shown in (c).}
}
\label{multi-stage_vary_stage1}
\end{figure}

From  Figs. \ref{multi-stage_vary_R} (d) and (e) it is clear that in both multi-stage filtration protocols, the secondary filters $F_{0.7, 2}$ and $F_{0.5, 2}$ are only lightly-used at termination, and could be used to process more filtrate. 
Specifically, we could use two or more first-stage filters $F_{R ,1}$ in order to create sufficient once-filtered fluid to pass through the second stage filter $F_{R ,2}$ and foul it significantly.  
We anticipate that increasing the volume of filtrate collected from stage 1 of the filtration using multiple $F_{R ,1}$ should lead to higher mass yield per filter by more fully utilizing the filtration capacity of the stage 2 filter $F_{R ,2}$.
Before investigating this idea in detail we first test it using two, three and four stage-1 filters ($l_1=2,3,4$), which lead to 2-stage, 2-stage and 3-stage filtrations, respectively. 
The results are presented in Fig. \ref{multi-stage_vary_stage1}, using the filter $F_{0.5}$ optimized for $R =0.5$ as in Fig.
\ref{multi-stage_vary_R}.

In the first test example, in stage 1 we collect filtrate by exhausting two $F_{0.5,1}$ filters ($l_1=2$; the fouling plot for each of these $F_{0.5,1}$ is identical to Fig. \ref{multi-stage_vary_R} (c) so is omitted in Fig. \ref{multi-stage_vary_stage1}) and then in stage 2, we send the combined filtrate repeatedly through an initially clean $F_{0.5,2}$. Figure \ref{multi-stage_vary_stage1}(a) shows the subsequent fouling of $F_{0.5,2}$, with alternating blue and red color indicating particle deposition and filter reuse as before. After passing the filtrate through $F_{0.5,2}$ three times, the final particle 1 removal requirement is met (so $l_2=1$ and $M=l_1+l_2=3$).

In the second example, in stage 1 we collect filtrate by exhausting three $F_{0.5, 1}$ filters ($l_1=3$; again see Fig.~\ref{multi-stage_vary_R}(c) for this stage). In stage 2 we pass the filtrate from stage 1 through an initially clean $F_{0.5,2}$. Fig.~\ref{multi-stage_vary_stage1}(b) shows the fouling of this second-stage filter $F_{0.5,2}$. It is used twice, but during the second use becomes completely fouled before all the filtrate can be filtered. Leaving aside for the moment the question of whether a second stage-2 filter should be introduced to deal with the leftover twice-filtered fluid, we check the (thrice filtered) filtrate from this second stage-2 filtration and find that it meets the final particle 1 removal requirement. In this example, $l_2=1$ and $M=1_1+l_2=4$. 

In the third example, at stage 1 we collect filtrate by exhausting four $F_{0.5,1}$ filters ($l_1=4$). This combined filtrate is then passed through a clean second-stage filter, $F_{0.5,2}$, whose fouling is shown in Fig. \ref{multi-stage_vary_stage1}(c). This $F_{0.5,2}$ filter is completely fouled after one use ($l_2=1$). Again, we defer the question of whether a second stage-2 filter would be cost-effective to deal with the remaining once-filtered feed, and check the particle-1 removal requirement of the twice-filtered feed. It is not yet satisfied, so we need a third stage of filtration with a new filter $F_{0.5, 3}$. Fig. \ref{multi-stage_vary_stage1}(d) shows the fouling of this $F_{0.5,3}$ filter, which is used three times before the final particle 1 removal requirement is satisfied ($l_3=1$). Here $M=l_1+l_2+l_3=6$.

We find that in the first example, when we collect filtrate from two $F_{0.5,1}$ filters (Fig. \ref{multi-stage_vary_stage1}(a)),
the mass yield of type 2 particles per filter is $0.012$, which is indeed greater than the value $0.0091$ obtained with the original two-stage filtration of Fig. \ref{multi-stage_vary_R}. However, with three stage-1 $F_{0.5,1}$ filters, the second example of Fig. \ref{multi-stage_vary_stage1}(b), the mass yield of type 2 particles per filter decreases to $0.0067$ (see table \ref{2t:mulltistage}), which may be explained by the fact that the second stage filter $F_{0.5,2}$ is completely fouled on its second use before all the filtrate obtained in stage 1 can be processed (the yield loss is due to the discarded filtrate). 
Similar loss of filtrate is observed in the third example, the 3-stage filtration of Figs. \ref{multi-stage_vary_stage1}(c) and (d), in which filtrate collected from four $F_{0.5,1}$ stage 1 filters was sent through a stage-2 filter $F_{0.5,2}$, which is exhausted before all of the stage-1 filtrate can be filtered a second time. 
Despite this loss, the mass yield per filter is 0.010, nearly as good as the first example of Fig.~\ref{multi-stage_vary_stage1}(a). Additional simulations of the second and third test scenarios, in which new filters were introduced to process the discarded filtrate, gave less favorable results than those presented here.
These three multi-stage filtration experiments suggest that a single stage-2 $F_{0.5,2}$ filter can process filtrate collected from three to four stage-1 $F_{0.5,1}$ filters, but no more.

The observations of Fig. \ref{multi-stage_vary_stage1}, though preliminary, indicate there may be an optimal ratio between the number of filters to use at different stages, which would utilize each filter's filtration capacity as fully as possible, and minimize the loss of filtrate at each stage, ultimately maximizing the mass yield per filter. 
We used our model to conduct such an investigation, the details of which are provided in Appendix \S \ref{optimal_ratio}. We find (by trial and error) that the mass yield per filter can be as high as 0.013 using a four-stage filtration, with the following numbers of filters per stage: $l_1=18, l_2=6, l_2=3, l_4=1$. This four-stage filtration is illustrated schematically in figure \ref{4stage}.

\begin{figure}
{\scriptsize}\rotatebox{0}{\includegraphics[scale=.5]{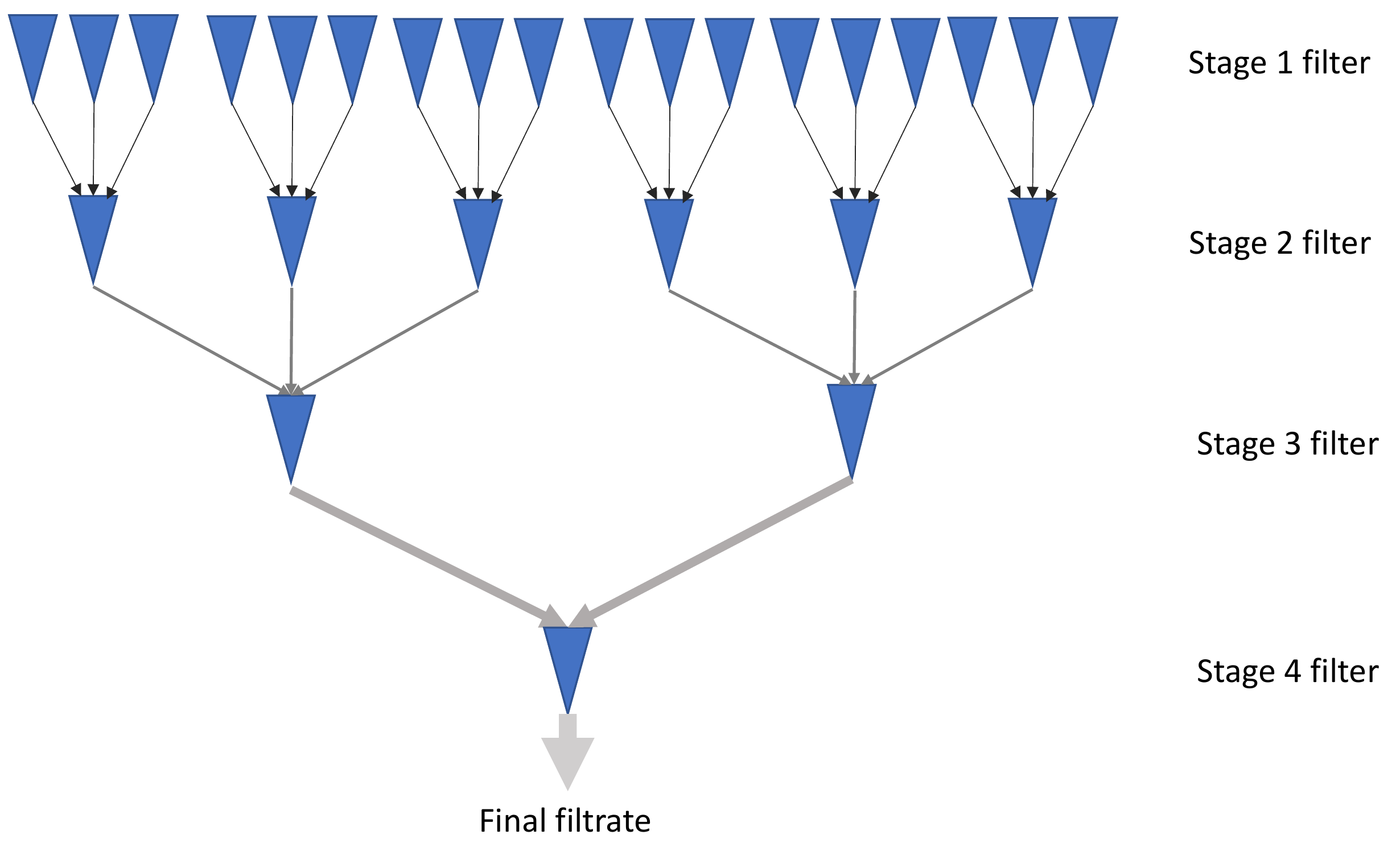}}

\caption{\footnotesize{4-stage filtration illustration, with eighteen stage 1 filters, six stage 2 filters, two stage 3 filters and one stage 4 filter.}
}
\label{4stage}
\end{figure}

\subsection{Optimization of constant flux filtration}\label{result_two_species_cflux}

In this section we briefly highlight results for the constant flux case, focusing on optimization {\bf Problem 3}, with objective function (\ref{max_Jtimes2acm_slow_cflux}). We study how the particle composition ratio in the feed ($\xi$) and the quantity of feed processed affect results. In all simulations we impose the additional constraints that the initial driving pressure $p(0,0) \le 100$ and the driving pressure at the end of the filtration should be no greater than 10 times the initial driving pressure, i.e.,  $p(0,{\rm t_f}) \le 10 p(0,0)$ (typical driving pressure at termination is about 1.5-2 times the initial pressure \cite{sun2008, blankert2006}). 
Simulations of our model at constant flux show that particle retention typically deteriorates over time, therefore, in addition to the initial particle removal requirement ($R_1(0) \ge \Rt$) we also impose that the accumulative particle removal should be greater than a fixed number at the end of the filtration ($\bar R_1({\rm t_f}) \ge \hat R_{}$, where $\hat{R}\leq \Rt$; this requirement also means that no fast optimization method is practicable for this case). The quantity of feed processed is fixed by specifying the number of time iterations $N$ (in all simulations presented here $N\le 1000$ and the time step is fixed).

Figure \ref{constant_flux_vary_xi} shows the evolution of the optimized pore profiles obtained with $\xi=0.9,0.5,0.1, \beta=0.1, \lambda_1=10, \hat R_{}=0.98, N=1000$. The figure shows: (a) driving pressure vs throughput, $(p(0,t),j)$-plot; (b) accumulative type 1 particle concentration in the filtrate vs throughput, $(c_{1\rm acm}, j)$-plot; (c) accumulative type 2 particle concentration in the filtrate vs throughput, $(c_{2\rm acm}, j)$-plot; and (d-f) show the pore profiles at the termination of the filtration for the three $\xi$-values, with blue color indicating deposited particles: (d) $\xi=0.9$; (e) $\xi=0.5$; (f) $\xi=0.1$. We observe that, as $\xi$ varies, the optimized pore profile changes significantly. For feed containing less impurity (the smallest value, $\xi=0.1$, Figure~\ref{constant_flux_vary_xi}(f)) the optimized pore profile is of $\Lambda$ shape (instead of the V shape we observed consistently in the constant pressure case) and the particle deposition is more evenly distributed over the length of the pore. Also, in contrast to the constant pressure case, we see that the particle type 1 retention capability of the filter decreases in time for all three $\xi$ values, with the most significant deterioration observed for the feed containing the highest fraction of impurity (the largest $\xi$-value, $\xi=0.9$, see Fig. \ref{constant_flux_vary_xi}(b)).  

\begin{figure}
{\scriptsize (a)}\rotatebox{0}{\includegraphics[scale=.38]{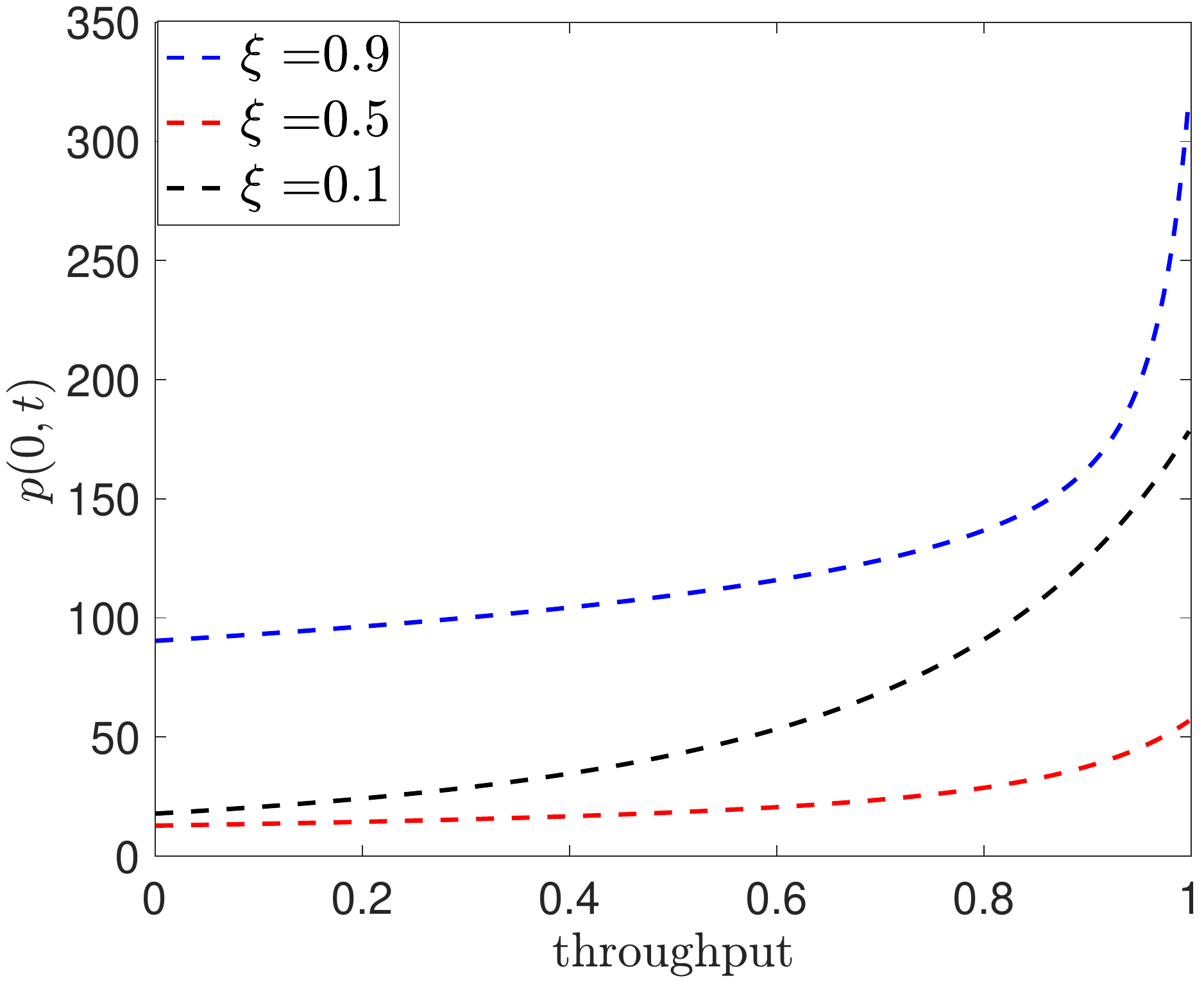}}
{\scriptsize (d)}\includegraphics[scale=.38]{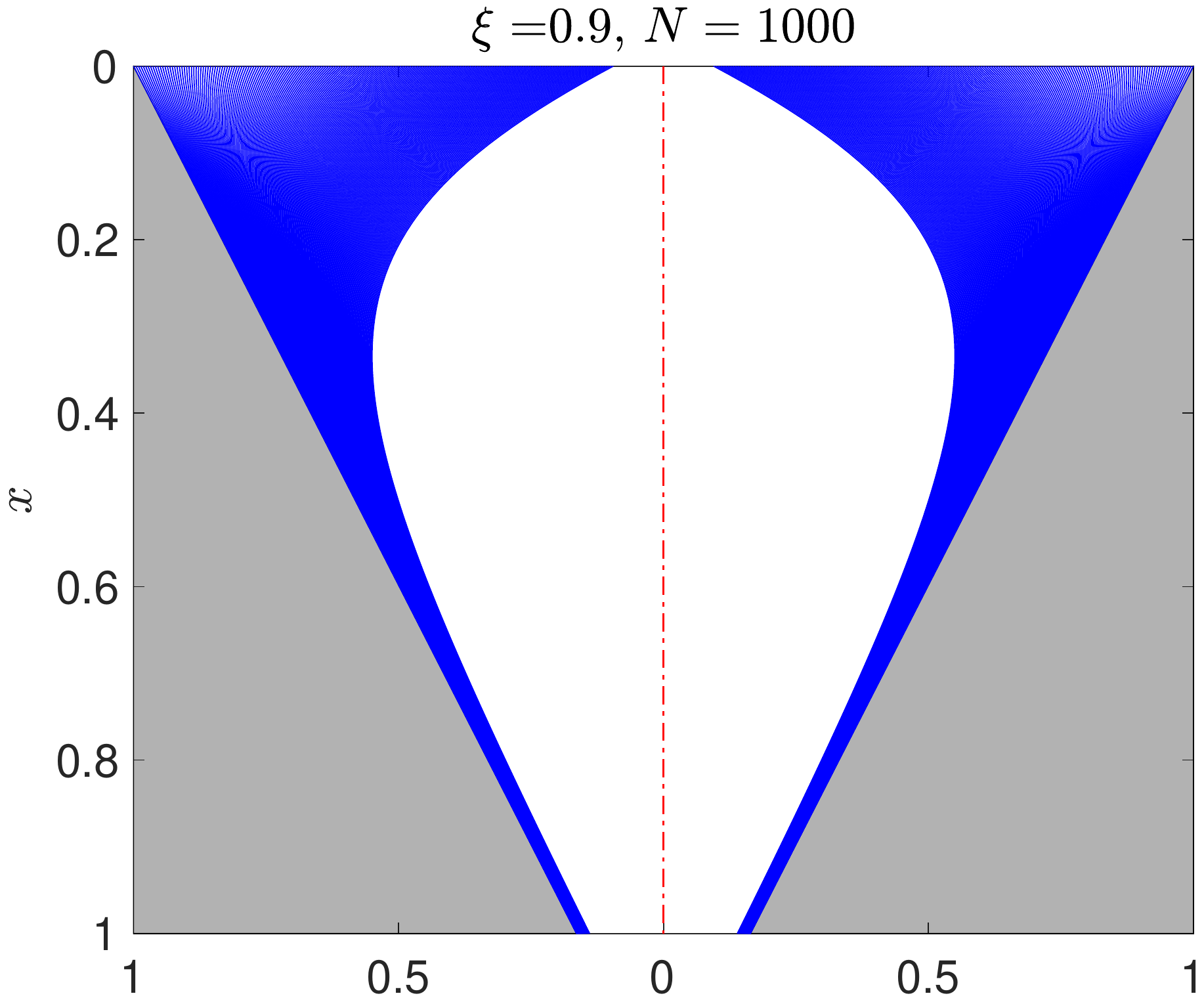}
{\scriptsize (b)}\includegraphics[scale=.37]{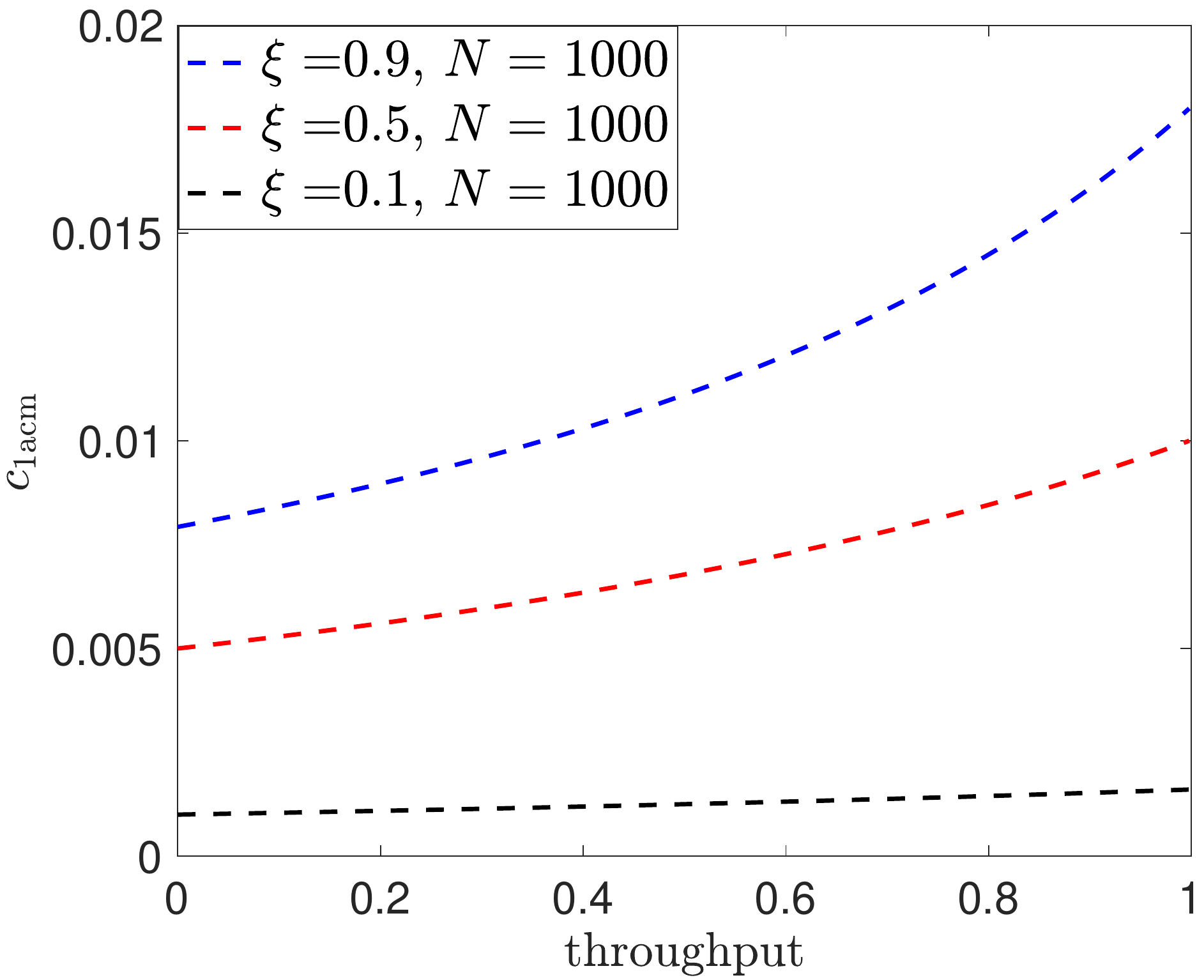}
{\scriptsize (e)}\includegraphics[scale=.38]{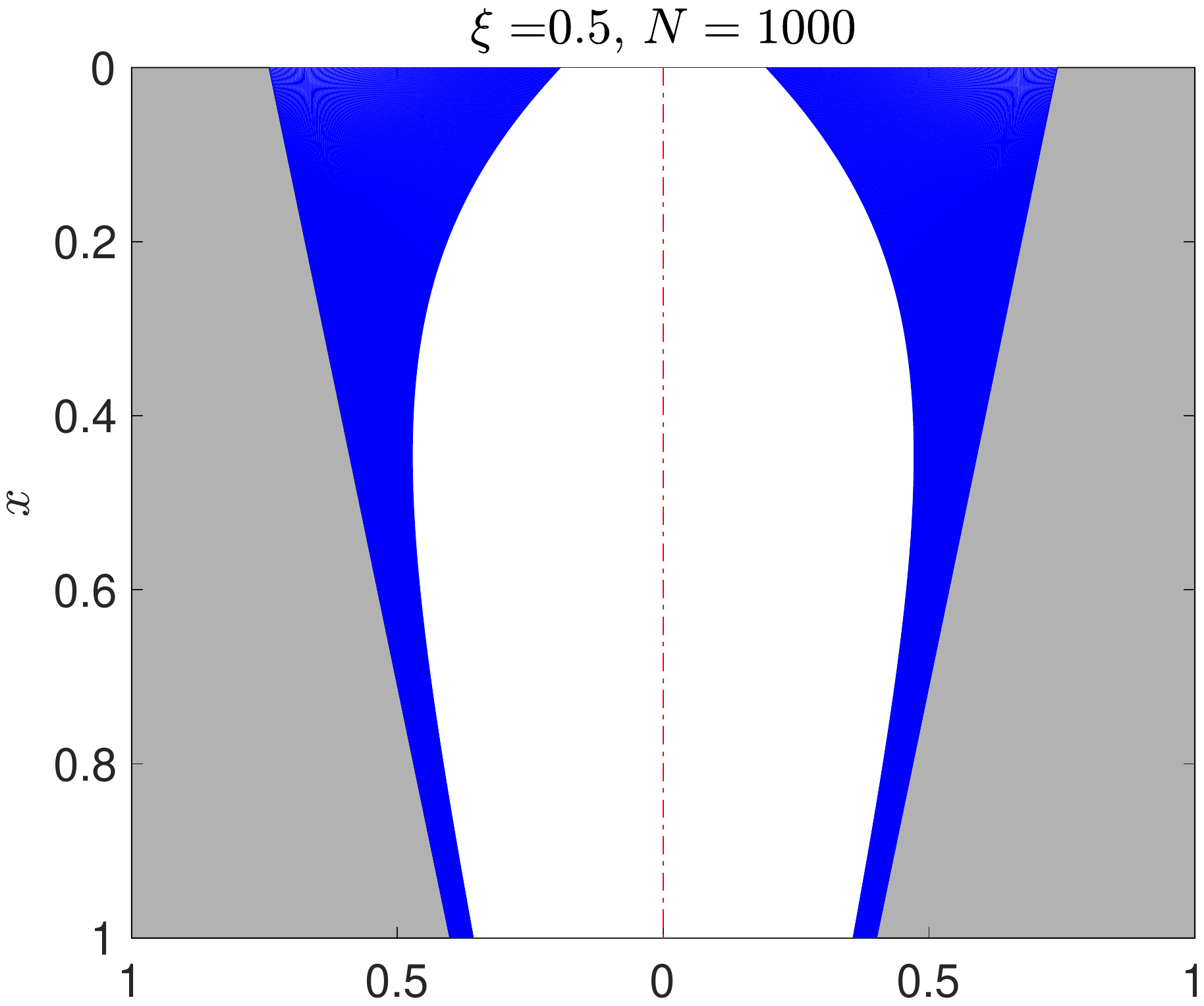}
{\scriptsize (c)}\includegraphics[scale=.37]{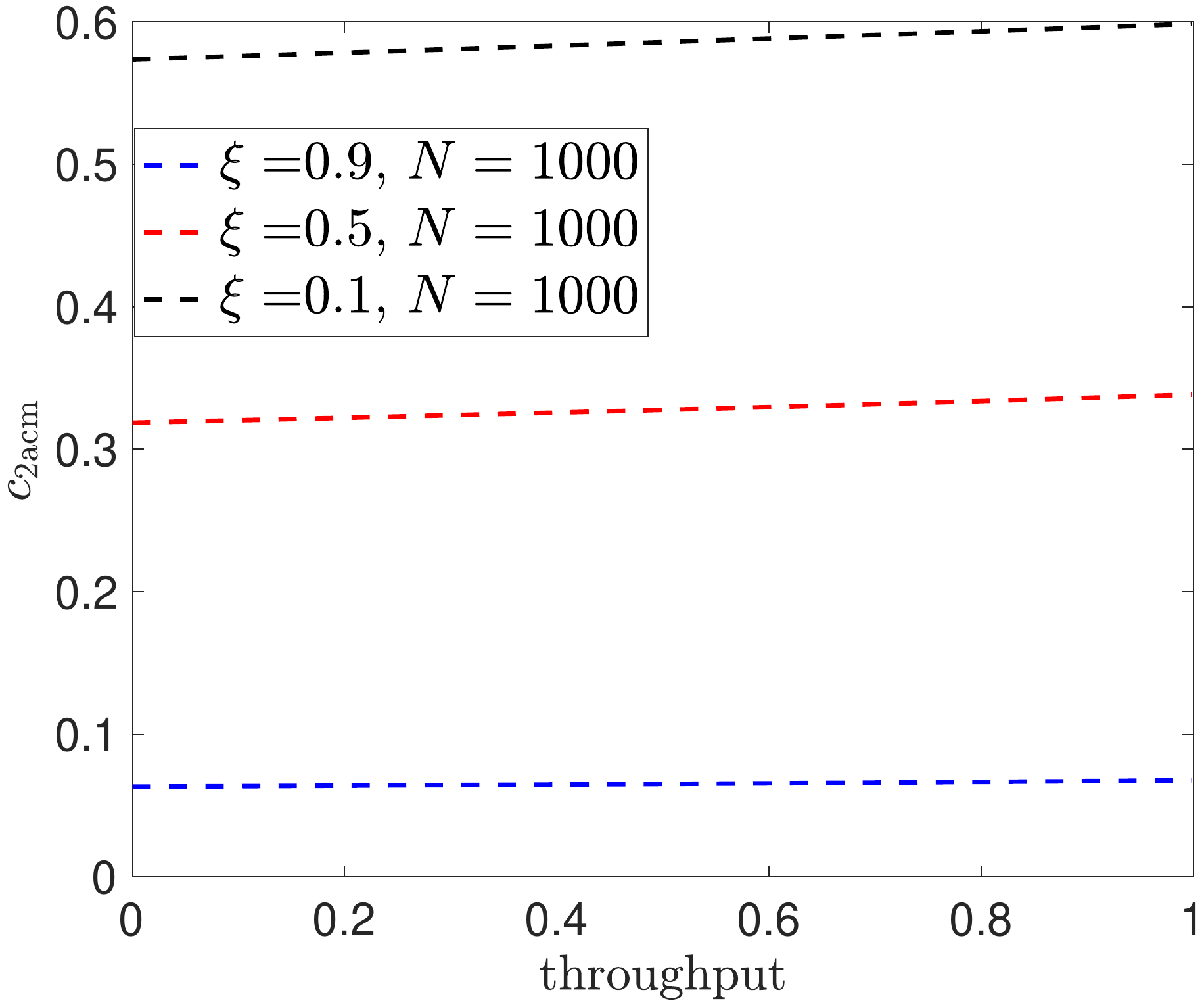}
{\scriptsize (f)}\includegraphics[scale=.38]{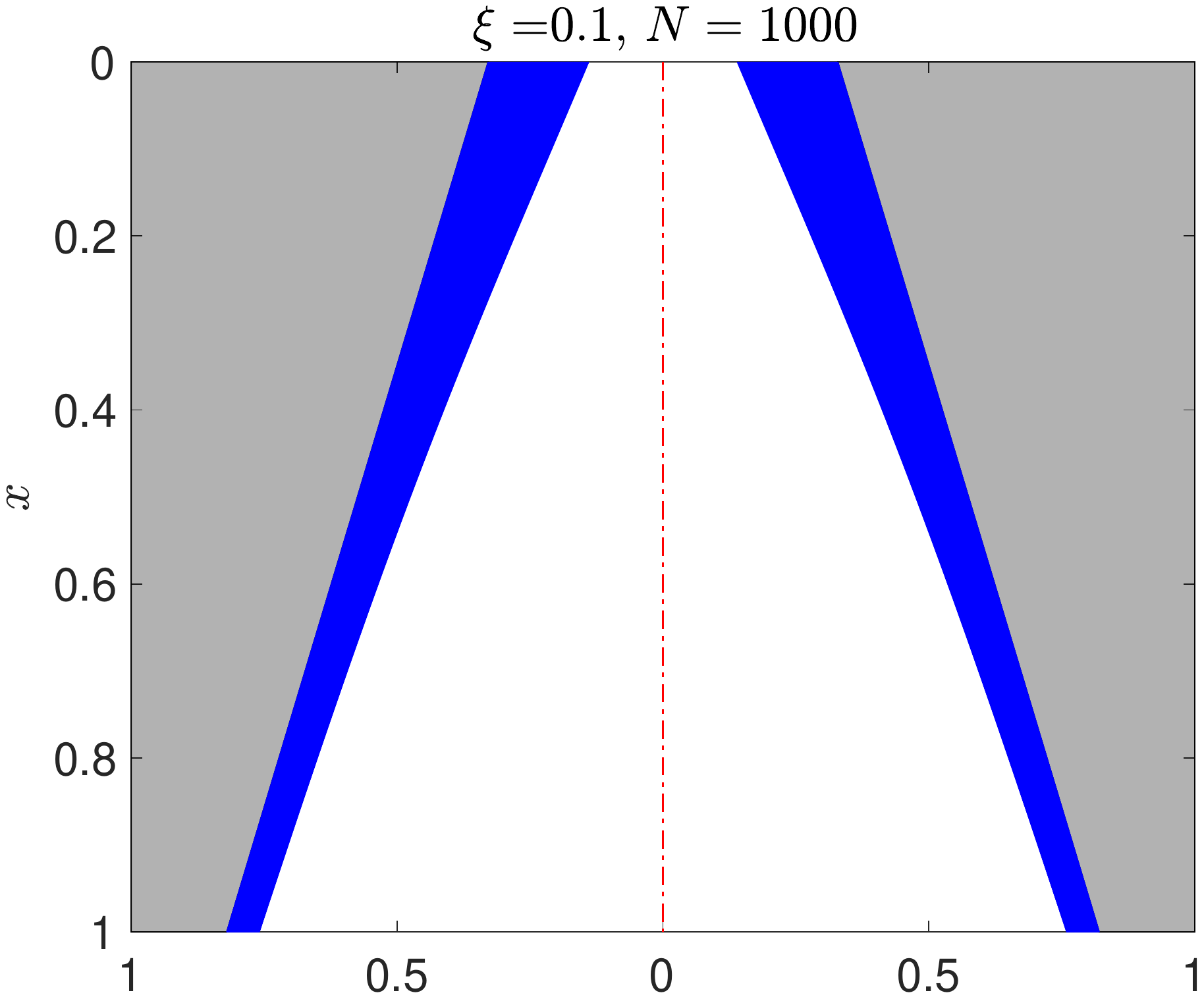}
\caption{\footnotesize{ (a-c) Evolution of optimized filters obtained for constant flux case with $\beta=0.1, \lambda_1=10, \hat R_{}=0.98, N=1000$ and $\xi=0.9,0.5,0.1$: (a) driving pressure vs throughput $(p(0,t),j)$ plot; (b) accumulative type 1 particle concentration in the filtrate vs throughput, $(c_{1\rm acm}, j)$ plot, and (c) $(c_{2\rm acm}, j)$ plot. (d-f) Pore profiles at the termination of filtration, with blue color indicating particle deposition: (d) $\xi=0.9$; (e) $\xi=0.5$; (f) $\xi=0.1$.
}}
\label{constant_flux_vary_xi}
\end{figure}

\begin{figure}
{\scriptsize (a)}\rotatebox{0}{\includegraphics[scale=.38]{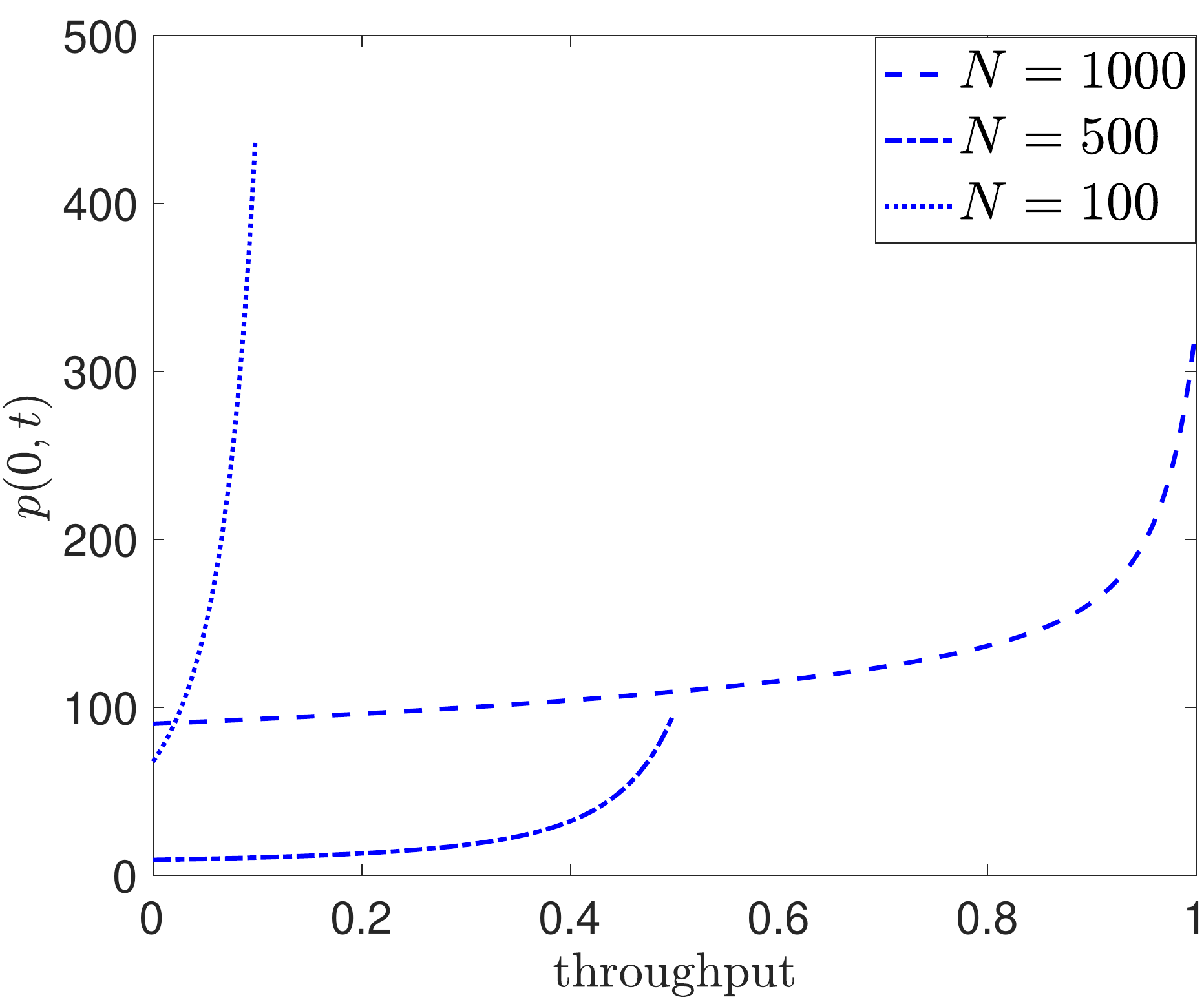}}
\hspace{5cm}
{\scriptsize (b)}\includegraphics[scale=.38]{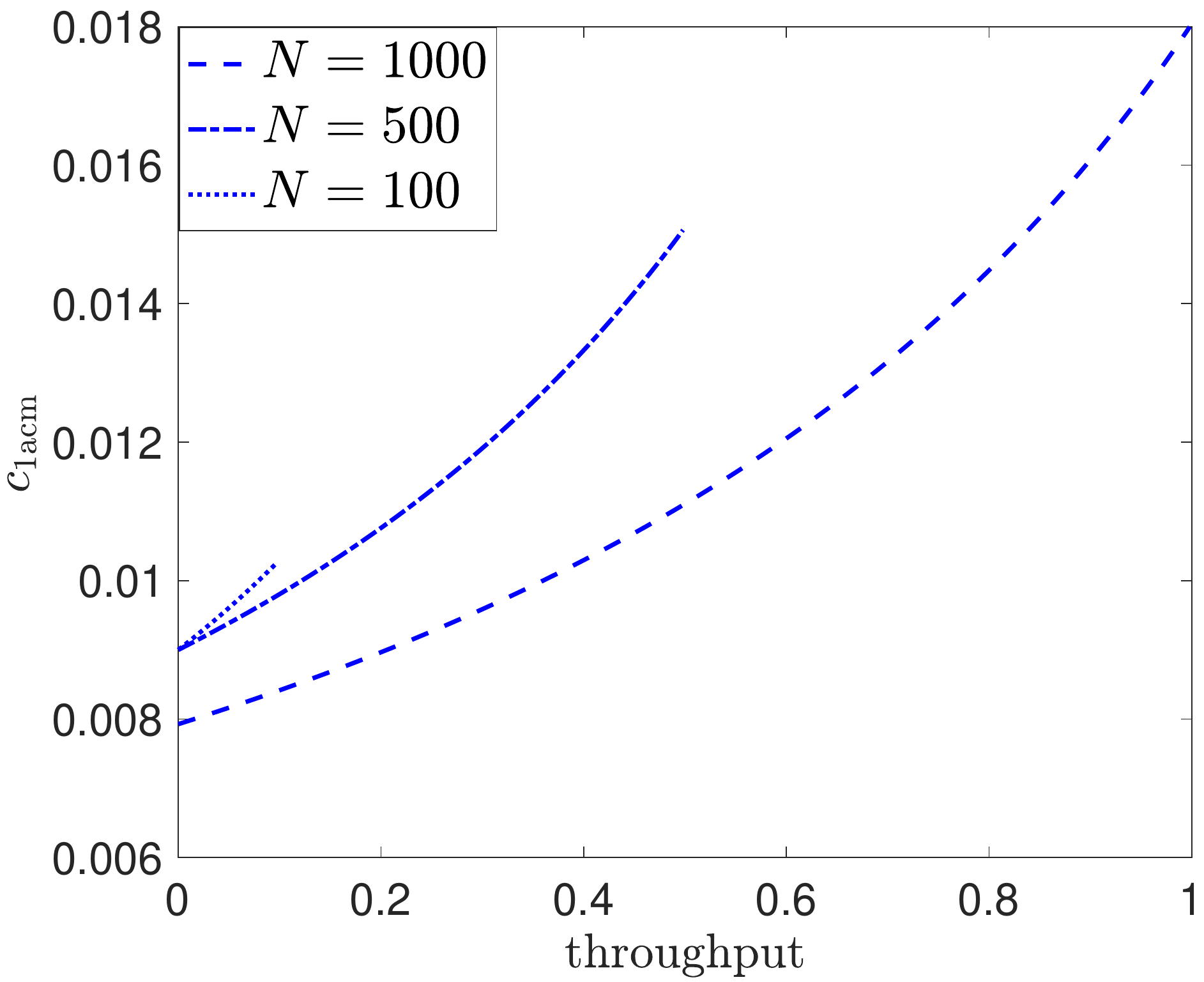}
{\scriptsize (c)}\includegraphics[scale=.38]{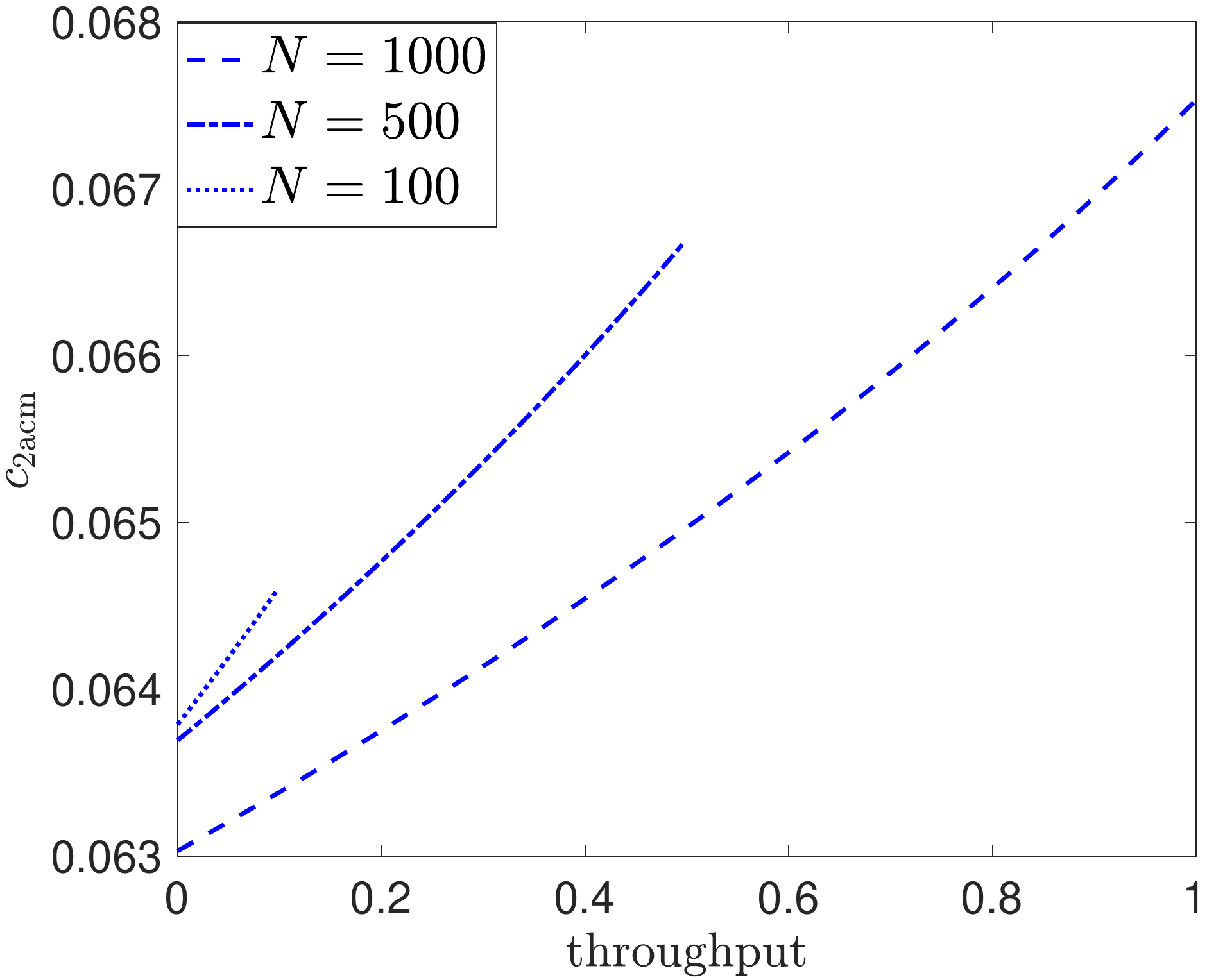}
{\scriptsize (d)}\includegraphics[scale=.39]{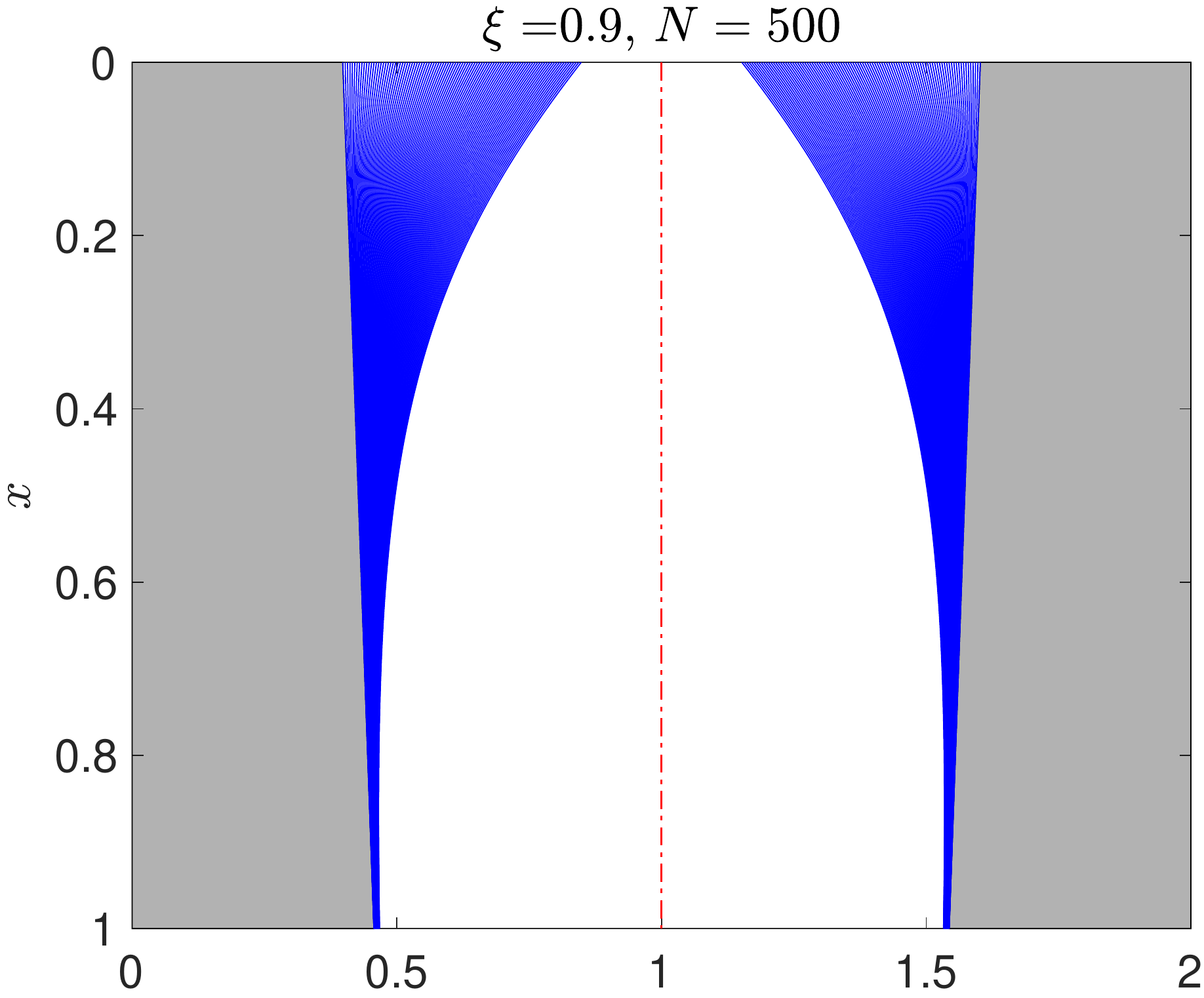}
{\scriptsize (e)}\includegraphics[scale=.39]{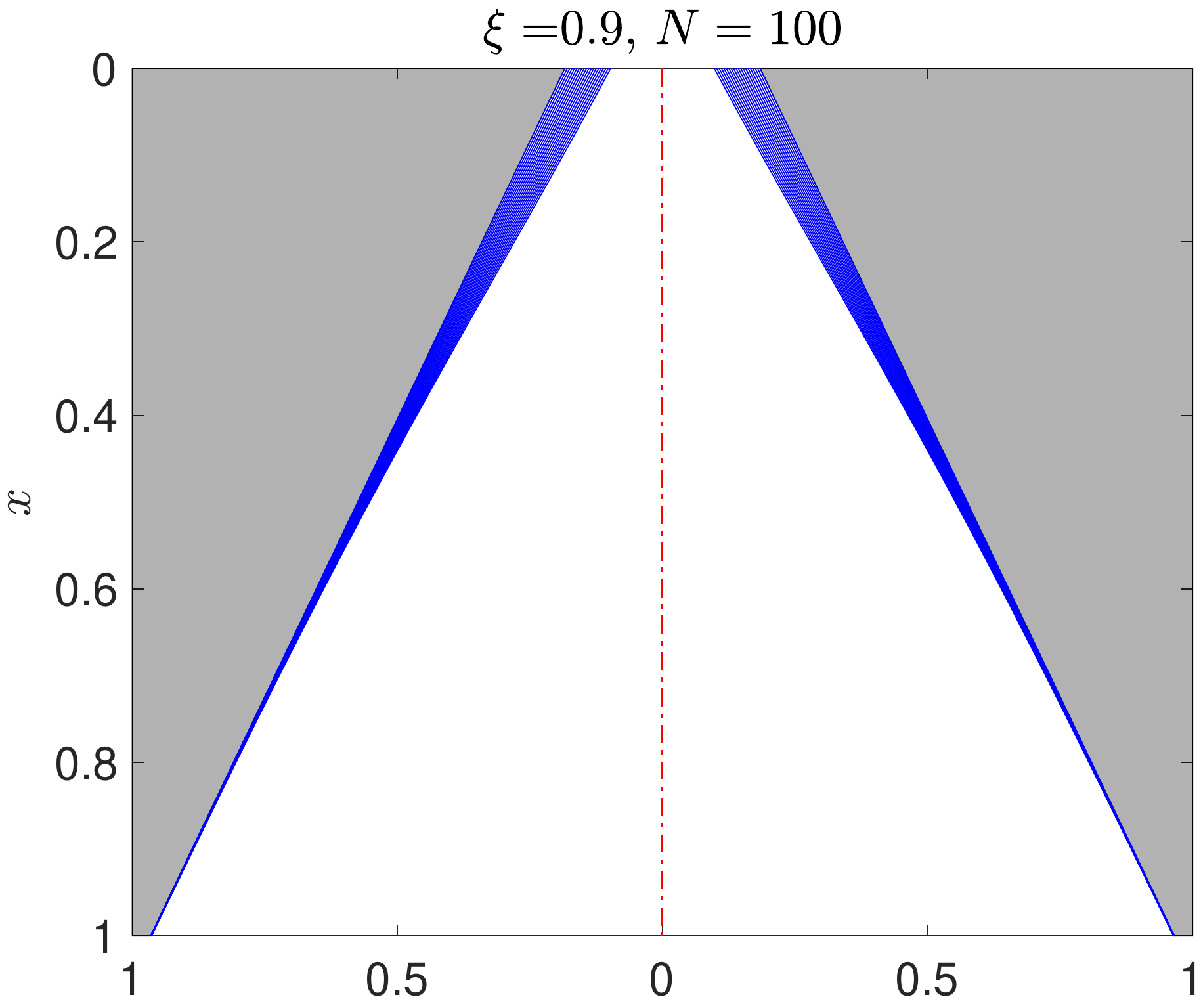}
\caption{\footnotesize{ (a-c) Evolution of optimized filters obtained for constant flux objective function (\ref{max_Jtimes2acm_slow_cflux}) with $\xi=0.9, \beta=0.1, \lambda_1=10, \hat R_{}=0.98$ and $N=1000, 500, 100$: (a) $(p(0,t),j)$ plot (b) $(c_{1\rm acm}, j)$ plot, and (c)  $(c_{2\rm acm}, j)$ plot. (d-e) Pore profiles at the termination of the filtration, with blue color indicating particle deposition: (d) $N=500$; (d) $N=100$.}}
\label{constant_flux_vary_N}
\end{figure}

Another question of interest for this constant flux case is: how does the amount of feed processed affect the optimization result? We illustrate this by considering three different values of $N$, the total number of timesteps in our simulations.  Figure \ref{constant_flux_vary_N} shows the evolution of pore profiles optimized for the constant flux objective function (\ref{max_Jtimes2acm_slow_cflux}), with $\xi=0.9, \beta=0.1, \lambda_1=10, \hat R_{}=0.98$ and $N=1000, 500, 100$. The figure shows: (a) $(p(0,t),j)$-plot; (b) $(c_{1\rm acm}, j)$-plot; (c) $(c_{2\rm acm}, j)$-plot; and (d-e) show the pore profiles at the termination of the filtration, with blue color indicating particle deposition, for (d) $N=500$ and (e) $N=100$ ($N=1000$ was shown earlier in Fig.~\ref{constant_flux_vary_xi}(d)). Collectively, Figs.~\ref{constant_flux_vary_xi}(d), \ref{constant_flux_vary_N}(d) and \ref{constant_flux_vary_N}(e) show that, as the quantity of feed decreases, the optimized pore profile changes from a V to a $\Lambda$ shape. Comparing the optimized pore profiles with Figs. \ref{constant_flux_vary_N}(a-c) we see that the $\Lambda$ shape is more prone to driving pressure increase and particle retention deterioration, as well as the more even distribution of fouling noted earlier; observations that we now explain. 

We deal first with the observation that for pores of $\Lambda$ shape, particles deposit more evenly along the pore depth compared to pores of V shape. Particle concentration is always highest at the pore entrance, which favors a high deposition rate; however, flux $u_p$ is also highest here for pores of $\Lambda$ shape, which is unfavorable for particle deposition (both observations follow from Eq. (\ref{eqn_nd_cflux_c_1})). On the other hand, at the pore exit, particle concentration is lowest (unfavorable for deposition); but flux $u_p$ is also lowest (favorable for particle deposition). Hence, for pores of $\Lambda$ shape, there is always a competition between particle concentration and flux, which leads to the observed even fouling distribution along the pore length.

We next argue heuristically that this more uniform particle deposition is responsible for the observed particle concentration increase as follows: From Eqs. (\ref{up_cflux}), (\ref{eqn_nd_cflux_c_1}) and (\ref{eqn_nd_cflux_c_2}) we obtain
\be
c_i(1,t)= \left\{
\begin{array}{c}
\xi \\
1-\xi 
\end{array}
\right\}
\exp \Big [-\frac{\lambda_i \pi}{4} \int_0^1 a(x,t)dx \Big], \qquad
 \left\{
\begin{array}{c}
i=1 \\
i=2
\end{array}
\right\},
\label{eqn_c_exp}
\ee
showing that the change in particle concentration at outlet for type $i$ particles depends on the change in the value of $\int_0^1 a(x,t)dx$. For $\Lambda$-shaped pores particle deposition is more even, thus $a(x,t)$ changes over the entire depth of the pore, with the consequence that $\int_0^1 a(x,t)dx$ changes more significantly than for V-shaped pores, where $a(x,t)$ changes significantly near the pore entrance, but on a region of small measure. The net effect for the $\Lambda$-shaped pore is the observed particle concentration increase in time in the filtrate.
The same argument may also explain the significant pressure change for pores of $\Lambda$ shape compared with pores of V shape, as the pressure change depends on the change of $\int_0^1 {a^{-4}(x,t)}{dx}$, see Eq. (\ref{eqn_nd_cflux_p}).

Collectively, these arguments suggest the following explanation for why the $\Lambda$ shape is selected for lower quantities of feed: For less feed, the filtration duration will be shorter; the significant particle concentration increase at the beginning of filtration observed with the $\Lambda$-shaped pore (see Fig. \ref{constant_flux_vary_N} (c))
is favorable for increasing the mass yield of type 2 particles, while the short filtration duration keeps the concentration increase for type 1 particles within the prescribed removal limit. 

\begin{figure}
{\scriptsize (a)}\rotatebox{0}{\includegraphics[scale=.37]{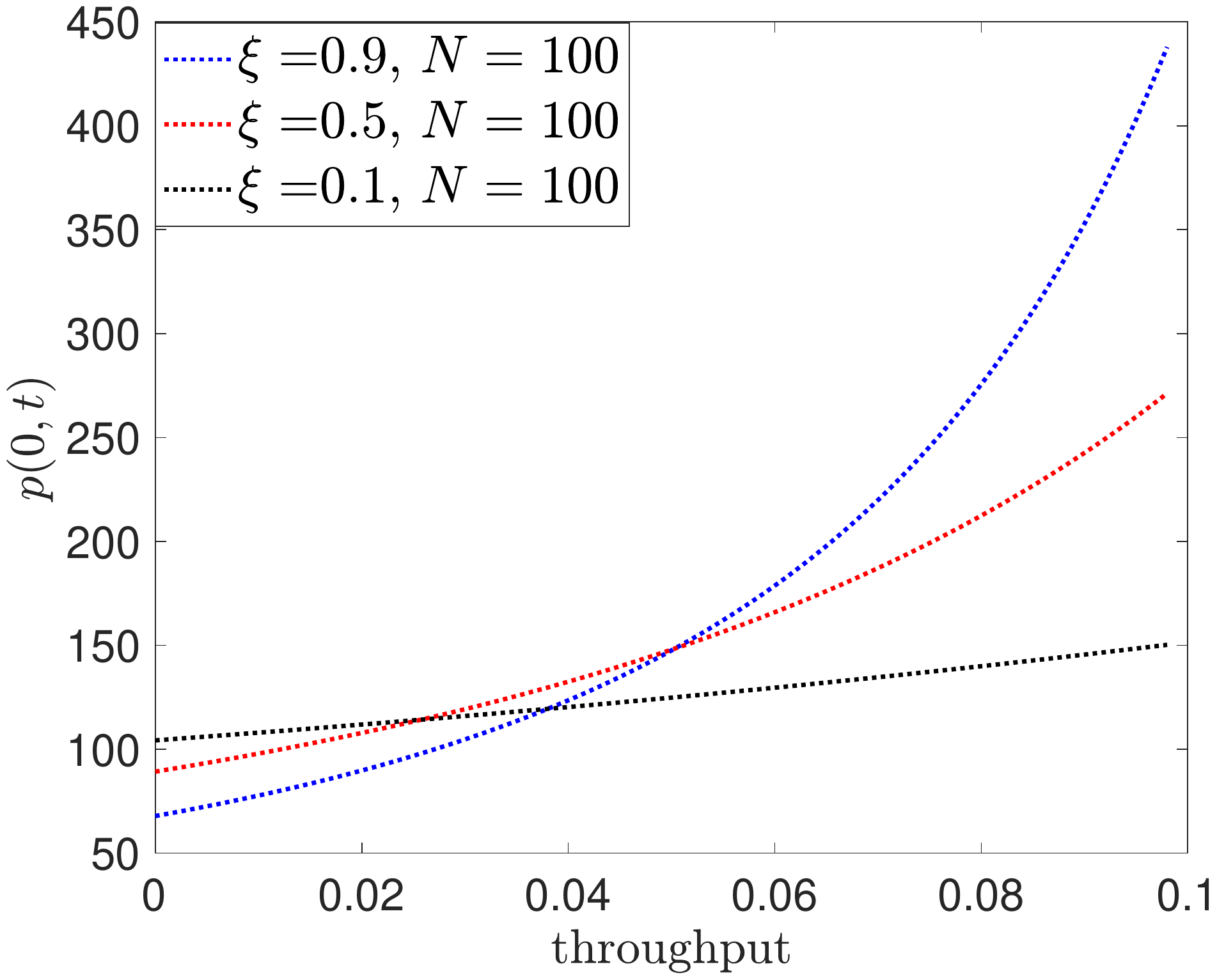}}
\hspace{5cm}
{\scriptsize (b)}\includegraphics[scale=.38]{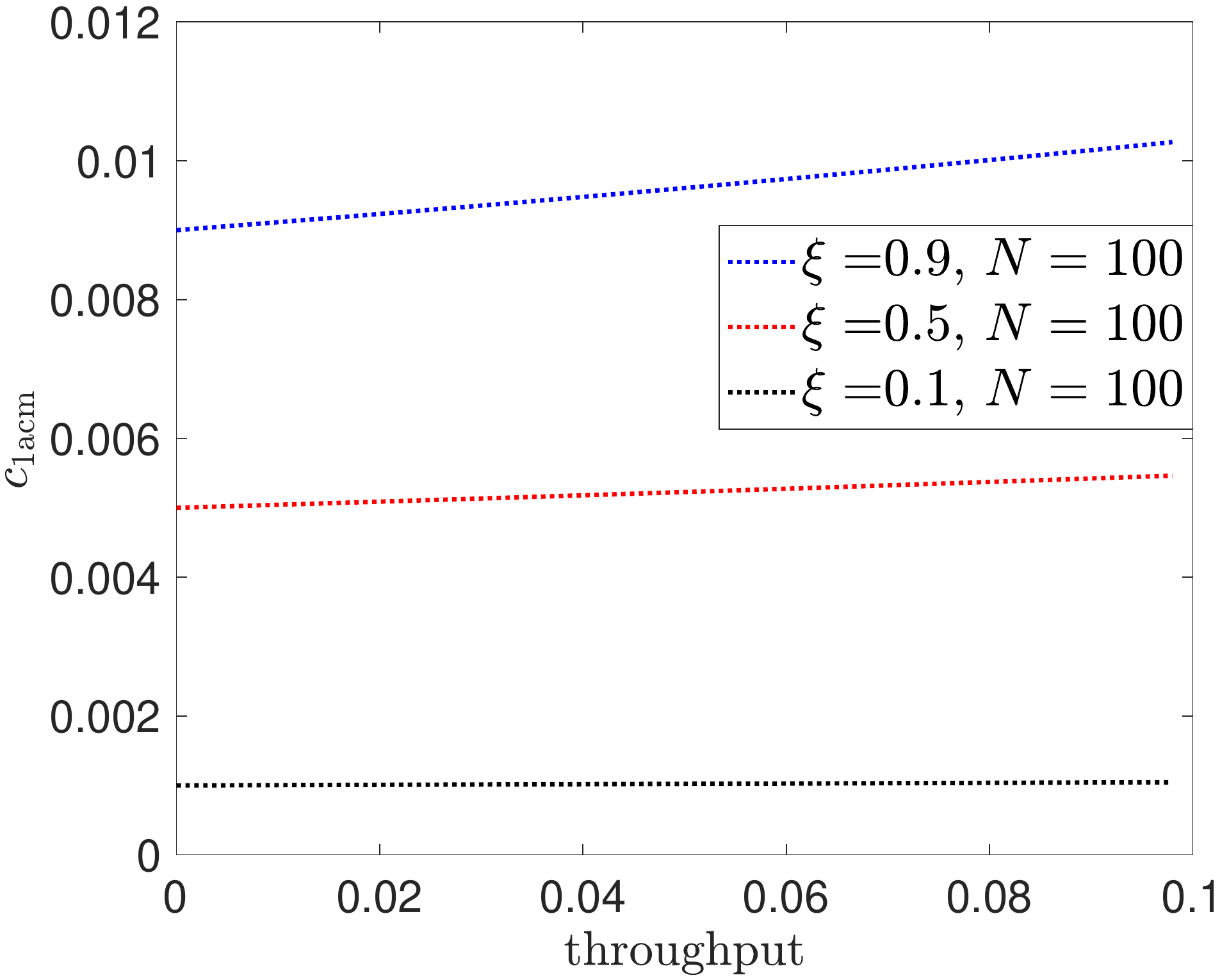}
{\scriptsize (c)}\includegraphics[scale=.38]{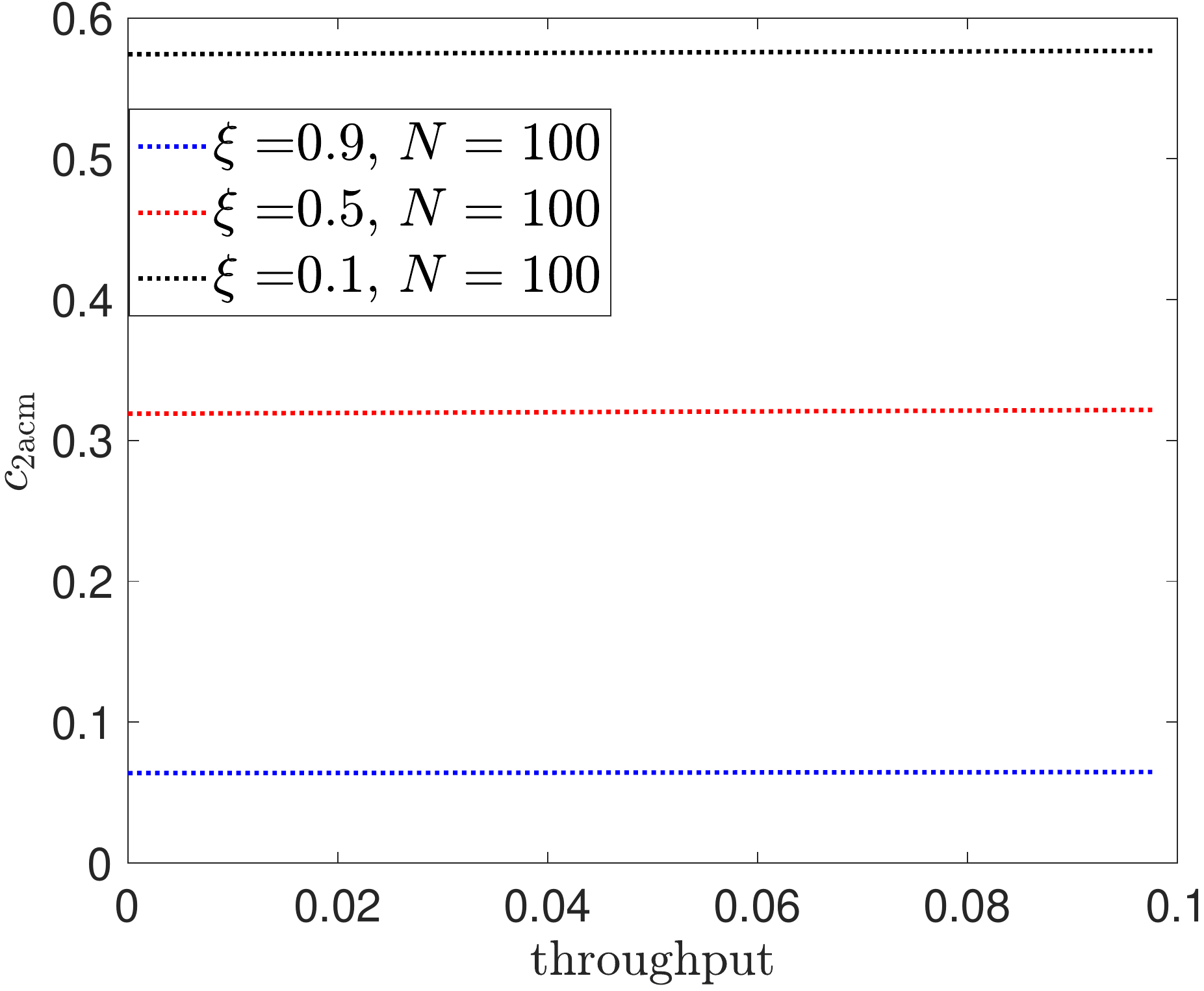}
{\scriptsize (d)}\includegraphics[scale=.38]{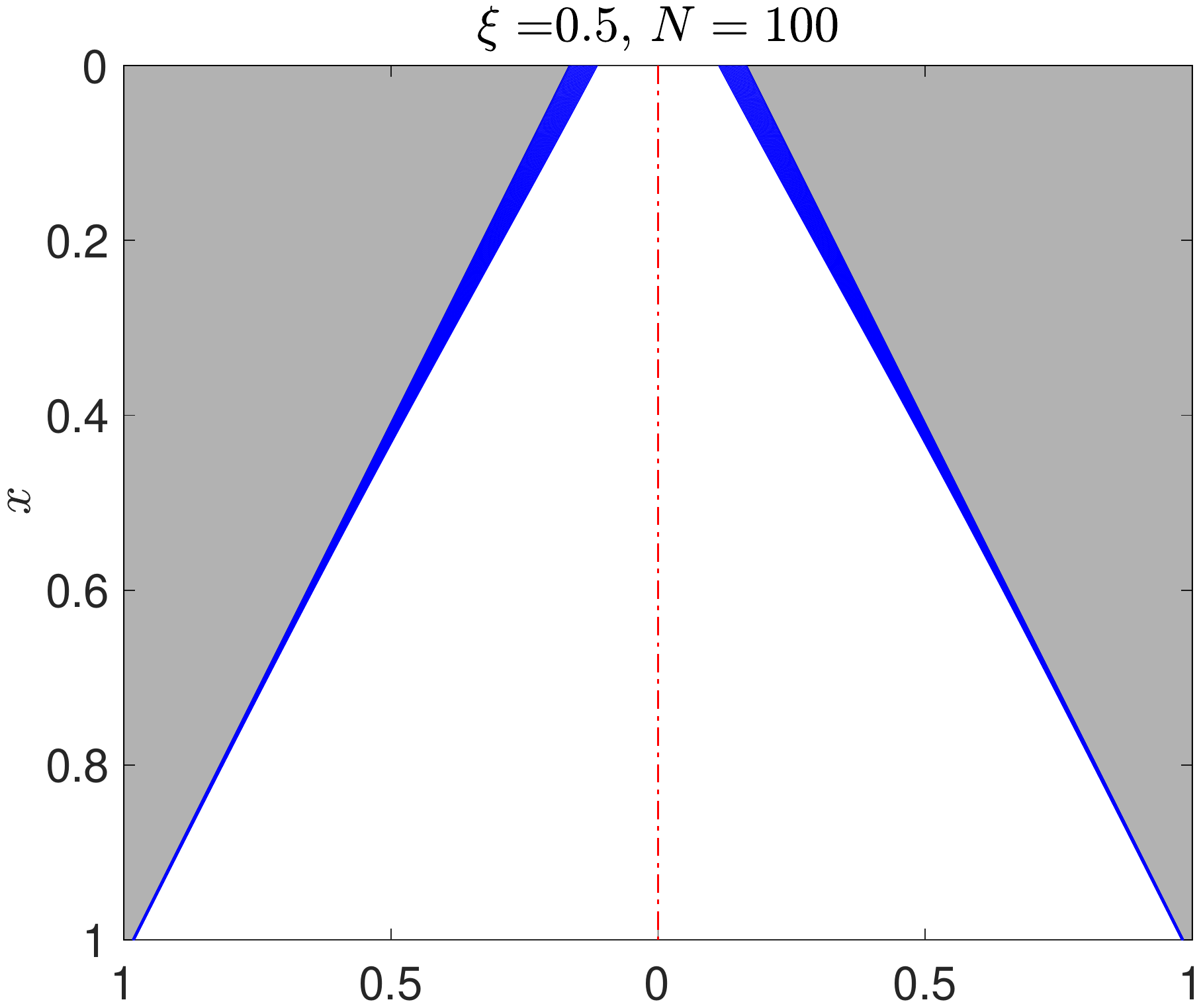}
{\scriptsize (e)}\includegraphics[scale=.38]{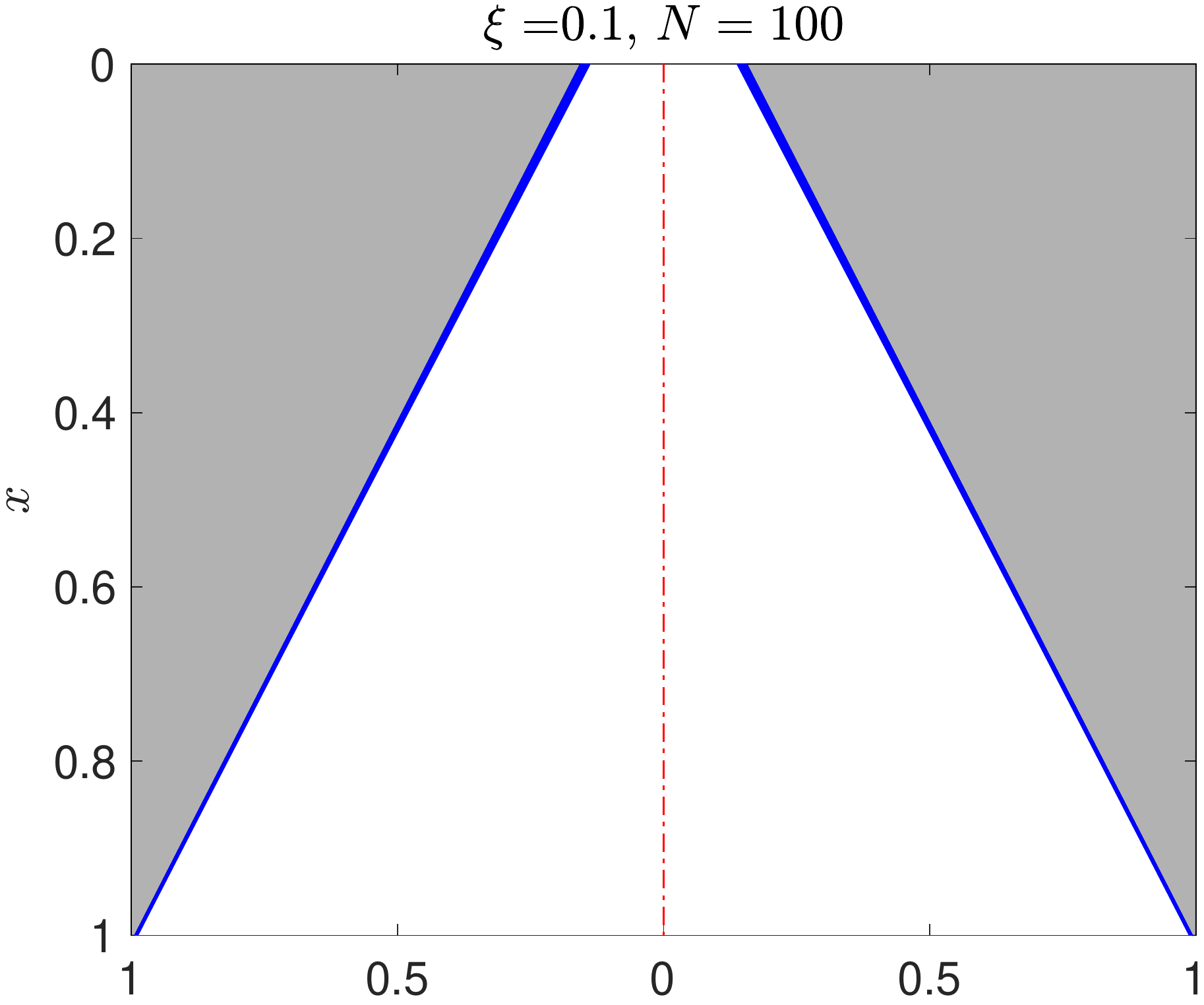}
\caption{\footnotesize{ (a-c) Evolution of optimized filters obtained for constant flux objective function (\ref{max_Jtimes2acm_slow_cflux}) with $\xi=0.9, 0.5, 0.1, \beta=0.1, \lambda_1=10, \hat R_{}=0.98$ and $N=100$: (a) $(p(0,t),j)$ plot (b) $(c_{1\rm acm}, j)$ plot, and (c)  $(c_{2\rm acm}, j)$ plot. (d-e) Pore profiles at the termination of the filtration, with blue color indicating particle deposition: (d) $\xi=0.5$; (d) $\xi=0.1$.
}}
\label{constant_flux_N=100}
\end{figure}

The results in Figure \ref{constant_flux_vary_N} raise the question of whether the optimized pore profile will take the $\Lambda$ shape more generally for sufficiently small feed quantity at constant flux regardless of particle composition ratio $\xi$. 
Figure \ref{constant_flux_N=100} shows a sequence of simulations with a small quantity of feed characterized by $N=100$, for different feed particle composition ratios $\xi=0.9, 0.5, 0.1$, with $\beta=0.1, \lambda_1=10, \hat R_{}=0.98$. 
The figure shows: (a) $(p(0,t),j)$-plot; (b) $(c_{1\rm acm}, j)$-plot; (c) $(c_{2\rm acm}, j)$-plot; while (d,e) show the pore profiles at the termination of the filtration, with blue color indicating deposited particles, for (d) $\xi=0.5$; (e) $\xi=0.1$ (the corresponding result for $\xi=0.9$ was shown in Fig. \ref{constant_flux_vary_N}(e)). 
Collectively, Fig.~\ref{constant_flux_vary_N}(e) and Figs.~\ref{constant_flux_N=100}(d) and (e) suggest that, for sufficiently small feed quantity, the optimized pore profile takes a $\Lambda$-shape regardless of feed particle composition. From Figs.~\ref{constant_flux_N=100}(b) and (c) we see the particle concentration changes are not significant for these short duration cases ($N=100$), and the particle removal requirement for type 1 particles is satisfied for all three feeds with different particle-composition ratios; however, the pressure increase is still visible, see Fig.~\ref{constant_flux_N=100}(a).


\section{Conclusions and Future Study\label{sec:conclusion}}

In this work we proposed a simplified mathematical model for filtration of feed containing multiple species of particles. Our focus in the main body of the paper was on a feed that contains just two particle species; a brief discussion of how the model extends to an arbitrary number of species is given in Appendix~\ref{results_multi_species}. For the two-species case, two important model parameters were identified and investigated to elucidate their effect on separation and optimal filter design: $\xi$, the concentration ratio of the two particle types in the feed, and $\beta =\Lambda_2\alpha_2/(\Lambda_1\alpha_1)$, the ratio of the effective particle deposition coefficients for the two particle types. A number of optimization problems for maximizing the mass yield of one particle species in the feed, while effectively removing the other, were considered, under both constant pressure and constant flux driving conditions. 
For filtration driven by a constant pressure drop, we found that the optimized pore profile is always of V-shape, which is in agreement with our earlier findings \cite{sun2020} for single-particle-species filtration (where the goal is to maximize total throughput of filtrate over the filter lifetime while removing a sufficient fraction of impurity). For filtration driven by a constant flux, the optimized pore profile may take either a V-shape or a $\Lambda$-shape depending on the particle composition ratio and the amount of feed considered for the optimization scenarios.

To increase the appeal and utility of our model for filter design applications, we proposed new objective functions (the fast optimization method) based on evaluating key quantities at the initial stage of the filtration. Due to the simpler forms of the proposed objectives, the fast method can be carried out with a relatively small number of initial search-points in design parameter space (compared with the slow method, which requires that a large number of simulations be run through to filter failure time).  The proposed fast method is approximately 100 times faster than the naive slow method. The ideas that motivated our fast method could potentially be usefully applied to other optimization problems that require evaluation of quantities at the end of the time evolution, provided those quantities exhibit some monotonicity over time. 

Observing that (based on our model predictions), effective separation in a single-stage filtration is usually achieved at the expense of short filter lifetime and inefficient filter use (most of the filter remaining only very lightly fouled), we also proposed an alternative approach for maximizing the mass yield per filter while achieving effective separation, using multi-stage filtration. 
With this approach we found that the mass yield per filter could be as much as two-and-a-half times that produced by the optimal single stage filtration, and surprisingly the purity of the final product is higher as well. In addition to the higher mass yield, the filter optimized for multi-stage filtration also requires less material to manufacture, due to its higher porosity.
Multi-stage filtration has been utilized in industry \cite{elsaid2020} and reported experimentally \cite{acheampong2014, lau2020}; however, to our best knowledge, little attention has been paid to optimizing this process from the theoretical side. 
We hope that our work will inspire further systematic studies into this promising approach.   
 
\subsection*{Acknowledgements} 

We thank Dr. Uwe Beuscher, Dr. Zhenyu He, and Dr. Vasu Venkateshwaran of W.L. Gore \& Associates for several useful discussions. 
All authors acknowledge financial support from the National Science Foundation under grant
NSF-DMS-1615719. 

\appendix

\section{Optimal ratio for multi-stage filtration\label{optimal_ratio}}

\begin{figure}
\includegraphics[scale=.83]{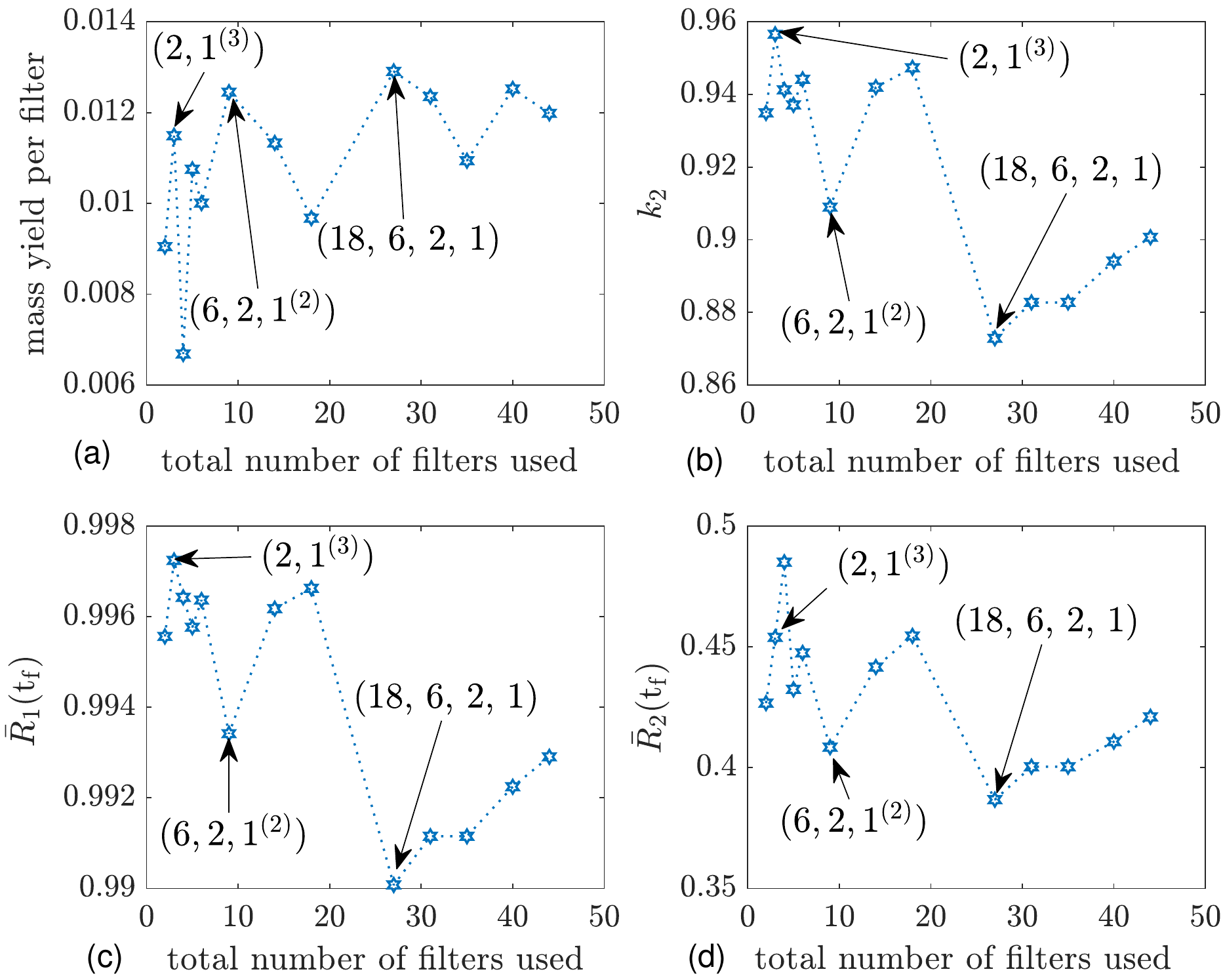}

\caption{\footnotesize{(a) Mass yield per filter, (b) purity of type 2 particles ($k_2$), (c,d) final cumulative particle removal ratios for particle types 1 and 2 ($\bar R_1(\rm t_f)$ and $\bar R_2(\rm t_f)$) are plotted against the total number of filters used in each multi-stage filtration. The local maximum mass yields per filter for 2-stage, 3-stage and 4-stage filtrations are labelled with a list of values $(l_m^{(k)})$, representing the number of filters $l_m$ used for stage $m$, and the number of times $k$ each filter is used, listed in order of increasing $m$.
}}
\label{fig9}
\end{figure}

\begin{table}
\centering
\begin{tabular}{|c|c|c|c|c|c|c|c|}
\hline
{\bf $l_1$} & {\bf $l_2$ }  & {\bf $l_3$ } & {\bf $l_4$ }  & {\bf $c_{\rm 1acm}(\rm t_f)$} & {\bf $c_{\rm 2acm}(\rm t_f)$}  & {\bf $j(\rm t_f)$}    & {\bf mass yield/filter}\\
\hline
$ 1 $ & $1^{(4)}$ &0 &0  & 0.00399 & 0.0573 & 0.316 & 0.0090\\
$ 2 $ & $1^{(3)}$ &0 &0  & 0.00247 & 0.0546 & 0.631 & 0.012\\
$ 3 $ & $1^{(2)}$ &0 &0  & 0.00321 & 0.0515 & 0.519 & 0.0067\\
$ 3$ & $1$ & $1^{(3)}$ &0 &0.00380 & 0.0568 & 0.947 & 0.011 \\
$ 4$ & $1$ & $1^{(3)}$ &0 &0.00327 & 0.0552 & 1.087 & 0.010 \\
$ 6$ & $2$ & $1^{(2)}$ &0 &0.00592 & 0.0592 & 1.894 & 0.012 \\
$ 9$ & $3$ & $1$ &$1^{(2)}$ &0.00343 & 0.0558 & 2.841 & 0.011 \\
$ 12$ & $4$ & $1$ &$1^2$ &0.00303 & 0.0545 & 3.193 & 0.0096 \\
$ 18$ & $6$ & $2$ &$1$   &0.00893 &  0.0589 & 5.683 & {\color{red}0.013} \\
$ 21$ & $7$ & $2$ &$1$ &0.00796 &0.0600 & 6.386 & 0.012 \\
$ 24$ & $8$ & $2$ &$1$ &0.00796 & 0.0600 & 6.386 & 0.011 \\
$ 27$ & $9$ & $3$ &$1$ &0.00698 &0.0589& 8.501 & 0.012 \\
$ 30$ & $10$ & $3$ &$1$ &0.00638 &  0.0579 & 9.108 & 0.011 \\
\hline
\end{tabular}
\caption{\footnotesize{Comparisons of multi-stage filtrations (up to four stages are considered) with differing ratios of the number ($l_m$) of filters $F_{0.5,m}$ used at stage $m$. The first four columns list values $l_m^{(k)}$, with superscripts $(k)$ indicating that each filter is used $k$ times. The remaining columns show final cumulative particle concentration for type 1 and type 2 particles, $c_{\rm 1acm}(\rm t_f)$ and $c_{\rm 2acm}(\rm t_f)$, total throughput $j(\rm t_f)$ and compound 2 mass yield per filter. The global maximum mass yield per filter is highlighted in red font. }}
\label{2t:mulltistage}
\end{table}

The observations of Fig. \ref{multi-stage_vary_stage1} indicate there may be an optimal ratio between the number of filters to use at different stages of a multi-stage filtration, which would utilize each filter's filtration capacity as fully as possible, and minimize the loss of filtrate at each stage, ultimately maximizing the mass yield per filter. 
We used our model to conduct such an investigation, and compiled our findings in table \ref{2t:mulltistage}, which is also presented graphically in Figure~\ref{fig9}. 
At all filtration stages filters $F_{0.5}$, optimized to maximize the mass yield of type 2 particles while meeting a particle type 1 removal threshold $R=0.5$, are used.  

Figure \ref{fig9} shows, for each multi-stage filtration considered, the total mass yield per filter (Fig.~\ref{fig9}(a)); the final purity $k_2$ of the filtrate  (Fig.~\ref{fig9}(b)); and the final cumulative particle removal ratios $\bar{R}_1 (t_{\rm f}), \bar{R}_2 (t_{\rm f})$ for the two particle types. We find that the maximal mass yield per filter is 0.013; the corresponding (four-stage) filtration is apparent as the global maximum of the  mass yield per filter in Fig.~\ref{fig9}(a). This maximum yield is also indicated in red font in table \ref{2t:mulltistage}, and is almost two and half times the yield obtained with a single-stage filtration optimized to maximize yield while immediately satisfying the purity constraint. 

This four-stage filtration is illustrated schematically in Fig. \ref{4stage}.  We note that the higher mass yield  per filter appears to be achieved at the expense of lowest purity $k_2=0.873$ among other multi-stage filtrations considered, with final cumulative particle 1 removal ratio $\bar R_1(\rm t_f)$ just above 0.99 (see Fig.~\ref{fig9}(b,c)), which makes sense as optimizers are generally found at the boundary of the feasible search space where one or more constraints are tight. However, the 2-stage local maximum simultaneously achieves high mass yield per filter together with the highest purity, $k_2=0.997$, which indicates that this two-stage filtration may be useful to achieve high mass yield without sacrificing the purity of the final product.

\section{Multiple species\label{results_multi_species}}

In this section, we present some sample  results for feed containing more than 2 species of particles. 
We non-dimensionalize our model (\ref{eqn_darcy})--(\ref{eqn_fouling}) using the same scalings as \S \ref{Nondimensionalization_cpressure} for most quantities, with the following variations: 
\be
c_i=\frac{C_i}{\sum_i C_{0i}} \label{eqn_scale_cpressure_x_c_m} ~\text{and}~
t=\frac{T}{T_0}, ~\text{with}~ T_0=\frac{W}{\Lambda_1 \alpha_1 \sum_i C_{0i}}. \label{eqn_scale_cpressure_a_t_m}
\ee
Eqs. (\ref{darcy_box})-(\ref{bc_box})  remain unchanged, and Eqs. (\ref{eqn_deposition})-(\ref{eqn_fouling}) take the form
\be
u_{\rm p} \frac{\pa c_i}{\pa x}=-\lambda_i \frac{c_i}{a}, \quad c_i(0,t)=\xi_i  ~\text{with}~ \sum_i \xi_i =1, \label{eqn_nd_cpressure_c_i_m}\\
\frac{\pa a}{\pa t}=-\sum_i \beta_i c_i, ~a(x,0)=a_0(x), \label{eqn_nd_cpressure_fouling_m}
\ee
where $\lambda_i={32\Lambda_i D^2 \mu }/({\pi W^3 P_0})$ is the deposition coefficient for particle type $i$, $\xi_i={C_{0i}}/\sum_i C_{0i}$ is the concentration ratio of type $i$ particles, $\beta_i=\Lambda_i \alpha_i / (\Lambda_1 \alpha_1)$ is the effective particle deposition coefficient for particle type $i$ (relative to particle type 1), and $a_0(x)$ is the pore profile at initial time $t=0$. 
To illustrate the model we consider an optimization problem similar to {\bf Problem 2}. 
For definiteness, we consider a feed containing three species of particles, with type 1 and type 3 the particles to be removed and type 2 the particles to be recovered from the feed. The goal is again to maximize the mass yield for type 2 particles ($c_{2 \rm acm}(t_{\rm f}) j(t_{\rm f})$) while achieving effective separation.
Similar to our earlier approach in \S \ref{result_multi-stage}, we define {\it effective separation} based on the final cumulative particle removal ratios as  ${\bar R}_1 (t_{\rm f}) \ge 0.99$, ${\bar R}_2 (t_{\rm f}) \le 0.5$, and ${\bar R}_3 (t_{\rm f}) \ge 0.9$. 

In Figure \ref{apfig1} we present optimization results for three feeds of different particle composition ratios (indicated by the different $\xi_i$ values for each curve) with $\beta_1=1$, $\beta_2=0.1, \beta_3=0.5$, $\lambda_1=1$: (a) flux vs throughput, ($u, j$) plot; (b) accumulative type 1 particle concentration vs throughput, ($c_{1\rm acm}, j$) plot; (c) accumulative type 2 particle concentration vs throughput, ($c_{2\rm acm}, j$) plot; (d) accumulative type 3 particle concentration vs throughput, ($c_{3\rm acm}, j$) plot. 
The optimization was carried out using the slow method (similar to the slow method outlined for {\bf Problem 2} in \S \ref{opt_method_slow}), and we find that the optimized pore profile takes a V-shape, as observed in our results for {\bf Problem 2} in \S \ref{result_two_species_cpressure}. The different feed composition does not significantly change the optimized pore profile. 
For higher concentrations of the heaviest-fouling particle (type 1 in this case) the pore closes faster with less total throughput observed in Figure \ref{apfig1} (a). 
We can see from Figures \ref{apfig1} (b) and (c) that the initial particle removal requirements are satisfied, though not sharp, in all three cases.
We note in Figure \ref{apfig1} (d) that the particle removal constraint for type 3 particles is tight at the start of filtration, which indicates that this constraint is the most demanding. 
This makes sense as we are requiring a relatively high removal ratio ($90\%$), with a much lower particle deposition coefficient compared to particle type 1 ($\lambda_3=\beta_3\lambda_1$ and $\beta_3=0.5$; 
in this simulation $\alpha_i$ is the same for each type of particles). In all three cases shown, effective separation is achieved by our definition, details are provided in Table \ref{tb:3species}.
Similar to the two species problem, the slow method for optimization takes about 40 minutes with 10000 searching points. Fast methods based on similar heuristics to those discussed in \S \ref{opt_method_fast} were explored but found to give unreliable results for three particle species; further investigation is needed to speed up the optimization.  
\begin{figure}
{\scriptsize (a)}\rotatebox{0}{\includegraphics[scale=.36]{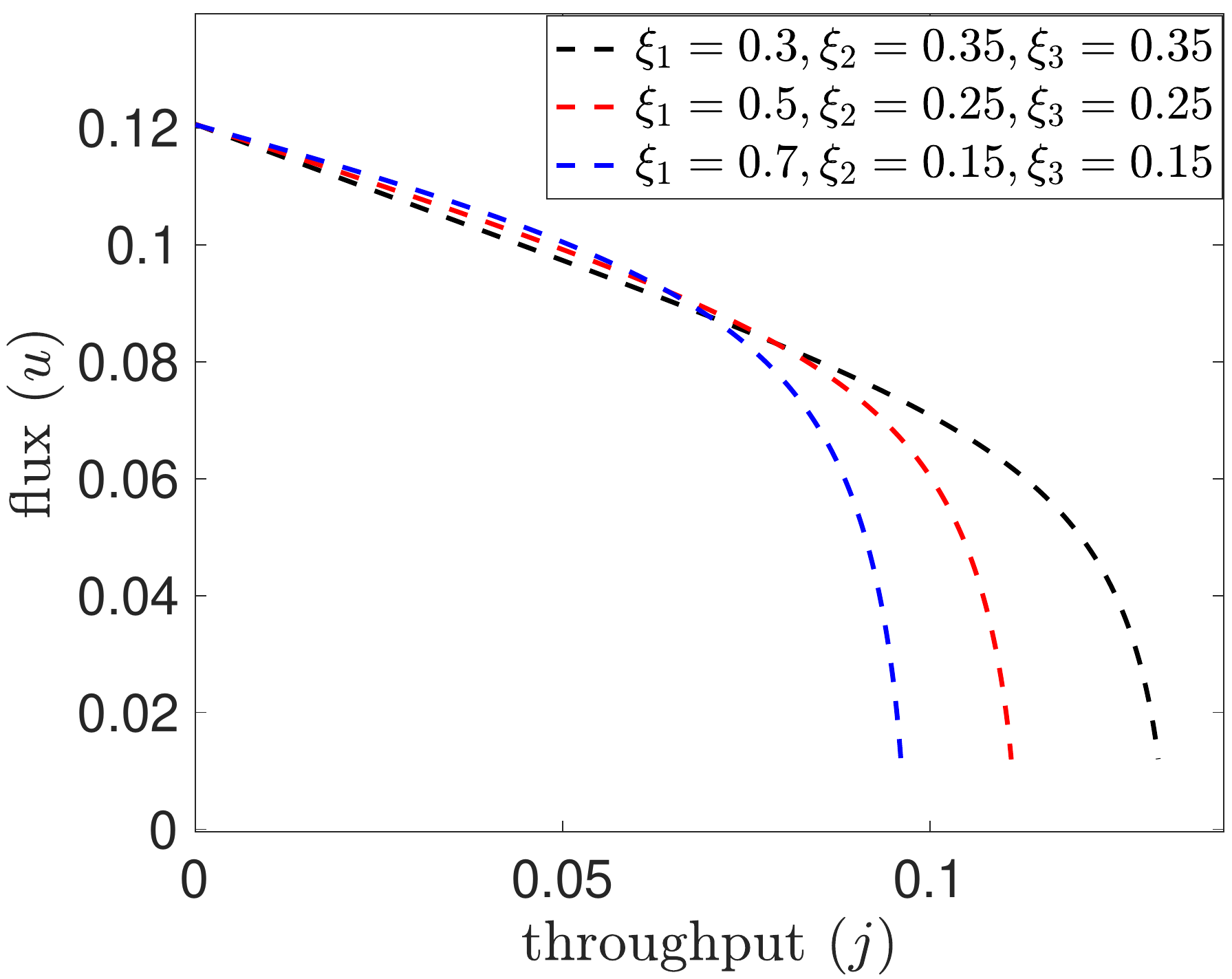}}
{\scriptsize (b)}\includegraphics[scale=.36]{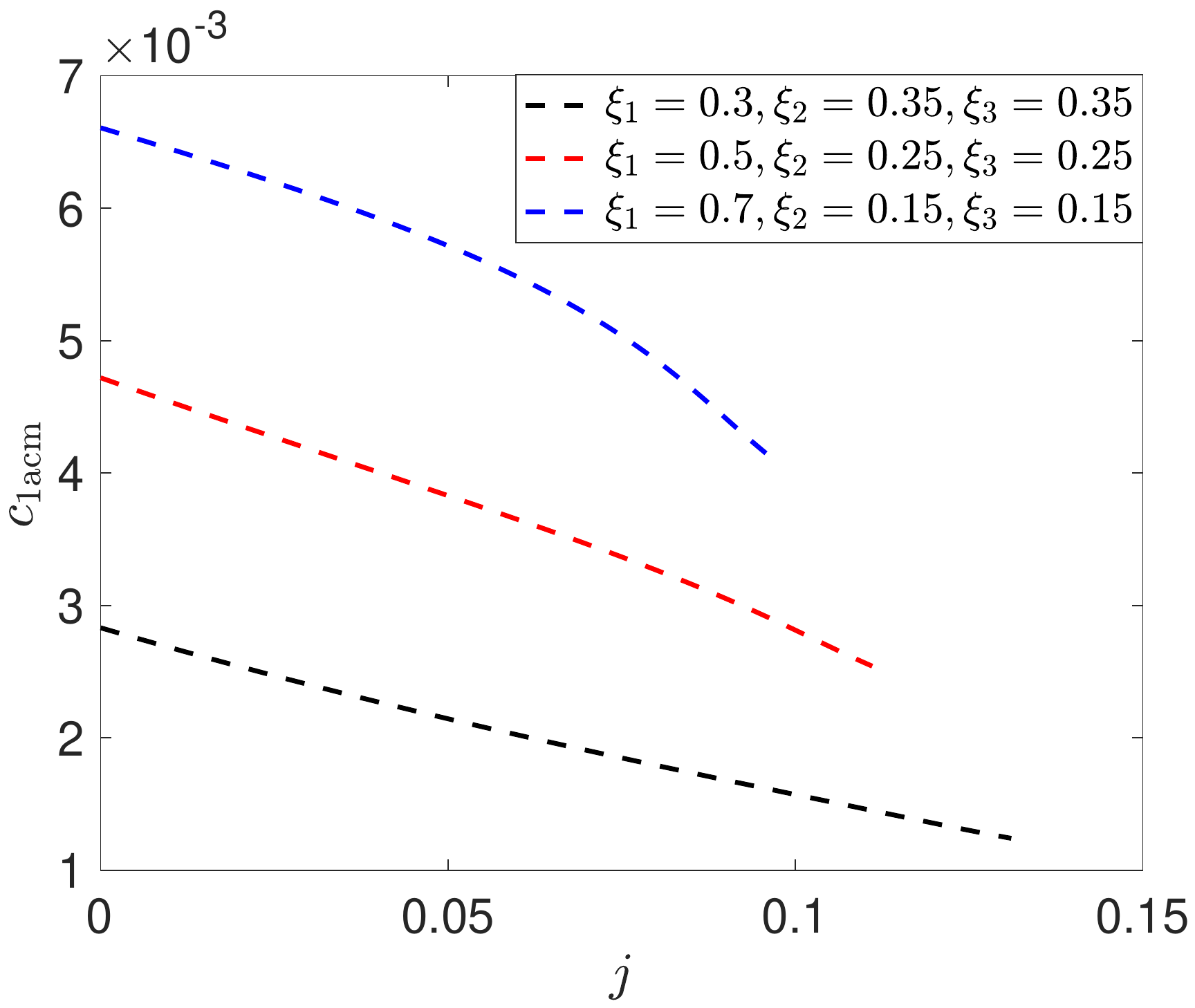}
{\scriptsize (c)}\includegraphics[scale=.36]{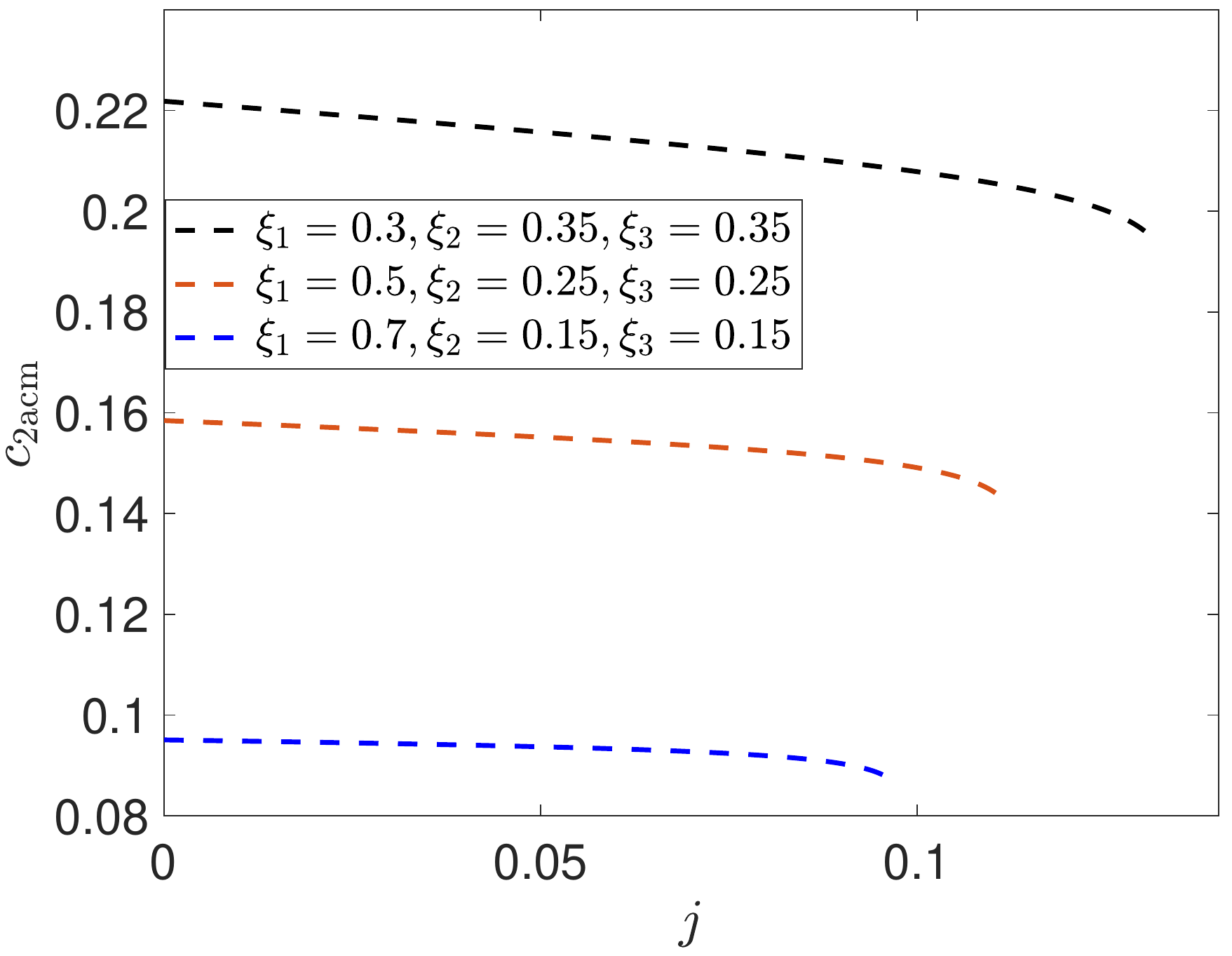}
{\scriptsize (d)}\includegraphics[scale=.36]{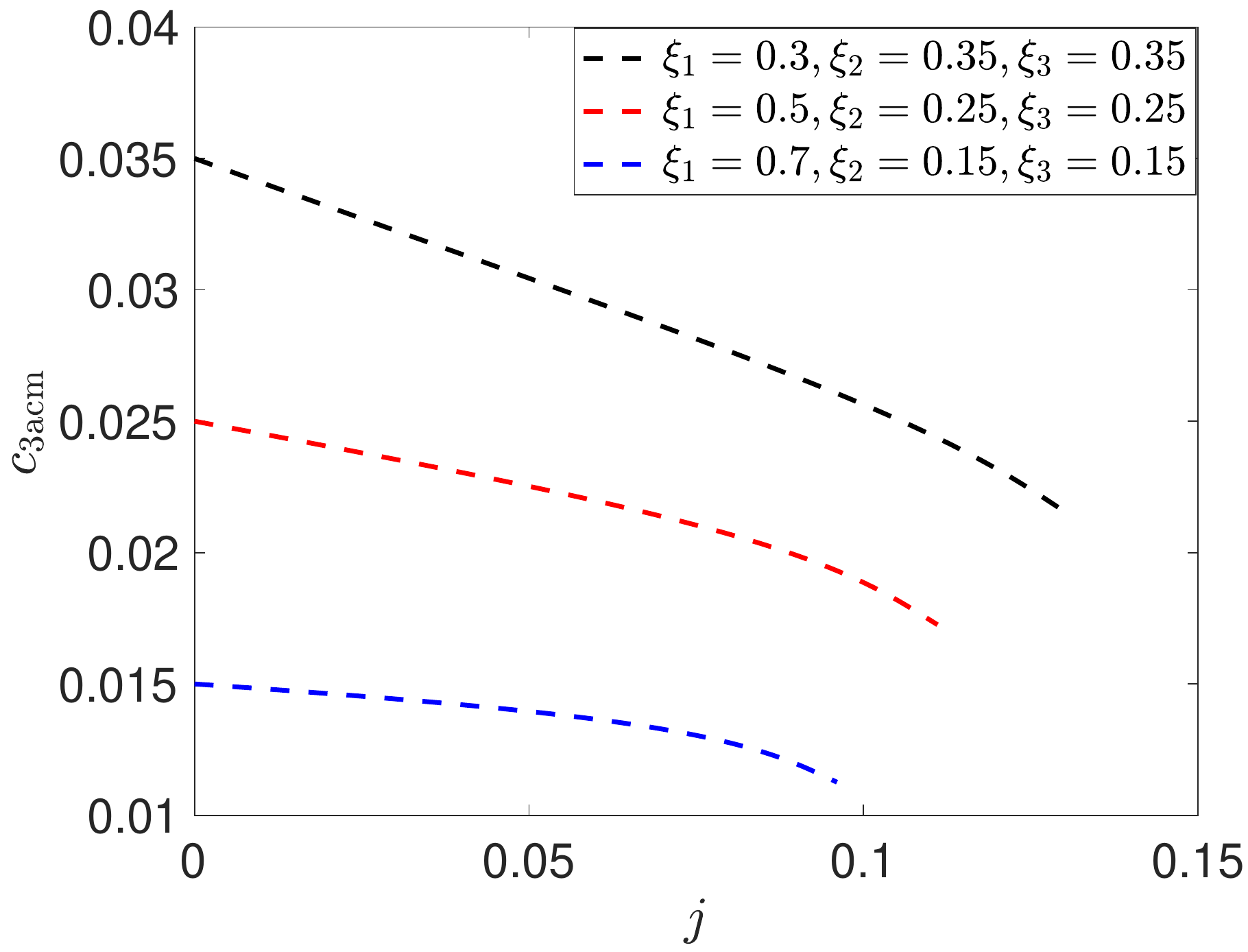}
\caption{
\footnotesize{ 
Evolution of optimized filters for three different feeds (1) $\xi_1=0.3, \xi_2=0.35, \xi_3=0.35$, (2) $\xi_1=0.5, \xi_2=0.25, \xi_3=0.25$ and (3) $\xi_1=0.7, \xi_2=0.15, \xi_3=0.15$, with $\beta_1=1$, $\beta_2=0.1, \beta_3=0.5$, $\lambda_1=1$, $R_1(0) \ge 0.99$, $R_2(0) \le 0.5$ and $R_3(0) \ge 0.9$: (a) ($u, j$) plot, (b) ($c_{1\rm acm}, j$) plot, (c) ($c_{2\rm acm}, j$) plot, (d) ($c_{3\rm acm}, j$) plot.
}}
\label{apfig1}
\end{figure}

\begin{table}
\centering
\begin{tabular}{|m{1cm}|m{1cm}|m{1cm}|c|c|c|c|c|}
\hline
{\bf $\xi_1$} & {\bf $\xi_2$}  & {\bf $\xi_3$} & {\bf $\bar R_1(\rm t_f)$} & {\bf $\bar R_2(\rm t_f)$} & {\bf $\bar R_3(\rm t_f)$} & {\bf $k_2$}  & {\bf $j(\rm t_f)$}\\
\hline
$ 0.3 $ & $0.35$ &0.35 &0.996 & 0.443 &0.939  &0.896 & 0.131\\
$ 0.5 $ & $0.25$ &0.25 &0.995 & 0.427 &0.931  &0.879 & 0.111\\
$ 0.7 $ & $0.15$ &0.15 &0.994 & 0.416 &0.925  &0.850 & 0.096 \\
\hline
\end{tabular}
\caption{\footnotesize{Three species feed filtration. We record $\xi_i$, particle ratios for each type of particles,  $\bar R_i(\rm t_f)$ the final cumulative particle removal ratio for particle type $i$, $k_2$ purity of type 2 particles in the final filtrate, $j(\rm t_f)$ total throughput.}}
\label{tb:3species}
\end{table}

\bibliography{filtration}   

\end{document}